\documentclass[letterpaper,11pt,oneside,reqno]{amsart}
\usepackage{bm}
\usepackage{mathrsfs}
\usepackage{amsfonts,amsmath, amssymb,amsthm,amscd,stmaryrd}
\usepackage[height=9.6in,width=5.95in]{geometry}
\usepackage[mpexclude,DIV13]{typearea}
\usepackage{verbatim}
\usepackage{hyperref}
\usepackage{graphicx}
\usepackage[latin1]{inputenc}
\usepackage{latexsym}
\usepackage{lscape}
\usepackage{epsfig}
\include{bibtex}

\usepackage{setspace}

\usepackage{parskip}
\begin{document}
\bibliographystyle{alpha}
\newcommand{\cn}[1]{\overline{#1}}
\newcommand{\e}[0]{\epsilon}
\newcommand{\bbf}[0]{\mathbf}

\newcommand{\Pfree}[5]{\ensuremath{\mathbb{P}^{#1,#2,#3,#4,#5}}}
\newcommand{\PfreeShort}{\ensuremath{\mathbb{P}^{BB}}}

%AAAAAAA
\newcommand{\WH}[8]{\ensuremath{\mathbb{W}^{#1,#2,#3,#4,#5,#6,#7}_{#8}}}
\newcommand{\Wfree}[5]{\ensuremath{\mathbb{W}^{#1,#2,#3,#4,#5}}}
\newcommand{\WHShort}[3]{\ensuremath{\mathbb{W}^{#1,#2}_{#3}}}
\newcommand{\WHShortCouple}[2]{\ensuremath{\mathbb{W}^{#1}_{#2}}}

\newcommand{\walk}[3]{\ensuremath{X^{#1,#2}_{#3}}}
\newcommand{\walkupdated}[3]{\ensuremath{\tilde{X}^{#1,#2}_{#3}}}
\newcommand{\walkfull}[2]{\ensuremath{X^{#1,#2}}}
\newcommand{\walkfullupdated}[2]{\ensuremath{\tilde{X}^{#1,#2}}}
%AAAAAAA

\newcommand{\PH}[8]{\ensuremath{\mathbb{Q}^{#1,#2,#3,#4,#5,#6,#7}_{#8}}}
\newcommand{\PHShort}[1]{\ensuremath{\mathbb{Q}_{#1}}}
\newcommand{\PHExp}[8]{\ensuremath{\mathbb{F}^{#1,#2,#3,#4,#5,#6,#7}_{#8}}}
%%%%

\newcommand{\D}[8]{\ensuremath{D^{#1,#2,#3,#4,#5,#6,#7}_{#8}}}
\newcommand{\DShort}[1]{\ensuremath{D_{#1}}}
\newcommand{\partfunc}[8]{\ensuremath{Z^{#1,#2,#3,#4,#5,#6,#7}_{#8}}}
\newcommand{\partfuncShort}[1]{\ensuremath{Z_{#1}}}
\newcommand{\bolt}[8]{\ensuremath{W^{#1,#2,#3,#4,#5,#6,#7}_{#8}}}
\newcommand{\boltShort}[1]{\ensuremath{W_{#1}}}
\newcommand{\boltNew}{\ensuremath{W}}
\newcommand{\QTLH}{\ensuremath{\mathfrak{H}}}
\newcommand{\QTLHgen}{\ensuremath{\mathfrak{L}}}

\newcommand{\whitenoise}{\ensuremath{\mathscr{\dot{W}}}}
\newcommand{\mf}{\mathfrak}

\newcommand{\EE}{\ensuremath{\mathbb{E}}}
\newcommand{\PP}{\ensuremath{\mathbb{P}}}
\newcommand{\var}{\textrm{var}}
\newcommand{\N}{\ensuremath{\mathbb{N}}}
\newcommand{\R}{\ensuremath{\mathbb{R}}}
\newcommand{\C}{\ensuremath{\mathbb{C}}}
\newcommand{\Z}{\ensuremath{\mathbb{Z}}}
\newcommand{\Q}{\ensuremath{\mathbb{Q}}}
\newcommand{\T}{\ensuremath{\mathbb{T}}}
\newcommand{\E}[0]{\mathbb{E}}
\newcommand{\OO}[0]{\Omega}
\newcommand{\F}[0]{\mathfrak{F}}
\def \Ai {{\rm Ai}}
\newcommand{\G}[0]{\mathfrak{G}}
\newcommand{\ta}[0]{\theta}
\newcommand{\w}[0]{\omega}
\newcommand{\ra}[0]{\rightarrow}
\newcommand{\vectoro}{\overline}
\newcommand{\crairy}{\mathcal{CA}}
\newcommand{\nc}{\mathsf{NoTouch}}
\newcommand{\ncf}{\mathsf{NoTouch}^f}
\newcommand{\wxy}{\mathcal{W}_{k;\bar{x},\bar{y}}}
\newcommand{\AP}{\mathfrak{a}}
\newcommand{\cm}{\mathfrak{c}}
\newtheorem{theorem}{Theorem}[section]
\newtheorem{partialtheorem}{Partial Theorem}[section]
\newtheorem{conj}[theorem]{Conjecture}
\newtheorem{lemma}[theorem]{Lemma}
\newtheorem{proposition}[theorem]{Proposition}
\newtheorem{corollary}[theorem]{Corollary}
\newtheorem{claim}[theorem]{Claim}
\newtheorem{experiment}[theorem]{Experimental Result}

\def\todo#1{\marginpar{\raggedright\footnotesize #1}}
\def\change#1{{\color{green}\todo{change}#1}}
\def\note#1{\textup{\textsf{\color{blue}(#1)}}}

\theoremstyle{definition}
\newtheorem{rem}[theorem]{Remark}

\theoremstyle{definition}
\newtheorem{com}[theorem]{Comment}

\theoremstyle{definition}
\newtheorem{definition}[theorem]{Definition}

\theoremstyle{definition}
\newtheorem{definitions}[theorem]{Definitions}

\theoremstyle{definition}
\newtheorem{conjecture}[theorem]{Conjecture}

\newcommand{\airysh}{\mathcal{A}}
\newcommand{\hfixed}{\mathcal{H}}
\newcommand{\afixed}{\mathcal{A}}
\newcommand{\canopynoarg}{\mathsf{C}}
\newcommand{\canopy}[3]{\ensuremath{\mathsf{C}_{#1,#2}^{#3}}}
\newcommand{\argmax}{x_{{\rm max}}}
\newcommand{\zmax}{z_{{\rm max}}}

\newcommand{\Rkle}{\ensuremath{\mathbb{R}^k_{>}}}
\newcommand{\Ronele}{\ensuremath{\mathbb{R}^k_{>}}}
\newcommand{\ewxy}{\mathcal{E}_{k;\bar{x},\bar{y}}}

\newcommand{\bxyf}{\mathcal{B}_{\bar{x},\bar{y},f}}
\newcommand{\bxyflr}{\mathcal{B}_{\bar{x},\bar{y},f}^{\ell,r}}

\newcommand{\bxyfone}{\mathcal{B}_{x_1,y_1,f}}

\newcommand{\ptac}{p}
\newcommand{\ptact}{v}

\newcommand{\fext}{\mathfrak{F}_{{\rm ext}}}
\newcommand{\gext}{\mathfrak{G}_{{\rm ext}}}
\newcommand{\xext}{{\rm xExt}(\mathfrak{c}_+)}

\newcommand{\dd}{\, {\rm d}}
\newcommand{\signc}{\Sigma}
\newcommand{\wxylr}{\mathcal{W}_{k;\bar{x},\bar{y}}^{\ell,r}}
\newcommand{\wxylrprime}{\mathcal{W}_{k;\bar{x}',\bar{y}'}^{\ell,r}}
\newcommand{\Rklezero}{\ensuremath{\mathbb{R}^k_{>0}}}
\newcommand{\XYfM}{\textrm{XY}^{f}_M}

\newcommand{\upright}{SC}
\newcommand{\staircase}{SC}
\newcommand{\energy}{E}
\newcommand{\xmax}{{\rm max}_1}
\newcommand{\ymax}{{\rm max}_2}
\newcommand{\lppls}{\mathcal{L}}
\newcommand{\lpplsre}{\mathcal{L}^{{\rm re}}}
\newcommand{\lpplsarg}[1]{\mathcal{L}_{n}^{\fa \to #1}}
\newcommand{\larg}[3]{\mathcal{L}_{n}^{#1,#2;#3}}
\newcommand{\BP}{M}
\newcommand{\weight}{\mathsf{Wgt}}
\newcommand{\properweight}{\mathsf{PropWgt}}
\newcommand{\pairweight}{\mathsf{PairWgt}}
\newcommand{\sumweight}{\mathsf{SumWgt}}
\newcommand{\mpgood}{\mathcal{G}}
\newcommand{\mpg}{\mathsf{Fav}}
\newcommand{\mcgone}{\mathsf{Fav}_1}
\newcommand{\radnik}[2]{\mathsf{RN}_{#1,#2}}
\newcommand{\size}[2]{\mathsf{S}_{#1,#2}}
\newcommand{\pdr}{\mathsf{PolyDevReg}}
\newcommand{\pwr}{\mathsf{PolyWgtReg}}
\newcommand{\lwr}{\mathsf{LocWgtReg}}
\newcommand{\hwp}{\mathsf{HighWgtPoly}}
\newcommand{\fbr}{\mathsf{ForBouqReg}}
\newcommand{\bbr}{\mathsf{BackBouqReg}}
\newcommand{\fsc}{\mathsf{FavSurCon}}
\newcommand{\maxpoly}{\mathrm{MaxDisjtPoly}}
\newcommand{\disjtindex}{\mathsf{DisjtIndex}}
\newcommand{\horsepindex}{\mathsf{HorSepIndex}}
\newcommand{\maxswf}{\mathsf{MaxScSumWgtFl}}
\newcommand{\emaxswf}{\e \! - \! \maxswf}
\newcommand{\minswf}{\mathsf{MinScSumWgtFl}}
\newcommand{\eminswf}{\e \! - \! \minswf}
\newcommand{\surreg}{\mathcal{R}}
\newcommand{\scf}{\mathsf{FavSurgCond}}
\newcommand{\disjtpoly}{\mathsf{DisjtPoly}}
\newcommand{\intint}[1]{\llbracket 1,#1 \rrbracket}
\newcommand{\maxsym}{*}
\newcommand{\polynum}{\#\mathsf{Poly}}
\newcommand{\dlp}{\mathsf{DisjtLinePoly}}
\newcommand{\lowb}{\underline{B}}
\newcommand{\highb}{\overline{B}}
\newcommand{\tottt}{t_{1,2}^{2/3}}
\newcommand{\tot}{t_{1,2}}
\newcommand{\btone}{t_1}
\newcommand{\bttwo}{t_2}
\newcommand{\formerE}{C}
\newcommand{\rcon}{r_0}
\newcommand{\para}{Q}
\newcommand{\Cstrong}{E}

\newcommand{\fluc}{\mathsf{Fluc}}

\newcommand{\mc}{\mathcal}
\newcommand{\vect}{\mathbf}
\newcommand{\bt}{\mathbf{t}}
\newcommand{\scB}{\mathscr{B}}
\newcommand{\scBres}{\mathscr{B}^{\mathrm{re}}}
\newcommand{\rightshadow}{\mathrm{RS}Z}
\newcommand{\dbm}{D}
\newcommand{\edgedbm}{D^{\rm edge}}
\newcommand{\gue}{\mathrm{GUE}}
\newcommand{\edgegue}{\mathrm{GUE}^{\mathrm{edge}}}
\newcommand{\eqdist}{\stackrel{(d)}{=}}
\newcommand{\geqdist}{\stackrel{(d)}{\succeq}}
\newcommand{\leqdist}{\stackrel{(d)}{\preceq}}
\newcommand{\scal}{{\rm sc}}
\newcommand{\fa}{x_0}
\newcommand{\hit}{H}
\newcommand{\scaledle}{\mathsf{Nr}\mc{L}}
\newcommand{\cenleup}{\mathscr{L}^{\uparrow}}
\newcommand{\cenledown}{\mathscr{L}^{\downarrow}}
\newcommand{\eln}{T}
\newcommand{\xmin}{{\rm Corner}^{\mfl,\mc{F}}}
\newcommand{\ymin}{{\rm Corner}^{\mfr,\mc{F}}}
\newcommand{\barxmin}{\overline{\rm Corner}^{\mfl,\mc{F}}}
\newcommand{\barymin}{\overline{\rm Corner}^{\mfr,\mc{F}}}
\newcommand{\qmin}{Q^{\mc{F}^1}}
%^{{\rm min}}}
\newcommand{\barqmin}{\bar{Q}^{\mc{F}^1}}
%^{{\rm min}}}
\newcommand{\test}{T}
\newcommand{\mfl}{\mf{l}}
\newcommand{\mfr}{\mf{r}}
\newcommand{\gfl}{\ell}
\newcommand{\gfr}{r}
\newcommand{\jre}{J}
\newcommand{\highfl}{{\rm HFL}}
\newcommand{\flyleap}{\mathsf{FlyLeap}}
\newcommand{\touch}{\mathsf{Touch}}
\newcommand{\notouch}{\mathsf{NoTouch}}
\newcommand{\close}{\mathsf{Close}}
\newcommand{\abovepar}{\mathsf{High}}
\newcommand{\vecint}{\bar{\iota}}
\newcommand{\cornthree}{{\rm Corner}^\mc{G}_{k,\mfl}}
\newcommand{\cornfour}{{\rm Corner}^\mc{H}_{k,\fa}}
\newcommand{\mpgg}{\mathsf{Fav}_{\mc{G}}}

\newcommand{\lefta}{M_{1,k+1}^{[-2\eln,\gfl]}}
\newcommand{\mida}{M_{1,k+1}^{[\gfl,\gfr]}}
\newcommand{\righta}{M_{1,k+1}^{[\gfr,2\eln]}}

\newcommand{\ipdval}{d}
\newcommand{\ctemp}{d_0}

\newcommand{\wien}{W}
\newcommand{\pole}{P}
\newcommand{\pp}{p}

\newcommand{\rmreg}{{\rm Reg}}

\newcommand{\const}{D_k}
\newcommand{\numcone}{14}
\newcommand{\numctwo}{13}
\newcommand{\numcthree}{6}
\newcommand{\rsC}{C}
\newcommand{\rsc}{c}
\newcommand{\cone}{c_1}
\newcommand{\Cone}{C_1}
\newcommand{\Ctwo}{C_2}
\newcommand{\smallc}{c_0}
\newcommand{\smallcprime}{c_1}
\newcommand{\smallcanother}{c_2}
\newcommand{\smallcnew}{c_3}
\newcommand{\Cda}{D}
\newcommand{\Kzero}{K_0}
\newcommand{\Rmac}{\varphi}
\newcommand{\rmac}{r}
\newcommand{\conseqmac}{H}
\newcommand{\constn}{C'}
\newcommand{\coninit}{\Psi}
\newcommand{\condee}{\hat{D}}
\newcommand{\conbrac}{\hat{C}}
\newcommand{\Cnew}{\tilde{C}}
\newcommand{\Cbig}{C^*}
\newcommand{\Ctbd}{C_+}
\newcommand{\Ctbs}{C_-}

\newcommand{\Cpop}{G}
\newcommand{\cpop}{g}
\newcommand{\correc}{Q}

\newcommand{\high}{{\rm High}}
\newcommand{\notlow}{{\rm NotLow}}

\newcommand{\maxmin}{\pwr}

\newcommand{\imax}{i_{{\rm max}}}

\newcommand{\wlp}{{\rm WLP}}

\newcommand{\canopynumber}{\mathsf{Canopy}{\#}}

\newcommand{\cannum}{{\#}\mathsf{SC}}

\newcommand{\boundgood}{\mathsf{G}}
\newcommand{\lshift}{\mc{L}^{\rm shift}}
\newcommand{\deltapi}{\theta}
\newcommand{\rootneigh}{\mathrm{RNI}}
\newcommand{\rootneighuse}{\mathrm{RNI}}
\newcommand{\manycan}{\mathsf{ManyCanopy}}
\newcommand{\specialpt}{\mathrm{spec}}

\newcommand{\dist}{\vert\vert}
\newcommand{\fik}{\mc{F}_i^{[K,K+1]^c}}
\newcommand{\mcfa}{\mc{H}[\fa]}
\newcommand{\tent}{{\rm Tent}}
\newcommand{\goodk}{\mc{G}_{K,K+1}}
\newcommand{\pairsep}{{\rm PS}}
\newcommand{\mbf}{\mathsf{MBF}}
\newcommand{\nbd}{\mathsf{NoBigDrop}}
\newcommand{\bd}{\mathsf{BigDrop}}
\newcommand{\jleft}{j_{{\rm left}}}
\newcommand{\jright}{j_{{\rm right}}}
\newcommand{\smalljfluc}{\mathsf{SmallJFluc}}
\newcommand{\mfone}{M_{\mc{F}^1}}
\newcommand{\mfthree}{M_{\mc{G}}}
\newcommand{\rhomac}{P}
\newcommand{\phimac}{\varphi}
\newcommand{\phimacone}{1/3}
\newcommand{\phimactwo}{1/9}
\newcommand{\phimacthree}{1/3}
\newcommand{\chimac}{\chi}
\newcommand{\xnmac}{z_{\mathcal{L}}}
\newcommand{\nmac}{N}
\newcommand{\Cwb}{E_0}
\newcommand{\initcond}{\mathcal{I}}
\newcommand{\neargeod}{\mathsf{NearPoly}}
\newcommand{\nearpoly}{\mathsf{NearPoly}}
\newcommand{\polyunique}{\mathrm{PolyUnique}}
\newcommand{\latecoal}{\mathsf{LateCoal}}
\newcommand{\nolatecoal}{\mathsf{NoLateCoal}}
\newcommand{\normalcoal}{\mathsf{NormalCoal}}
\newcommand{\regfluc}{\mathsf{RegFluc}}
\newcommand{\mdeltaweight}{\mathsf{Max}\Delta\mathsf{Wgt}}
\newcommand{\ovbar}[1]{\mkern 1.5mu\overline{\mkern-1.5mu#1\mkern-1.5mu}\mkern 1.5mu}

\newcommand{\paradelta}{\Delta^{\cup}\,}

\newcommand{\kay}{k}
\newcommand{\emm}{m}

\newcommand{\correctmac}{R}

%\ifx\hyperlink\undefined
\def\fff#1{&{{\pageref{#1}}}\cr}
\def\hfff#1{\label{#1}}
%\else
%\def\fff#1{&{\hyperlink{#1}{\pageref*{#1}}}\cr}
%\def\hfff#1{\label{#1}\hypertarget{#1}}
%\fi

\title[Rarity of disjoint polymers in last passage percolation]{Exponents governing  the rarity of disjoint polymers \\ in Brownian last passage percolation}

\author[A. Hammond]{Alan Hammond}
\address{A. Hammond\\
  Department of Mathematics and Statistics\\
 U.C. Berkeley \\
 899 Evans Hall \\
  Berkeley, CA, 94720-3840 \\
  U.S.A.}
  \email{alanmh@berkeley.edu}
  \thanks{The author is supported by NSF grant DMS-$1512908$.}

\begin{abstract} 
In last passage percolation models lying in the KPZ universality class, long maximizing paths have a typical deviation from the linear interpolation of their endpoints governed by the two-thirds power of the interpolating distance. This two-thirds power dictates a choice of scaled coordinates, in which these maximizers, now called polymers,  cross unit distances with unit-order fluctuations.
%The resulting field of polymers, indexed by a pair of planar endpoints, . 
In this article, we consider Brownian last passage percolation in these scaled coordinates, and prove that the probability of the presence of $k$ disjoint polymers crossing a unit-order region
while beginning and ending within a short distance $\e$ of each other is bounded above by~$\e^{(k^2 - 1)/2 \, + \,  o(1)}$.
% has a superpolynomial tail in $k$, uniformly in the scaling parameter. 
 This result, which we conjecture to be sharp, yields understanding of the uniform nature of the coalescence structure of polymers, and plays a foundational role in~\cite{Patch} in proving comparison on unit-order scales  to Brownian motion for polymer weight profiles from  general initial data. The present paper also contains an on-scale articulation of the two-thirds power law for polymer geometry: polymers fluctuate by $\e^{2/3}$ on short scales $\e$.
\end{abstract}

\maketitle

\tableofcontents

\section{Introduction}

\subsection{KPZ universality, last passage percolation models, and scaled coordinates}
The $1 + 1$ dimensional Kardar-Parisi-Zhang (KPZ) universality class includes a wide range of random interface models suspended over a one-dimensional domain, in which growth in a direction normal to the surface competes with a smoothening surface tension in the presence of a local randomizing force that roughens the surface.
These characteristic features are evinced by many last passage percolation models. Such an LPP model comes equipped with a planar random environment, which is independent in disjoint regions. 
Directed paths, that are permitted say to move only in a direction in the first quadrant, are then assigned energy via this randomness, by say integrating the environment's value along the path.
For a given pair of planar points, the path attaining the maximum weight over directed paths with such endpoints is called a geodesic.

For LPP models lying in the KPZ class, a geodesic that crosses a large distance $n$ (in say a northeasterly direction) has an energy that grows linearly in $n$ with a standard deviation of order~$n^{1/3}$. If the lower geodesic endpoint is held fixed, and the higher one is permitted to vary horizontally, then the geodesic energy as a function of the variable endpoint plays the role of the random interface mentioned at the outset. (With the first endpoint thus fixed, the energy profile goes by the name `narrow wedge'.)
 Non-trivial correlations in the geodesic energy are witnessed when this horizontal variation has order $n^{2/3}$. 
 These assertions have been rigorously demonstrated for only a few LPP models, each of which enjoys an integrable structure: the seminal work of Baik, Deift and Johansson~~\cite{BDJ1999}
 rigorously established the one-third exponent, and the GUE Tracy-Widom distributional limit, for the case of Poissonian last passage percolation, while the two-thirds power law for maximal transversal fluctuation was derived for this model by Johansson~~\cite{Johansson2000}.

 In seeking to understand the canonical structures of KPZ universality, we are led, in view of these facts, to represent the field of geodesics in a scaled system of coordinates, under which a northeasterly displacement of order~$n$  becomes a vertical displacement of one unit, and a horizontal displacement of order $n^{2/3}$ becomes a unit horizontal displacement. Moreover, the accompanying system of energies also transfers to scaled coordinates, with the scaled geodesic energy being specified by centring about the mean value and division by the typical scale of $n^{1/3}$.

 \subsection{Probabilistic proof techniques for KPZ, and scaled Brownian LPP}

The theory of KPZ universality has advanced through physical insights, numerical analysis, and several techniques of integrable or algebraic origin. We will not hazard a summary of literature to support this one-sentence history, but refer to the reader to~\cite{IvanSurvey} for a KPZ survey from 2012; in fact, integrable and analytic approaches to KPZ have attracted great interest around and since that time. 
Now, it is hardly  deep or controversial to say that  many problems and models in KPZ are intrinsically random. This fact may suggest that it would be valuable to approach the problems of KPZ universality from a predominately probabilistic perspective.

An important illustration is offered by 
 Brownian last passage percolation. This LPP model has very attractive probabilistic features: for example, in the narrow wedge case, the scaled geodesic energy profile may be embedded as the uppermost curve in a Dyson diffusion, namely a system of one-dimensional Brownian motions conditioned on mutual avoidance (with a suitable boundary condition). 
If we depict Brownian LPP in the scaled coordinates that have been described, may we analyse it with probabilistic tools and thus gain some insight into universal KPZ structures? 

The present article forms part of a four-paper answer to this question. 
The companion papers are~\cite{BrownianReg},~\cite{ModCon} and~\cite{Patch}.
Our focus is on scaled coordinates, and we will adopt the terms {\em polymer} and {\em weight} to refer to scaled geodesics and their scaled energy. 
The reader may glance ahead to the right sketch in Figure~\ref{f.scaling} for a depiction of a polymer, bearing in mind that, when the scaling parameter $n$ is high, the microscopic backtracking apparent in the illustration becomes negligible, so that 
the lifetime~$[t_1,t_2]$ of a polymer may be viewed as a vertical interval, with the polymer being a random real-valued function defined on that interval; the range of the polymer thus crosses the planar strip $\R \times [t_1,t_2]$ from the lower to the upper side.

\subsection{Principal conclusions and themes in overview}
% Brownian last passage percolation is an LPP model with attractive integrable and probabilistic features whose definition we will recall shortly.
 %Notable among these is that the locally Brownian nature of scaled energy profiles, which for other LPP models should emerge only in the high~$n$ limit, may be expected to be present for finite~$n$. 
% It is a central goal of the present article to 
% rigorously implement ideas in the preceding paragraph, for Brownian LPP. 
In this article, we will reach three principal conclusions concerning scaled Brownian LPP. As we informally summarise them now, we omit mention of the scaling parameter $n \in \N$: roughly, our assertions should be understood uniformly in high choices of this parameter.
 \begin{itemize}
 \item Consider the event that there exist $k$ polymers of unit lifetime~$[0,1]$,
  where each begins and ends in a given interval of  some given {\em small} scaled length~$\e$. 
 %This theorem's proof depends on our second result, 
Theorem~\ref{t.disjtpoly.pop},
 states an upper bound of $\e^{(k^2 - 1)/2 \, + \, o(1)}$
 on the probability of this event. The left sketch of the upcoming Figure~\ref{f.triple} depicts the event with $k=3$.
 \item Derived as a consequence, Theorem~\ref{t.maxpoly.pop} asserts that the probability that $\emm$ polymers coexist disjointly in a unit-order region has a superpolynomial decay in $\emm$. 
 \item We have mentioned that
 the maximal geodesic fluctuation in Poissonian last passage percolation has been shown in~\cite{Johansson2000}
 to be governed by an exponent of two-thirds. Here is a further expression of this two-thirds power law, showing that this exponent also governs local behaviour of the scaled geodesic:  a polymer fluctuates by more than $\e^{2/3} r$ on a short duration~$\e$ with probability at most~$\exp \big\{ - O(r^{\alpha}) \big\}$ (for a broad range of values of~$r$).  In Theorem~\ref{t.polyfluc}, we will prove a Brownian LPP version of this assertion, with $\alpha = 3/4$.  
 %The result, which is needed to prove Theorem~\ref{t.disjtpoly.pop},
 %is a useful tool for understanding the geometry of scaled geodesics and the forests
 %formed by their coalescence. 
 \end{itemize}
 
It is the first of these results, namely the exponent bound of the form $(k^2 - 1)/2 \, + \, o(1)$ on the rarity of~$k$ disjoint polymers with nearby endpoints, which we regard as the most fundamental.  
Four reasons, discussed in turn next, will help to explain why.

\begin{figure}[ht]
\begin{center}
\includegraphics[height=12cm]{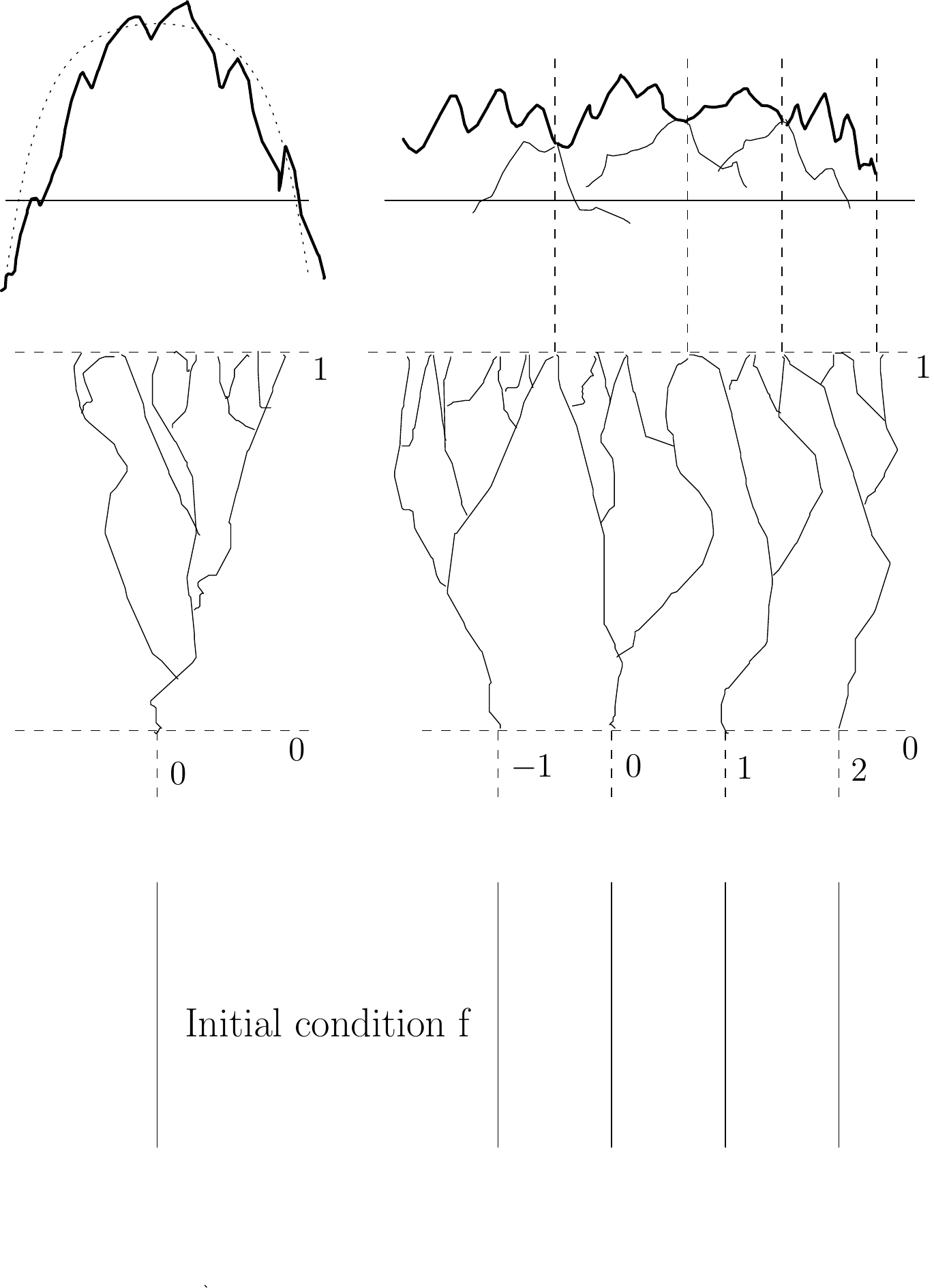}
\caption{
The left and right sketches each depict three aspects of scaled random growth in Brownian LPP, corresponding to two different initial conditions. 
These initial conditions $f$ are depicted at the bottom. On the left, the narrow wedge case of growth from zero is depicted. The polymer weight profile, top left, is known to be locally Brownian (and globally parabolic). In the middle left, we see the associated system of polymers, forming a tree with root at the origin. On the right, growth is instead initiated from the set of integers. The field of polymers now forms a forest, rather than a tree, and the polymer weight profile, which is depicted in bold, is the maximum of the profiles associated to the roots of the various trees. Indeed, a general polymer weight
follows the profile associated to the consecutive trees in the associated polymer forest. This profile may be understood to be locally Brownian as an inference from the Brownian resemblance enjoyed by the `narrow-wedge' profiles associated to the individual trees. To make this idea work, we need to understand that there are not too many trees: their number per unit length must be shown to be tight in the scaling parameter $n$. It is here that Theorem~\ref{t.disjtpoly.pop}, and Theorem~\ref{t.maxpoly.pop}, provide the necessary input.}
\label{f.manytrees}
\end{center}
\end{figure}

\noindent{\em 1. Validating the Brownian regularity of polymer weight profiles begun with general initial data.}
Theorem~\ref{t.disjtpoly.pop} is very useful. In fact, Theorem~\ref{t.maxpoly.pop} will be derived as a rather direct consequence. More fundamentally, Theorem~\ref{t.disjtpoly.pop} plays a key role in the principal highway of implications in our four-paper study. In order to indicate why briefly, we mention first a key theme concerning scaled Brownian LPP: that {\em narrow wedge polymer weight profiles resemble Brownian motion}. This assertion may be understood in a {\em local limit},
when the interval on which comparison is made between the weight profile and Brownian motion  shrinks to zero, in which case~\cite{Hagg} and~\cite{CatorPimentel,Pimentel18} offer rigorous expressions of this statement; they do so respectively by analysing the largest particle in the extended Airy point process and by means of a technique of local comparision to an equilibrium regime. But the assertion may also be understood on a {\em unit scale}: indeed, Theorem~$2.11$ and Proposition~$2.5$ of~\cite{BrownianReg} demonstrate that, after an affine adjustment, narrow wedge profiles withstand a very demanding comparison to Brownian bridge when the profile is restricted to any given compact interval.
In the narrow wedge case that these works address, random growth is initiated from a point. Such growth may also begin from much more general initial data.
The main conclusion of~\cite{Patch}, that article's Theorem~$1.2$, asserts that the corresponding polymer weight profiles
enjoy in a uniform sense a strong resemblance to Brownian motion on unit scales.
Theorem~\ref{t.disjtpoly.pop} is an engine that drives the derivation of this conclusion.
 Figure~\ref{f.manytrees} gives a very impressionistic account of why the theorem is valuable for this purpose.

\noindent{\em 2. Elucidating fractal geometry in  the  Airy sheet.}  For suitable LPP models, the polymer weight profile in the upper-left sketch of Figure~\ref{f.manytrees} converges distributionally~\cite{PrahoferSpohn,Johansson2003} in the high scaling parameter limit $n \to \infty$
to the Airy$_2$ process after   the subtraction of a parabola. This process, which has finite-dimensional distributions specified by Fredholm determinants, is a central object in KPZ universality. It has been expected that models in the KPZ universality class share richer universal structure than the Airy$_2$ process, and two significant recent advances have rigorously validated this expectation.

The first of these advances is 
the recent construction~\cite{MQR17} of the {\em KPZ fixed point}. The narrow wedge polymer weight profile may be viewed as a time-one snapshot of a scaled random growth process initiated from the origin at time zero. 
Growth may be initiated from a much more general initial condition. In~\cite{MQR17}, Matetski, Quastel and Remenik
have utilized a biorthogonal ensemble representation found by~\cite{Sas05,BFPS07} associated to the totally asymmetric exclusion process in order to  find Fredholm determinant formulas for the multi-point distribution of the height function of this growth process begun from an arbitrary initial condition. Using these formulas to take the KPZ scaling limit, the authors construct a scale invariant Markov process that lies at the heart of the KPZ universality class. The time-one evolution of this Markov process may be applied to very general initial data, and the result is the scaled profile begun from such data, which generalizes the~${\rm Airy}_2$ process seen in the narrow wedge case. 
Were~\cite{MQR17} to be adapted to Brownian last passage percolation, the present article's theorems and consequences might become applicable to universal KPZ structure. In particular, this outcome should be rather direct in the case of the consequence~\cite[Theorem~$1.2$]{Patch} concerning the Brownian nature of general initial condition profiles.

The second recent advance~\cite{DOV18} is perhaps of even greater relevance to the results that we present.
The {\em space-time Airy sheet} is a rich shared scaled structure in the KPZ universality class.
It simultaneously encodes the collection of weights of polymers running between any pair of planar endpoint locations. Its existence was mooted in~\cite{CQR2015}. In~\cite{DOV18}, the space-time Airy sheet is constructed  
 by use of an extension of the Robinson-Schensted-Knuth correspondence which expresses the construction in terms of a last passage percolation problem whose underlying environment is itself a copy of the distributional limit of the narrow wedge profile in Brownian last passage percolation. 

These advances~\cite{MQR17} and~\cite{DOV18} offer new prospects for probabilistic inquiry into KPZ, including applications of the present article. For example, fractal geometry embedded in the Airy sheet has been thus elucidated in~\cite{BGH18}. To explain this application, note that  
many copies of the Airy$_2$ process are coupled together in the Airy sheet. It is already interesting to consider how two such processes rooted at distinct planar points may be coupled: the random process sending $y \in \R$ to the {\em difference} in weight between polymers beginning at $(-1,0)$ and at $(1,0)$ and ending together at $(y,1)$ is such an example. 
This {\em Airy difference process}, mapping $\R$ to $\R$, is a random fractal, being the distribution function of a random Cantor set of Hausdorff dimension one-half. In proving this, \cite{BGH18} provides the upper bound on Hausdorff dimension as a consequence of Theorem~\ref{t.disjtpoly.pop} applied in the case of a pair of polymers.

% $\R \to \R: z \to W_n (0,z)$. The analysis of~\cite{DOV18} is assisted by~\cite{DV18}, an article making a Brownian Gibbs analysis of scaled Brownian  LPP in order to provide valuable estimates for the study of the very novel LPP problem introduced in~\cite{DOV18}.   

\noindent{\em 3. Polymer surgery as a KPZ tool.} 
 The proof of Theorem~\ref{t.disjtpoly.pop} has conceptual interest.
Two principal probabilistic techniques are at work in the four-paper study, namely {\em Brownian Gibbs resamplings} and {\em polymer surgery}. The latter is introduced in order to prove Theorem~\ref{t.disjtpoly.pop}. In order to explain how, it is useful first to discuss the former. 
This first, resampling, tool is employed to analyse the narrow wedge polymer weight profile by viewing it as the top curve in a Dyson diffusion. Allied with a simple observation from integrable probability, the Karlin-McGregor formula, the technique has been used in~\cite{BrownianReg} to solve a cousin of the problem addressed in Theorem~\ref{t.disjtpoly.pop}. The solution will later be presented as Theorem~\ref{t.neargeod}.   
Crudely, in place of considering, as  Theorem~\ref{t.disjtpoly.pop} does, $k$ disjoint polymers with nearby endpoints, we may instead consider the event that there exists a system of $k$ scaled paths, crossing between two {\em given} points at a unit vertical displacement, the paths pairwise disjoint except for these shared endpoints, and with each path coming close to being a polymer, in the sense that each path weight is close to maximal for the given endpoints. Theorem~\ref{t.neargeod} identifies in a sharp way the decay rate in the probability of this event as a natural parameter that measures this closeness tends to zero: the theorem identifies the exponent~$k^2 -1$ as governing this rate.

How then will polymer surgery help to prove Theorem~\ref{t.disjtpoly.pop}?
Faced with $k$ polymers with nearby endpoints, we seek to deform them surgically near their endpoints in order to typically achieve a shared-endpoint outcome whose probability is determined by Theorem~\ref{t.neargeod}: Figure~\ref{f.triple} is a one-sketch summary. A probabilistic bridge is being built into an understanding of integrable origin, and the proof technique is thus emblematic of a broader aim of expanding KPZ horizons by probabilistic means.

\noindent{\em 4. The exponent $(k^2 - 1)/2$.} In its key application in~\cite{Patch}, Theorem~\ref{t.disjtpoly.pop} is applied when $k=3$, for a triple of polymers. The form of the exponent is important in this application, dictating that a certain Radon-Nikodym derivative that describes the strength of Brownian comparison lies in~$L^{3-}$. Merely knowing that the exponent is at least a positive constant, for example, would not be adequate even for deriving $L^1$-membership.

We believe that the $(k^2 - 1)/2$ exponent in Theorem~\ref{t.disjtpoly.pop} is sharp and will present a conjecture to this effect in Subsection~\ref{s.maxpoly}. It is rigorously known that, in Theorem~\ref{t.disjtpoly.pop}'s proof, nothing is lost in the exponent after the reduction to the shared endpoint Theorem~\ref{t.neargeod} is made, because the latter result provides a sharp resolution of the cousin problem discussed above in the third point. Moreover, the conjecture of Theorem~\ref{t.disjtpoly.pop}'s sharpness has recently been validated in~\cite[Theorem~$2.4$]{BGH18} for the special case $k=2$ as a corollary of the identification of the Hausdorff dimension of the points of increase of the Airy difference process that we have mentioned.

%More generally, the form of the exponent reflects the unison of probabilistic and integrable techniques in this four-paper study; 

\subsection{A suggestion for further reading of overview}
In seeking to explain the progress made in this article in a few paragraphs, we have
also offered a whistle-stop tour of ideas in a broader study of scaled Brownian LPP. 
Such a brief presentation may well prompt questions, and we refer the reader who wants a more detailed but still informal overview of the main results and ideas in the four-paper study to~\cite[Section~$1.2$]{BrownianReg}. 

For the reader who wishes to read the broader study, the present article would be read after~\cite{BrownianReg} and~\cite{ModCon} but before~\cite{Patch}. However, the present article has been written so that it may be read on its own. As such, consulting the other articles, including  the further reading just suggested, is an optional extra. 

%two results from ModCon are used here

%Specifically, the role of Theorem~\ref{t.disjtpoly.pop}  

% An overview of the four-paper study appears in~\cite[Section~$1.2$]{BrownianReg}:  
% the ideas of forests of coalescing scaled geodesics, their relation to mutual coexistence of disjoint geodesics and their application to scaled energy profiles, are explained in a more leisurely way there.

%Several elements will be specified precisely later. 

% For example, in Section~\ref{s.multiweight}, elements of a theory of Brownian Gibbs line ensembles, which has been developed in~\cite{BrownianReg}, will be introduced. We will also employ two results from~\cite{ModCon} concerning scaled Brownian LPP.
% All three of this paper's main conclusions, mentioned above, will be then applied  in~\cite{Patch} in order to prove that scaled geodesic energy profiles are Brownian on unit-order scales, for very general initial data.    

\subsection{Connections: probabilistic tools, geodesics and coalescence}

The idea of taking limited integrable inputs, such as control of narrow wedge polymer profiles at given points, and exploiting them via probabilistic means to reach much stronger conclusions is central to this article and the broader four-paper Brownian LPP study. A similar conceptual approach governs many ideas  in~\cite{SlowBondSol} and~\cite{SlowBond}, where the slow bond conjecture for the totally asymmetric exclusion process, concerning the macroscopic effects of slightly attenuating passage of particles through the origin, is proved, and the associated invariant measures determined. Such a probabilistic technique is also encountered in~\cite{BSS17}. Indeed, in \cite[Theorem~$2$]{BSS17}, 
 an assertion similar to, and in fact slightly stronger than, Theorem~\ref{t.polyfluc} has been proved for exponential last passage percolation: in essence, the assertion in the third bullet point in the preceding discussion has been verified for~$\alpha = 1$.
 
Roughly opposite to the phenomenon of geodesic disjointness is the process of geodesic coalescence. An interesting arena of study of this coalescence is the class of first passage percolation models, which are variants of last passage percolation that are not integrable. Although scaling exponents are not rigorously known in this setting, it has been heuristically appreciated for some time that certain circumstances should coincide in this context: the non-existence of bi-infinite geodesics; an exponent of one-third in the standard deviation  of geodesic energy; and smoothness of the limiting shape of energy-per-unit-length as a function of angle. A summary of progress until 1995 regarding coalescence of geodesics is offered by Newman in~\cite{Newman95}; this progress included geodesic uniqueness for given direction, subject to a curvature assumption on the limiting shape. 
The theory of semi-infinite geodesics was developed by~\cite{FerrariPimentel05} and~\cite{Coupier11}. Building on this theory, Pimentel adapted Newman's technique for uniqueness to the context of exponential LPP, in which the curvature assumption is known to be verified. In this work, it is also proved that, for two geodesics of given direction that begin at points separated by distance $\ell$, it is typical that coalescence of the geodesic occurs by a time of  order $\ell^{3/2}$. Moreover, Pimentel conjectured that 
the probability that coalescence has failed to occur by time $r \ell^{3/2}$ has tail $r^{-2/3}$. In the context of finite geodesics,~\cite{BSS17} has proved an upper bound on this failure probability of the form $r^{-c}$ for some positive $c > 0$, and there seems to be promise in these techniques to obtain a sharper form of the result. Hoffman~\cite{Hoffman08} initiated very fruitful progress on infinite geodesics in first passage percolation
by using Busemann functions, with significant geometric information emerging in~\cite{DamronHanson14,AhlbergHoffman16} and~\cite{DamronHanson17}.
Busemann functions have also been studied in non-integrable last passage percolation contexts: see~\cite{GRS15a} and~\cite{GRS15b}. 
%[GRASepp Geodesics and the competition interface, GRASepp Stationary cocycles].

\subsection{Brownian last passage percolation}\label{s.brlpp}

In the two remaining introductory sections, we specify this model precisely, and set up notation for the use of scaled coordinates; and we state our main results.

\subsubsection{The model's definition.}
This model was introduced by~\cite{GlynnWhitt} and further studied in~\cite{O'ConnellYor};
we will call it Brownian LPP.
On a probability space carrying a law labelled~$\PP$, we let $B:\Z \times \R \to \R$ denote an ensemble of independent  two-sided standard Brownian motions $B(k,\cdot):\R\to \R$, $k \in \Z$.

Let $i,j \in \Z$ with $i \leq j$. 
We denote the integer interval $\{ i,\cdots  , j \}$ by $\llbracket i,j \rrbracket$.
Further let $x,y \in \R$ with $x \leq y$.
With these parameters given, we consider the collection of  non-decreasing lists 
 $\big\{ z_k: k \in \llbracket i+1,j \rrbracket \big\}$ of values $z_k \in [x,y]$. 
With the convention that $z_i = x$ and $z_{j+1} = y$,
we associate an energy to any such list, namely   $\sum_{k=i}^j \big( B ( k,z_{k+1} ) - B( k,z_k ) \big)$.
We may then define  the maximum energy, 
$M^1_{(x,i) \to (y,j)}$,  to be the supremum of the energies of all such lists. 
%(The superscript `$1$' anticipates a generalization that we will specify presently under which the maximum over several staircases is instead considered.)

\subsubsection{A geometric view: staircases.}\label{s.staircases}
In order to make a study of those lists that attain this maximum energy, 
we begin by noting that the lists are in bijection with certain subsets of $[x,y] \times [i,j] \subset \R^2$ that we call {\em staircases}.
Staircases offer a geometric perspective on Brownian LPP and perhaps help in visualizing the problems in question.

%To introduce them ,we first mention  some basic notation: for  $i,j \in \Z$ with $i \leq j$,
%we denote the interval interval $\{ i,\cdots  , j \}$ by $\llbracket i,j \rrbracket$.

 The staircase~\hfff{staircase} associated to the non-decreasing list $\big\{ z_k: k \in \llbracket i+1,j \rrbracket \big\}$ is specified as the union of certain horizontal planar line segments, and certain vertical ones.
The horizontal segments take the form $[ z_k,z_{k+1} ] \times \{ k \}$ for $k \in \llbracket i , j \rrbracket$.
Here, the convention that $z_i = x$ and  $z_{j+1} = y$ is again adopted. 
The right and left endpoints of each consecutive pair of horizontal segments are interpolated by a vertical planar line segment of unit length. It is this collection of vertical line segments that form
the vertical segments of the staircase.

The resulting staircase may be depicted as the range of an alternately rightward and upward moving path from starting point $(x,i)$ to ending point $(y,j)$. 
The set of staircases with these starting and ending points will be denoted by $\staircase_{(x,i) \to (y,j)}$.
Such staircases are in bijection with the collection of non-decreasing lists considered above. Thus, any staircase $\phi \in \staircase_{(x,i) \to (y,j)}$
is assigned an energy~\hfff{energy} $E(\phi) = \sum_{k=i}^j \big( B ( k,z_{k+1} ) - B( k,z_k ) \big)$ via the associated $z$-list.

\subsubsection{Energy maximizing staircases are called geodesics.}
A staircase  $\phi \in \staircase_{(x,i) \to (y,j)}$ whose energy  attains the maximum value $M^1_{(x,i) \to (y,j)}$ is called a geodesic~\hfff{geodesic} from $(x,i)$ to~$(y,j)$.
It is a simple consequence of the continuity of the constituent Brownian paths $B(k,\cdot)$
that such a geodesic exists for all choices of $(x,y) \in \R^2$ with $x \leq y$.
As we will later explain, in Lemma~\ref{l.severalpolyunique}, the geodesic with given endpoints is almost surely unique.

\subsubsection{The scaling map.}
For $n \in \N$, consider the $n$-indexed {\em scaling} map\hfff{scalingmap} $R_n:\R^2 \to \R^2$ given by
$$
 R_n \big(v_1,v_2 \big) = \Big( 2^{-1} n^{-2/3}( v_1 - v_2) \, , \,   v_2/n \Big) \, .
$$ 
 
The scaling map acts on subsets $C$ of $\R^2$ by setting
$R_n(C) = \big\{ R_n(x): x \in C \big\}$.

(Clearly, $n$ must be positive. In fact, $\N$ will denote $\{1,2,\cdots \}$ throughout.)

\subsubsection{Scaling transforms staircases to zigzags.}
The image of any staircase under $R_n$
will be called an $n$-zigzag~\hfff{zigzag}. The starting and ending points of an $n$-zigzag $Z$ are defined to be the image under $R_n$
of such points for the staircase $S$ for which $Z = R_n(S)$.
 
Note that the set of horizontal lines is invariant under $R_n$, while vertical lines are mapped to lines of gradient  $- 2 n^{-1/3}$.
As such, an $n$-zigzag is the range of a piecewise affine path from the starting point to the ending point which alternately moves rightwards  along horizontal line segments  and northwesterly along sloping line segments, where each sloping line segment has gradient  $- 2 n^{-1/3}$; the first and last segment in this journey may be either horizontal or sloping.

 \subsubsection{Scaled geodesics are called polymers.}
 For $n \in \N$, the image of any geodesic under the scaling map $R_n$ will be called an $n$-polymer, or often simply a polymer~\hfff{polymer}.  This usage of the term `polymer' for `scaled geodesic' is apt for our study, due to the central role played by these objects. The usage is not, however, standard: the term `polymer' is often used to refer to typical realizations of the path measure in LPP models at positive temperature.

\subsubsection{Some basic notation}
For $\ell \geq 1$, we write $\R^\ell_\leq$
for the subset of $\R^\ell$ whose elements $(z_1,\cdots,z_\ell)$
are non-decreasing sequences. When the sequences are increasing, we instead write $\R^\ell_<$. We also use the notation $A^\ell_\leq$ and $A^\ell_<$.
Here, $A \subset \R$ and the sequence elements are supposed to belong to $A$.
We will typically use this notation when $\ell=2$.
 
 \subsubsection{Compatible triples.}
 Let $(n,t_1,t_2) \in \N \times \R^2_<$, which is to say that $n \in \N$ and $t_1,t_2 \in \R$ with $t_1 < t_2$.
 Taking $x,y \in \R$, does there exist an $n$-zigzag from $(x,t_1)$ to $(y,t_2)$?
 As far as the data $(n,t_1,t_2)$ is concerned, such an $n$-zigzag may exist only if 
 \begin{equation}\label{e.ctprop}
     \textrm{$t_1$ and $t_2$ are integer multiplies of $n^{-1}$} \, .
\end{equation}
We say that data $(n,t_1,t_2)  \in \N \times \R^2_<$ is a {\em compatible triple} \hfff{comptriple} if it verifies the last condition.

An important piece of notation associated to a compatible triple is $\tot$~\hfff{tot}, which we will use to denote the difference $t_2 - t_1$. The law of the underlying Brownian ensemble $B: \Z \times \R \to \R$ is invariant under integer shifts in the first, curve indexing, coordinate. This translates to a distributional invariance of scaled objects under vertical shifts by multiples of $n^{-1}$, something that makes the parameter $\tot$
of far greater significance than $t_1$ or $t_2$.

\begin{figure}[ht]
\begin{center}
\includegraphics[height=7cm]{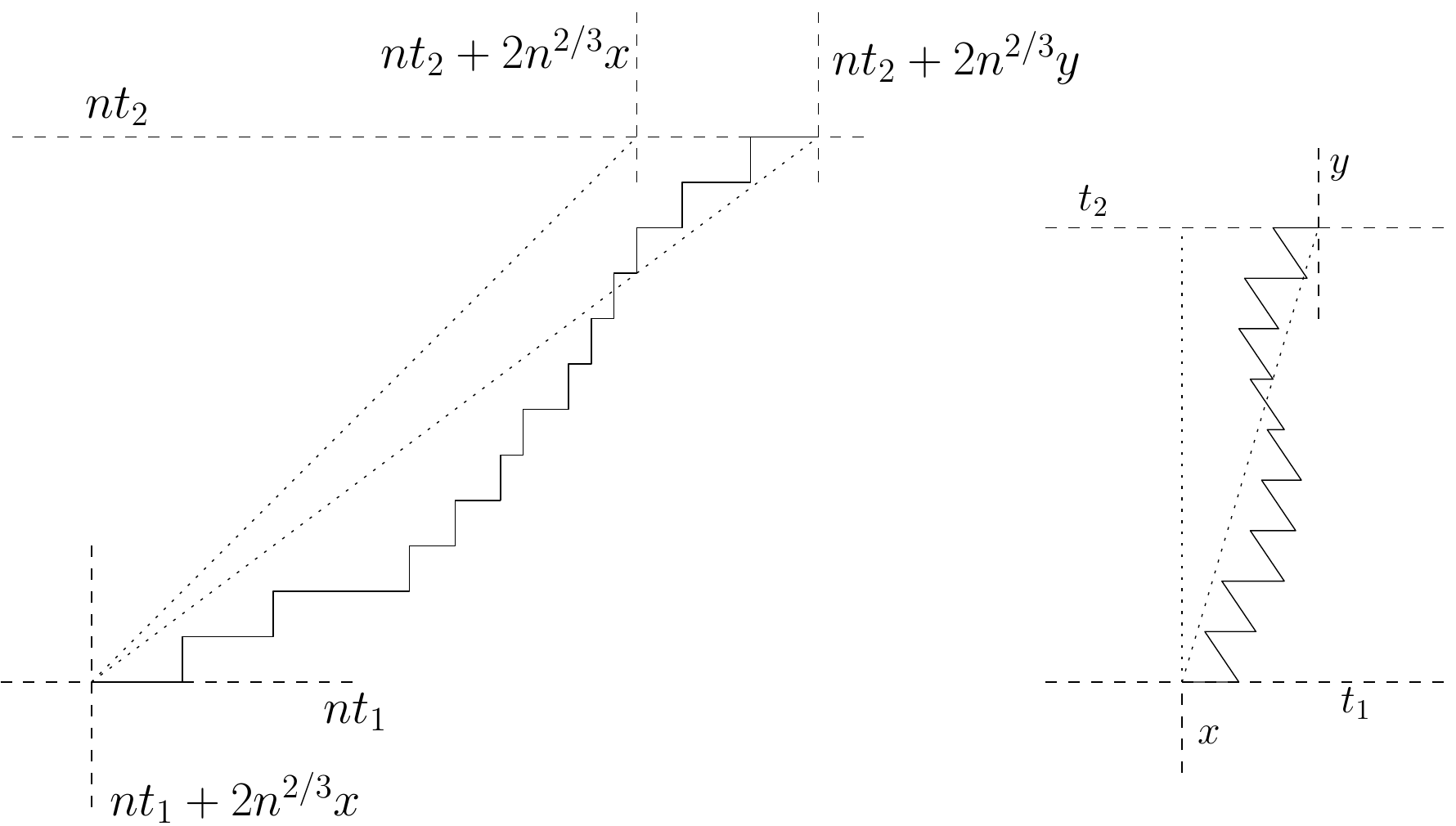}
\caption{Let $(n,t_1,t_2)$ be a compatible triple and let $x, y \in \R$. The endpoints of the geodesic in the left sketch have been selected so that, when the scaling map~$R_n$ is applied to produce the right sketch,
the $n$-polymer $\rho_{n;(x,t_1)}^{(y,t_2)}$ results.}
\label{f.scaling}
\end{center}
\end{figure}

Supposing now that $(n,t_1,t_2)$ is indeed a compatible triple, 
the condition 
 $y \geq x - 2^{-1} n^{1/3} \tot$ ensures that the preimage of $(y,t_2)$ under the scaling map $R_n$ lies northeasterly of the preimage of $(x,t_1)$.
 Thus, an  $n$-zigzag from $(x,t_1)$ to $(y,t_2)$ exists in this circumstance. We have mentioned that geodesics exist uniquely for given endpoints almost surely.
 Taking the image under scaling, this translates to the almost sure existence and uniqueness of the $n$-polymer from  $(x,t_1)$ to $(y,t_2)$. This polymer~\hfff{polynot} will be denoted $\rho_{n;(x,t_1)}^{(y,t_2)}$: see Figure~\ref{f.scaling}. 
This notation has characteristics in common with several later examples, in which objects are described in scaled coordinates.
Round bracketed expressions in the subscript or superscript will refer to  
a space-time pair, with the more advanced time in the superscript.
Typically some aspect of the $n$-polymer from $(x,t_1)$ to $(y,t_2)$ is being described when this ${\mathsf{TBA}}_{n;(x,t_1)}^{(y,t_2)}$ notation is used.  Appendix~\ref{s.glossary}, which offers a glossary of the principal notation in this article, records several examples in rough accordance with this convention.

\subsection{Main results, and a conjecture}

This article reaches conclusions in two main directions: it establishes upper bounds on the probability of the coexistence of disjoint polymers with nearby endpoints, and it gives expression to the two-thirds power law that dictates polymer geometry, showing that polymers fluctuate by order~$\e^{2/3}$ on a short scale~$\e$. In the next two subsections, the principal conclusions in these directions are respectively stated. 

\subsubsection{The rarity of many disjoint polymers}\label{s.maxpoly}

Let $(n,t_1,t_2) \in \N \times \R^2_<$ be a compatible triple, and let $I,J \subset \R$ be intervals. Set $\maxpoly_{n;(I,t_1)}^{(J,t_2)}$ \hfff{maxcard} equal to  the maximum cardinality of a {\em pairwise disjoint} set of   $n$-polymers each of whose starting and ending points have the respective forms $(x,t_1)$ and $(y,t_2)$ 
where $x$ is some element of $I$ and $y$ is some element of $J$.
% in this way, the event that $\maxpoly_{n;(I,t_1)}^{(J,t_2)} \geq k$ for given $k \in \intint{n}$ is characterized by the existence of vectors   $\bar{u} \in I^k_<$ and $\bar{v} \in J^k_<$ for which the polymer collection $\big\{ \rho_{n;(u_i,t_1)}^{(v_i,t_2)}: i \in \llbracket 1,k \rrbracket  \big\}$ is disjoint. (But polymer uniqueness should be addressed.)  

%For $k \in \N$, $\bar{u} \in [0,1]^k_<$
%and $\bar{v} \in [0,1]^k_<$, consider the collection of $n$-polymers $\rho_{n,(u_i,0)}^{(v_i,1)} \subset \R \times [0,1]$ for $i \in \intint{k}$.
%The maximum disjoint $n$-polymer cardinality $M_n$ is defined to be the maximal value of $k \in \N$
%for which such vectors $\bar{u}$ and $\bar{v}$
%exist such that this collection is pairwise disjoint.

(We make two parathetical comments. First, our notation here is an extension of the usage described a moment ago. In this case, the first element in the space-time pairs $(I,t_1)$ and $(J,t_2)$ is an interval, rather than a point. Second, the condition of pairwise disjointness in this definition can be slightly weakened so that our assertions remain valid. We prefer to defer elaborating on this.
A remark in Section~\ref{s.closure} will offer an altered definition of  $\maxpoly_{n;(I,t_1)}^{(J,t_2)}$ whose use results in slightly stronger results.)
% In specifying an event with such notation, our convention is to denote an intersection of events associated to the pairs $(x,t_1)$ and $(y,t_2)$ as $(x,y)$ varies over $I \times J$.

Our first  main conclusion treats the  problem of polymers  of shared unit-order lifetimes that cross between intervals that are of a small length~$\e$, albeit one that is independent of the scaling parameter~$n$. Theorem~\ref{t.disjtpoly.pop} provides an upper bound of $\e^{(k^2 - 1)/2 \, + \, o(1)}$ on the probability that~$k$ such polymers may coexist disjointly; here, $k \in \N$ is fixed, and the result is asserted uniformly in high~$n$.

%There exist two sequences $\big\{ C_k : k \in \N \big\}$
%and  $\big\{ c_k : k \in \N \big\}$ of positive constants that satisfy $\limsup C_k^{1/k^2} < \infty$ and $\liminf c_k > 0$.

\begin{theorem}\label{t.disjtpoly.pop}
There exists a positive  constant $\Cpop$ such that the following holds.
Let  $(n,t_1,t_2) \in \N \times \R^2_<$ be a compatible triple. 
Let  $k \in \N$, $\e > 0$ and $x,y \in \R$ satisfy the conditions that $k \geq 2$,
$$
 \e \leq    \Cpop^{-4k^2}   \, , \, n \tot \geq   \Cpop^{k^2}  \big( 1 +   \vert x - y \vert^{36}  \tot^{-24}   \big) \e^{-\Cpop} 
$$
and  
$\tot^{-2/3} \vert x - y \vert  \leq \e^{-1/2} \big( \log \e^{-1} \big)^{-2/3} G^{-k}$.

Setting $I = [x-\tot^{2/3}\e,x+\tot^{2/3}\e]$ and $J = [y-\tot^{2/3}\e,y+\tot^{2/3}\e]$, we have that  
$$
\PP \Big( \maxpoly_{n;(I,t_1)}^{(J,t_2)}  \geq k \Big) 
  \leq  
 \e^{(k^2 - 1)/2} \cdot \correctmac \, ,
$$
where $\correctmac$ is a positive correction term that is bounded above by $\Cpop^{k^3}  \exp \big\{ \Cpop^k \big( \log \e^{-1} \big)^{5/6} \big\}$.
\end{theorem}
Theorem~\ref{t.disjtpoly.pop} is a slightly simplified form of a more detailed result that will be presented later, Theorem~\ref{t.disjtpoly}.
Despite this simpler form, it is worth mentioning that, since a two-thirds power governs the horizontal coordinate under the scaling transformation, there is no generality lost by considering the special case $t_1 = 0$ and $t_2 = 1$, and thus $\tot = 1$, in Theorem~\ref{t.disjtpoly.pop}. (This point is explained further when the  {\em scaling principle} is discussed in Section~\ref{s.scalingprinciple}.)
The reader is thus encouraged to read this result with $t_1 =0$ and $t_2 = 1$, as well as with $x$ and $y$ both supposed to be at most one say. 
%Fixing $k \in \N$, we may take $\e > 0$ small but fixed in terms of $k$, and find that the result holds uniformly in high~$n$.
Indeed, doing so shows that
the next result follows immediately.
% from Theorem~\ref{t.disjtpoly.pop}.
\begin{corollary}\label{c.forconjecture}
For each $k \in \N$, $k \geq 2$, we have that
$$
 \liminf_{\e \searrow 0} \, \liminf_n \, \frac{\log \PP \Big( \maxpoly_{n;([- \e,\e],0)}^{([-\e,\e],1)}  \geq k \Big)}{\log \e}  \, \geq \, \frac{k^2 - 1}{2} \, .
 $$
 \end{corollary}
 Note that, despite the $\geq$ symbol, this result is an {\em upper bound} on the concerned probability.
The formulation permits us to express a conjecture in a symmetric form which indicates a belief that the exponent~$(k^2-1)/2$ in Theorem~\ref{t.disjtpoly.pop} is sharp.
\begin{conjecture}\label{c.two}
For each $k \in \N$, $k \geq 2$, we have that
$$
 \limsup_{\e \searrow 0} \, \limsup_n \, \frac{\log \PP \Big( \maxpoly_{n;([- \e,\e],0)}^{([-\e,\e],1)}  \geq k \Big)}{\log \e} \, \leq \, \frac{k^2 - 1}{2} \, .
 $$
\end{conjecture}
If the conjecture is valid, then equality will hold in the corollary and the conjecture. As we have mentioned, Conjecture~\ref{c.two} has recently been validated in~\cite[Theorem~$2.4$]{BGH18} for the case that $k=2$.

The second main conclusion is a rather direct consequence of the first. It asserts that, for polymers that cross between unit-length intervals separated by unit-order times, the maximum number of such disjoint polymers has a tail that decays super-polynomially. 
\begin{theorem}\label{t.maxpoly.pop}
There exist  constants $\Cpop \geq \cpop > 0$ such that the following holds.
Let
$(n,t_1,t_2) \in \N \times \R^2_<$ be a compatible triple. Further let $x,y \in \R$, $h \in \N$ 
%$h \geq 1$, 
and $\emm \in \N$. Suppose that 
$$
 \emm \geq  (\Cpop  h)^{2\cpop ^{-1}\Cpop } \vee  \big( \vert x - y \vert \tot^{-2/3} + 2h \big)^3 
$$
 and $n \tot  \geq  \max \big\{ 1, ( \vert x - y \vert \tot^{-2/3} + 2h )^{36} \big\} \Cpop  \emm^\Cpop$.
Then
$$
 \PP \Big(  \maxpoly_{n;(I,t_1)}^{(J,t_2)} \geq \emm \Big)  \, \leq \,  \emm^{-  \cpop (\log \log \emm)^2}     \, ,
$$
where $I$ denotes the interval $[x,x+h\tot^{2/3}]$ and $J$ the interval $[y,y+h \tot^{2/3}]$.
\end{theorem}
 Focus on the meaning of this theorem is brought by considering $t_1 =0$, $t_2 = 1$, $h=1$, with $x$ and $y$ chosen to lie in a unit interval about the origin. In this case, the demanded lower bound on $\emm$ merely excludes an initial bounded interval, and the lower bound on $n$ takes the form $n \geq \Theta (\emm^G)$. The latter condition is not demanding given that much of the interest in the use of scaled coordinates lies either in the high $n$ limit or in statements made uniformly for high enough $n$. Theorem~\ref{t.maxpoly.pop} is such a statement; as seen by these parameter choices, it asserts that there is at most probability $\emm^{-g (\log \log \emm)^2}$ of $\emm$ disjoint polymers crossing a bounded region, uniformly in high $n$.

\subsubsection{Polymer fluctuation}\label{s.polyflucintro}

Our principal conclusion in this regard, Theorem~\ref{t.polyfluc}, will here be stated, after a few paragraphs in which we attend to to setting up and explaining some necessary notation and concepts. 
The theorem  is used to prove Theorem~\ref{t.disjtpoly.pop} and it is also needed in~\cite{Patch}.

Let $(n,t_1,t_2) \in \N \times \R^2_<$ be a compatible triple, and let $x,y \in \R$.
The polymer $\rho_{n;(x,t_1)}^{(y,t_2)}$ has been defined to be a subset of $\R \times [t_1,t_2]$ containing $(x,t_1)$ and $(y,t_2)$, but really as $n$ rises towards infinity, it becomes more natural to seek to view it as a random function that maps its lifetime $[t_1,t_2]$ to the real line.
In choosing to adopt this perspective, we will abuse notation: taking $t \in [t_1,t_2]$, we will speak of the value  $\rho_{n;(x,t_1)}^{(y,t_2)}(t) \in \R$, as if the polymer were in fact a function of $[t_1,t_2]$.
Some convention must be adopted to resolve certain microscopic ambiguities as we make use of this new notation, however.
First, we will refer to  $\rho_{n;(x,t_1)}^{(y,t_2)}(t)$ only when $t \in [t_1,t_2]$ satisfies $nt \in \Z$, a condition that ensures that the intersection of the set $\rho_{n;(x,t_1)}^{(y,t_2)}$  with the line at height $t$ takes place along a horizontal planar interval.

Second, we have to explain which among the points in this interval $\rho_{n;(x,t_1)}^{(y,t_2)} \cap \{  (\cdot,t): \cdot \in \R  \}$
 we wish to denote by $\rho_{n;(x,t_1)}^{(y,t_2)}(t)$.
 To present and explain our convention in this regard, we let
   $\ell_{(x,t_1)}^{(y,t_2)}$ denote the planar line segment whose endpoints are $(x,t_1)$ and $(y,t_2)$. Adopting the same perspective as for the polymer, we abuse notation to view  $\ell_{(x,t_1)}^{(y,t_2)}$ as a function from $[t_1,t_2]$ to $\R$, so that  $\ell_{(x,t_1)}^{(y,t_2)}(t) = \tot^{-1} \big( (t_2 - t) x + (t - t_1)y \big)$.
   
Our convention will be to set  $\rho_{n;(x,t_1)}^{(y,t_2)}(t)$ equal to~$z$ where $(z,t)$ is that point in the horizontal segment  $\rho_{n;(x,t_1)}^{(y,t_2)} \cap \{  (\cdot,t): \cdot \in \R  \}$ whose distance from   $\ell_{(x,t_1)}^{(y,t_2)}(t)$ is maximal. (An arbitrary tie-breaking rule, say $\rho_{n;(x,t_1)}^{(y,t_2)}(t) \geq  \ell_{(x,t_1)}^{(y,t_2)}(t)$, resolves the dispute if there are two such points.)
The reason for this very particular convention is that our purpose in using it is to explore, in the soon-to-be-stated Theorem~\ref{t.polyfluc},
upper bounds on the probability of large fluctuations between the polymer  $\rho_{n;(x,t_1)}^{(y,t_2)}$ and the line segment  $\ell_{(x,t_1)}^{(y,t_2)}$ that interpolates the polymer's endpoints.
Our convention ensures that the form of the theorem would remain valid were any other convention instead adopted.     

In order to study the intermediate time $(1-a)t_1 + at_2$ (in the role of $t$ in the preceding), we now  let $a \in (0,1)$ and impose that $a \tot n \in \Z$: doing so ensures that, as desired, $t \in n^{-1} \Z$, where $t = (1-a)t_1 + at_2$.

Consider also $r > 0$. Define the {\em polymer deviation regularity} \hfff{polydevreg} event 
\begin{equation}\label{e.pdr}
\pdr_{n;(x,t_1)}^{(y,t_2)}\big(a,r\big) \, = \, \bigg\{  \,  \Big\vert \, \rho_{n;(x,t_1)}^{(y,t_2)} \big(  (1-a) t_1 +  a t_2 \big) - \ell_{(x,t_1)}^{(y,t_2)}  \big(  (1-a) t_1 +  a t_2 \big) \, \Big\vert \, \leq \, r t_{1,2}^{2/3} \big(  a \wedge (1-a) \big)^{2/3} \bigg\} \, ,
\end{equation}
where $\wedge$ denotes minimum. For example, if $a \in (0,1/2)$,
the polymer's deviation from the interpolating line segment, at height $(1-a)t_1 + at_2$ (when the polymer's journey has run for time $a\tot$), is measured in the natural time-to-the-two-thirds  scaled units obtained by division by $(a \tot)^{2/3}$, and compared to the given value $r > 0$.

For intervals $I,J \subset \R$, we extend this definition by setting
$$
 \pdr_{n;(I,t_1)}^{(J,t_2)}\big(a,r\big)   = \bigcap_{x \in I, y \in J}  \pdr_{n;(x,t_1)}^{(y,t_2)}\big(a,r\big) \, .
$$ 
The perceptive reader may notice a problem with the last definition. 
The polymer $\rho_{n;(x,t_1)}^{(y,t_2)}$
is well defined almost surely for given endpoints, but this property is no longer assured as the parameters vary over $x \in I$ and $y \in J$.
In the case of exceptional $(x,y)$ where several $n$-polymers
move from $(x,t_1)$ to $(y,t_2)$, we interpret $\rho_{n;(x,t_1)}^{(y,t_2)}$ as the union of all these polymers, for the purpose of defining $\rho_{n;(x,t_1)}^{(y,t_2)}(t)$. This convention permits us to identify worst case behaviour, so that the event $\neg \, \pdr_{n;(I,t_1)}^{(J,t_2)}\big(a,r\big)$ is triggered by a suitably large fluctuation on the part of any concerned polymer. Here, and later, $\neg \, A$ denotes the complement of the event $A$.

We are ready to state our conclusion concerning polymer fluctuation. 
 Bounds in this and many later results have been expressed explicitly up to two positive constants $c$ and $C$. See Section~\ref{s.useful} for an explanation of how the value of this pair of constants is fixed. We further set  $c_1 = 2^{-5/2} c \wedge 1/8$.
\begin{theorem}\label{t.polyfluc}
Let $(n,t_1,t_2) \in \N \times \R^2_<$ be a compatible triple, and let $x,y \in \R$. 
\begin{enumerate}
\item
Let  $a \in \big[1 - 10^{-11} c_1^2 , 1 \big)$ satisfy $a \tot \in n^{-1} \Z$.  
 Suppose that $n \in \N$ satisfies 
$$
 n \tot   \geq  \max \bigg\{ 
 10^{32} (1-a)^{-25} c^{-18} \, \, , \, \, 
 10^{24} c^{-18} (1-a)^{-25} \vert x - y \vert^{36} \tot^{-24} 
 \bigg\} \, .
$$
Let $r > 0$ be a parameter that satisfies
$$ 
 r  \geq 
 \max \bigg\{  10^9 c_1^{-4/5} \, \, , \, \,  15 C^{1/2} \, \, , \, \,  87(1-a)^{1/3} \tot^{-2/3} \vert x - y \vert   \bigg\}
 $$
and
$r \leq 3  (1-a)^{25/9} n^{1/36} \tot^{1/36}$. 

Writing $I = \big[ x, x +  t_{1,2}^{2/3} (1-a)^{2/3} r \big]$
and  $J = \big[ y, y +  t_{1,2}^{2/3} (1-a)^{2/3} r \big]$, we have that
$$
\PP \Big( \neg \, \pdr_{n;(I,t_1)}^{(J,t_2)}\big(a,2r\big) \, \Big) \leq  
 44  C r    \exp \big\{ - 10^{-11} c_1  r^{3/4} \big\} 
 \, .
$$
%where  $R =  t_{1,2}^{2/3} (1-a)^{2/3} r$.
\item The same statement holds verbatim when appearances of $a$ are replaced by $1-a$.
\end{enumerate}
\end{theorem}

One aspect of our presentation may be apparent from the form of this theorem: we have chosen to be fairly explicit in recording hypothesis bounds in our results. This approach can provide quite lengthy formulas when hypotheses are recorded, and we encourage the reader not to be distracted by this from the essential meaning of results.
In the present case, for example, we may set $t_1 = 0$ and $t_1 = 1$, so that $\tot = 1$ results. We may note that $a$ must be a certain small distance from either zero or one, and that the condition $a \in n^{-1}\Z$ is a negligible constraint, given that we are interested in high choices of $n$. Working in a case where $x$ and $y$ are bounded above by one, we are stipulating that $n \geq n_0(a)$.
The parameter $r$, whose role is to gauge scaled fluctuation, is then constrained to lie between a universal constant and a multiple of $n^{1/36}$. We see then that the theorem is asserting that the maximum scaled fluctuation during the first, or last, duration $a$ of time among all polymers beginning and ending in a given unit interval, separated at a unit duration, exceeds $r \in [\Theta(1),\Theta(1)n^{1/36}]$
with probability at most $\exp \{ - O(1) r^{3/4} \}$.
Since we are interested in high $n$, and the bound so available when $r = \Theta(1) n^{1/36}$ is rapidly decaying in $n$, the condition imposed on $r$ is quite weak. 

If the reader is ever disconcerted by the various conditions imposed in results, it may be helpful to consider the formal case where $n = \infty$, in which upper bounds, such as $r \leq \Theta(1) n^{1/36}$ above,
become obsolete. 

As a matter of convenience for the upcoming proofs, we also state a version of Theorem~\ref{t.polyfluc} in which the intervals $I$ and $J$ are singleton sets.

\begin{proposition}\label{p.polyfluc}
Let $(n,t_1,t_2) \in \R^2_<$ be a compatible triple.  Let $x,y \in \R$ and   let $a \in \big[1 - 10^{-11} c_1^2 , 1 \big)$ satisfy $a \tot \in n^{-1} \Z$. Suppose that 
$$
 n \tot   \geq  \max \bigg\{ 
 10^{32} (1-a)^{-25} c^{-18} \, \, , \, \, 
 10^{24} c^{-18} (1-a)^{-25} \vert x - y \vert^{36} \tot^{-24} 
 \bigg\} \, .
$$
Let $r > 0$ be a parameter that satisfies
$$ 
 r  \geq 
 \max \bigg\{  10^9 c_1^{-4/5} \, \, , \, \,  15 C^{1/2} \, \, , \, \,  87(1-a)^{1/3} \tot^{-2/3} \vert x - y \vert   \bigg\}
 $$
and
$r \leq 3  (1-a)^{25/9} n^{1/36} \tot^{1/36}$. 
Then
$$
 \PP \Big(  \neg \, \pdr_{n;(x,t_1)}^{(y,t_2)}\big(a,r\big) \, \Big) \leq 
22  C r    \exp \big\{ - 10^{-11} c_1  r^{3/4} \big\}    \,  .
$$
\end{proposition}
\noindent{\em Remark.}
The argument leading to the  proposition will show that the result equally applies when it is instead supposed that $a \in (0,10^{-11}c_1^2]$, provided that the instances of $1-a$
in the hypothesis conditions are replaced by $a$.

\subsubsection{Acknowledgments.}
%The author thanks Riddhipratim Basu, Shirshendu Ganguly and Jeremy Quastel for valuable conversations at many stages of this project. 
%submit 
The author thanks Riddhipratim Basu, 
Ivan Corwin,  Shirshendu Ganguly and Jeremy Quastel for valuable conversations. He thanks a referee for thorough and very useful comments.  
%He is supported by NSF grant DMS-$1512908$.

\section{A road map for proving the rarity of disjoint polymers with close endpoints}

In this section, we present a very coarse overview of the strategy of the proof of our result, Theorem~\ref{t.disjtpoly.pop}, establishing the improbability of the event that several disjoint polymers begin and end in a common pair of short intervals. (Theorem~\ref{t.polyfluc}'s proof depends on different ideas and will appear towards the end of the paper.)
Before we begin, we first present a key concept: the scaled energy, or {\em weight}~\hfff{weight}, of an $n$-zigzag.

\subsection{Polymer weights}

\subsubsection{Staircase energy scales to zigzag weight.}
Let $n \in \N$ and $(i,j) \in \N^2_<$.
Any $n$-zigzag $Z$ from $(x,i/n)$ to $(y,j/n)$  is ascribed a scaled energy, which we will refer to as its weight, 
$\weight(Z) = \weight_n(Z)$, given by 
\begin{equation}\label{e.weightzigzag}
 \weight(Z) =  2^{-1/2} n^{-1/3} \Big( E(S) -  2(j - i)  - 2n^{2/3}(y-x) \Big) 
\end{equation}
where $Z$ is the image under $R_n$ of the staircase $S$.

\subsubsection{Maximum weight.}  Let $(n,t_1,t_2) \in \N \times \R^2_<$ be a compatible triple. 
Suppose that $x,y \in \R$ satisfy 
 $y \geq x - 2^{-1} n^{1/3} \tot$. Define~\hfff{maxweight}
%\begin{equation}\label{e.weightm}
$$
  \weight_{n;(x,t_1)}^{(y,t_2)} \,     =  \,   2^{-1/2} n^{-1/3} \Big(  M^1_{(n t_1 + 2n^{2/3}x,n t_1) \to (n t_2 + 2n^{2/3}y,n t_2)} - 2n \tot -  2n^{2/3}(y-x) \Big) \, .
$$
%\end{equation}
Thus, $\weight_{n;(x,t_1)}^{(y,t_2)}$ equals the maximum weight of any $n$-zigzag from $(x,t_1)$ to $(y,t_2)$. When the polymer $\rho_{n;(x,t_1)}^{(y,t_2)}$ is well defined, as it is almost surely, 
we have that $\weight_{n;(x,t_1)}^{(y,t_2)} = \weight \big( \rho_{n;(x,t_1)}^{(y,t_2)} \big)$. It is also worth noting, however, that the system of weights  $\weight_{n;(x,t_1)}^{(y,t_2)}$ is well defined almost surely, even as the five concerned parameters vary over admissible choices, despite the possible presence of exceptional parameter choices at which the corresponding polymer is not unique.

\subsection{The road map}\label{s.roadmap}

We begin by revisiting a theme from Subsection~\ref{s.polyflucintro}. When $n$ is large, an $n$-polymer such as $\rho_{n;(x,t_1)}^{(y,t_2)}$ is rather naturally viewed as a random function of its lifetime $[t_1,t_2]$, one that begins at $x$ and ends at $y$. The microscopic ambiguities in the specification of this random function become vanishingly small in the limit of high $n$. Indeed, for the purposes of these paragraphs of overview, we yield to the temptation to set $n=\infty$, and discuss polymers as objects arising after a scaling limit has been taken. In this way, $\rho_{\infty;(x,t_1)}^{(y,t_2)}$ is interpreted as the scaled limit polymer between space-time locations $(x,t_1)$ and $(y,t_2)$. The polymer may be viewed as a subset of $\R^2$, living inside the strip $\R \times [t_1,t_2]$, or it may be viewed as a random, real-valued, function of the lifetime $[t_1,t_2]$.  Speaking of such objects as $\rho_{\infty;(x,t_1)}^{(y,t_2)}$ raises substantial questions about the uniqueness of a limiting description as $n$ tends to infinity. 
In this non-rigorous overview, we have no intention of trying to address these questions. We simply use the framework of $n = \infty$ as a convenient device for heuristic discussion, since this framework is unencumbered by microscopic details (such as the structure of zigzags). The term zigzag has become rather inapt in this context; it may be replaced by the term {\em path}, a path being any continuous real-valued function defined on a given interval $[t_1,t_2]$. Each path has a weight and a path of maximum weight given its endpoints is a polymer. When these endpoints are $(x,t_1)$ and $(y,t_2)$, this maximum weight is called $\weight_{\infty;(x,t_1)}^{(y,t_2)}$. 

\begin{figure}[ht]
\begin{center}
\includegraphics[height=9cm]{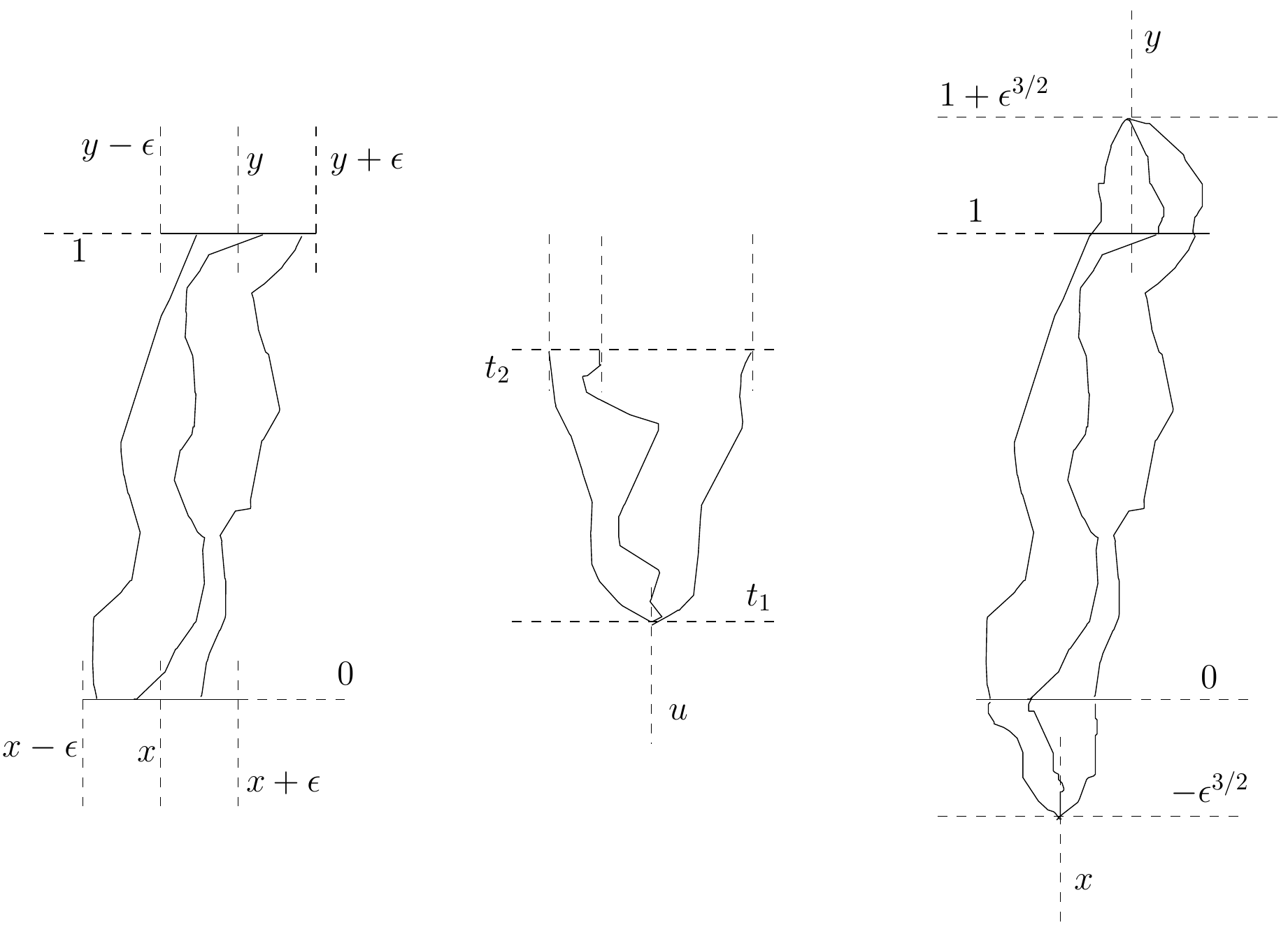}
\caption{{\em Left:}  The three polymers concerned in the formally specified event $\maxpoly_{\infty;([x-\e,x+\e],0)}^{([y-\e,y+\e],1)} \geq 3$. {\em Middle:}  A bouquet of three paths of lifetime $[t_1,t_2]$. {\em Right:} Two bouquets, facing in opposite directions, are employed to extend the original three polymers into longer paths, each beginning and ending at $(x,-\e^{3/2})$ and $(y,1 +\e^{3/2})$.}
\label{f.triple}
\end{center}
\end{figure}

Viewed through the prism of $n = \infty$, the event that Theorem~\ref{t.disjtpoly.pop} discusses is depicted in the left sketch of Figure~\ref{f.triple}. 
Here we take $(n,t_1,t_2) = (\infty,0,1)$. Depicted with $\kay=3$ is the event  
 $\maxpoly_{\infty;([x-\e,x+\e],0)}^{([y-\e,y+\e],1)}  \geq \kay$ that there exist $\kay$ disjoint polymers with lifetime $[0,1]$ each of which begins at distance at most $\e$ from $x \in \R$ and ends at such a distance from $y \in \R$. 
 %(It would have been in keeping with the notation of Theorem~\ref{t.disjtpoly.pop} to have used $k$ in place of $m$ here. We will in fact reserve $k$ for its use in Theorem~\ref{t.maxpoly.pop} and will henceforth use $m$ in place of $k$ in regard to Theorem~\ref{t.disjtpoly.pop}.)
 
 It is valuable to bear in mind three general features of polymer geometry and weight that correspond to powers of $2/3$, $1/3$ and $1/2$.
 The first of these three principles has been articulated rigorously by Theorem~\ref{t.polyfluc}. For the latter two, we will recall in Section~\ref{s.input} corresponding results for Brownian LPP from~\cite{ModCon}.

 \noindent{\em A power of two-thirds dictates polymer geometry.}  A polymer  whose lifetime is $[t_1,t_2]$, and thus has duration $\tot$, 
 fluctuates laterally, away from the planar line segment that interpolates its endpoints, by an order of $\tot^{2/3}$.

 \noindent{\em  A power of one-third dictates polymer weight.} A polymer whose lifetime is $[t_1,t_2]$ has a weight of order $\tot^{1/3}$. 
 Actually, this is only true if the polymer makes no significant lateral movement. For example, $\rho_{\infty;(x,t_1)}^{(y,t_2)}$ may be expected to have a weight of order $\tot^{1/3}$
 provided that the endpoints verify $\vert y - x \vert = O(\tot^{2/3})$. The last condition indicates that the endpoint discrepancy is of the order of the polymer fluctuation, so that such a polymer is not deviating more than would be expected by a polymer whose starting and ending points coincide.
 
 \noindent{\em A power of one-half dictates polymer weight differences.} Consider two polymers of unit duration both of whose endpoints differ by a small quantity $\e$. 
 Then the polymers' weights typically differ by an order of $\e^{1/2}$. For example, $\weight_{\infty;(x + \e,0)}^{(y + \e ,1)} - \weight_{\infty;(x ,0)}^{(y ,1)}$ typically has order $\e^{1/2}$, at least when $\vert y - x \vert = 0(1)$.

We now introduce an event that has much in common with the event  that $\maxpoly_{\infty;([x-\e,x+\e],0)}^{([y-\e,y+\e],1)}$  is at least $\kay$. 
Let $\eta > 0$. The {\em near polymer}~\hfff{nearpoly} event  $\nearpoly_{\infty,\kay;(x,0)}^{(y,1)}(\eta)$ is said to occur when there exist $\kay$ paths mapping $[0,1]$ to $\R$, each beginning at $x \in \R$ and ending at $y \in \R$, 
but which are pairwise disjoint except at the endpoint locations,
such that the sum of the paths' weights exceeds $\kay \cdot \weight_{\infty;(x,0)}^{(y,1)} \, - \, \eta$.  To understand the meaning of this definition, note that, since each of these $\kay$ paths has a weight bounded above by that of the polymer with these given endpoints, the sum of the paths' weights may be at most $\kay \cdot \weight_{\infty;(x,0)}^{(y,1)}$. Since weights are of unit order, due to a unit duration being considered, we see that, when $\eta > 0$ is fixed to be a small value, the event  $\nearpoly_{\infty,\kay;(x,0)}^{(y,1)}(\eta)$ expresses the presence of a system of $\kay$ near polymers, with shared endpoints but otherwise disjoint, whose sum weight misses the maximum value attainable in principle by a smaller than typical discrepancy of $\eta$.

The backbone of our derivation of Theorem~\ref{t.disjtpoly.pop} consists of arguing that the occurrence of the event $\maxpoly_{\infty;([x-\e,x+\e],0)}^{([y-\e,y+\e],1)} \geq \kay$
typically entails that  $\nearpoly_{\infty,\kay;(x,0)}^{(y,1)}(\eta)$ also occurs, where here $\eta$ is set to be of the order of $\e^{1/2}$. Before we offer an overview of why this assertion may be expected to be true, we point out that the order of probability of the latter event is understood. Indeed, the use of shared endpoints in the definition of  $\nearpoly_{\infty,\kay;(x,0)}^{(y,1)}(\eta)$ permits the use of integrable techniques arising from the Robinson-Schensted-Knuth correspondence: this event is shown in~\cite{BrownianReg} to have probability of order $\eta^{\kay^2 - 1 + o(1)}$. The result in question will shortly be reviewed as Corollary~\ref{c.neargeod.t}. In reality, the result actually concerns the events   $\nearpoly_{n,k;(x,0)}^{(y,1)}(\eta)$ when $n \in \N$ is finite, and makes an assertion uniformly in high $n$ about them. 
Anyway, since we are taking $\eta = \Theta(\e^{1/2})$, we see that the form of Theorem~\ref{t.disjtpoly.pop} may be expected to follow from such a bound.

Now, how is it that we may argue that the $\kay$ disjoint polymer event $\maxpoly_{\infty;([x-\e,x+\e],0)}^{([y-\e,y+\e],1)} \geq \kay$ does indeed typically entail 
 $\nearpoly_{\infty,\kay;(x,0)}^{(y,1)}(\eta)$ with $\eta = \Theta \big( \e^{1/2} \big)$? The first event provides much of the structure of the $\kay$ near polymers that would realize the second. 
 But in the second, the near polymers begin and end at precisely the same points, rather than merely nearby one another. What is needed is surgery to alter the $\kay$ polymers that begin and end in $[x-\e,x + \e] \times \{ 0 \}$ and $[y - \e,y + \e] \times \{ 1 \}$ so that their endpoints are actually equal. 
 
 Our device enabling this surgery to effect local correction of endpoint locations is called a {\em multi-polymer bouquet}. Let $[t_1,t_2]$ be a short real interval, and let $u \in \R$. Fixing an increasing list of $\kay$ real values, $(v_1,\cdots,v_\kay)$, each of unit order, we may consider the problem of constructing a system of $\kay$ paths, each with lifetime $[t_1,t_2]$, with each path beginning at $(u,t_1)$, and with the $i\textsuperscript{th}$ path ending at $\big(u + \tot^{2/3} v_i , t_2 \big)$; crucially, we insist that the paths be disjoint, except for their shared starting location. Such a bouquet of paths is depicted in the middle sketch in Figure~\ref{f.triple}. These paths are moving on the natural two-thirds power fluctuation scale given their endpoints, so it is a natural belief that they may be selected so that the sum of their weights has order $\tot^{1/3}$, which is the natural polymer weight scale for polymers of duration $\tot$. 
  It is one of the upcoming challenges for the proof of Theorem~\ref{t.disjtpoly.pop} to prove that the construction of such a path system may be undertaken. We will refer to this challenge as the problem of  {\em bouquet construction}.

 Accepting for now that such a bouquet of $\kay$ paths, tied together at a common base point, may be constructed, it is of course also natural to believe that such a construction may be undertaken where instead it is at the endpoint time that the various paths share a common location. The two constructions may be called a forward and a backward bouquet. 
 
 The two bouquets may be used to tie together the nearby endpoints of the $\kay$ polymers that arise in the event  $\maxpoly_{\infty;([x-\e,x+\e],0)}^{([y-\e,y+\e],1)} \geq \kay$. 
 Rather than altering the paths so that they begin and end at $(x,0)$ and $(y,1)$, we will prefer to extend the duration of the paths a little in each direction. Because the spatial discrepancy of endpoint locations is of order $\e$, the two-thirds power law for polymer geometry indicates that the natural time scale for our polymer bouquets will be $\e^{3/2}$. 
 A forward bouquet of $\kay$ paths will be used at the start and a backward one at the end. The $\kay$ paths in the forward bouquet will each begin at $(x,-\e^{3/2})$; their other endpoints, at time zero, will be the $\kay$ starting points of the duration $[0,1]$ polymers. The backward bouquet will take the $\kay$ ending points of the polymers and tie them together at the common location $(y,1+\e^{3/2})$.
 
 In this way, a system of $\kay$ paths from $(x,-\e^{3/2})$ to $(y,1+\e^{3/2})$ will be obtained: see the right sketch in Figure~\ref{f.triple}. Each is split into three pieces: a short $\e^{3/2}$-duration element of the forward bouquet, a long $[0,1]$-time piece given by one of the polymers, and then a further short element of the backward bouquet. This system will be shown to typically  realize  
 $\nearpoly_{\infty,\kay;(x,-\e^{3/2})}^{(y,1 + \e^{3/2})}(\eta)$ with $\eta = \Theta \big( \e^{1/2} \big)$. Note that the starting and ending times are $-\e^{3/2}$ and $1 + \e^{3/2}$, rather than zero and one. 
 %(as we at first suggested).
 
 We now summarise the challenges involved in proving this. The value $n$ is finite, with statements to be understood uniformly in high choices of this parameter. 
 \begin{enumerate}
 \item Polymer weight similarity: if the construction is to result in a collection of $\kay$ paths from $(x,-\e^{3/2})$ to $(y,1 + \e^{3/2})$ each of whose weights reaches within order $\e^{1/2}$
 of the maximum achievable value   $\weight_{n;(x,-\e^{3/2})}^{(y,1 + \e^{3/2})}$, then we will want to begin by confirming that the difference in weight between any pair of the $\kay$ polymers moving between $[x - \e,x+\e] \times \{0\}$ and $[y - \e,y+\e] \times \{ 1 \}$ itself has order $\e^{1/2}$. This is an instance of the one-half power law principle enunciated above.
 \item Bouquet construction: the forward and backward bouquets, defined on the time intervals $[-\e^{3/2},0]$ and $[1,1+\e^{3/2}]$, must be constructed so that each bouquet has cumulative weight of the desired order~$\e^{1/2}$.
 \item Final polymer comparison: even if our original polymers have similar weights, and they can be altered to run from $(x,-\e^{3/2})$ to $(y,1+\e^{3/2})$ with weight changes of order $\e^{1/2}$, there could still be trouble in principle. We need to know that each of the extended paths has a weight suitably close to  $\weight_{n;(x,-\e^{3/2})}^{(y,1 + \e^{3/2})}$. Thus, we need to verify that the polymer running from  $(x,-\e^{3/2})$ to $(y,1+\e^{3/2})$ has a weight that exceeds the extended paths' weights by an order of at most $\e^{1/2}$. 
 \end{enumerate} 
 In a sense, the third of these challenges is simply a restatement of the entire problem. When we discuss it later, however, we will resolve the problem in light of the solutions of the first two difficulties.
 
 This completes our first overview of the road map to the proof of Theorem~\ref{t.disjtpoly.pop}. We have been left with three challenges to solve. 
 In the next three sections, we will present definitions and tools and cite results that will be used to implement the road map, as well as to prove the polymer fluctuation Theorem~\ref{t.polyfluc}.
 After these tools have been described, we will describe at the end of Section~\ref{s.input} the structure of the rest of the paper. In essence, we will then return to the road map, elaborate further in outline how to resolve its three challenges, and implement it rigorously.

\section{Tools: staircase collections and multi-polymers}

The definition and study of the event $\nearpoly_{n,\kay;(x,0)}^{(y,1)}(\eta)$ (for $\kay$ fixed and $n$, although now finite, high in applications), and the problem of bouquet construction, will both involve the use of systems of staircases or zigzags. We begin this discussion of tools by setting up definitions in this regard.

\subsection{Staircase collections}\label{s.staircase}

Two staircases are called horizontally separate if there is no planar horizontal interval of positive length that is a subset of a horizontal interval in both staircases.

We now introduce notation for collections of pairwise horizontally separate staircases. 
For $k \in \N$, let $(x_i,s_i)$ and $(y_i,f_i)$, $i \in \intint{k}$, be a collection of pairs of elements of $\R \times \N$. (The symbols $s$ and $f$ are used in reference to the staircases' heights at the {\em start} and {\em finish}.)

Let $\staircase^k_{(\bar{x},\bar{s}) \to (\bar{y},\bar{f})}$
denote the set of $k$-tuples $(\phi_1,\cdots,\phi_k)$,
where $\phi_i$ is a staircase from $(x_i,s_i)$ to $(y_i,f_i)$
and each pair $(\phi_i,\phi_j)$, $i \not= j$, is  horizontally separate. 
(This set may be empty; in order that it be non-empty, it is necessary that  $y_i \geq x_i$ and $f_i \geq s_i$
for  $i \in \intint{k}$.) 
Note also that 
$\staircase^1_{(x_1,s_1) \to (y_1,f_1)}$ equals 
$\staircase_{(x_1,s_1) \to (y_1,f_1)}$. 

We also associate an energy to each member of  $\staircase^k_{(\bar{x},\bar{s}) \to (\bar{y},\bar{f})}$.
Each of the $k$ elements of any $k$-tuple in $\staircase^k_{(\bar{x},\bar{s}) \to (\bar{y},\bar{f})}$ has an energy, as we described in Subsection~\ref{s.staircases}. Define the energy $E\big( \phi \big)$ of any $\phi = \big( \phi_1,\cdots,\phi_k \big) \in \staircase^k_{(\bar{x},\bar{s}) \to (\bar{y},\bar{f})}$ to be $\sum_{j=1}^k E(\phi_j)$.

When  $\upright^k_{(\bar{x},\bar{s}) \to (\bar{y},\bar{f})} \not= \emptyset$, we further define the maximum $k$-tuple energy
\begin{equation}\label{e.mell}
   M^k_{(\bar{x},\bar{s}) \to (\bar{y},\bar{f})}  = \sup \Big\{ E(\phi): \phi \in  \upright^k_{(\bar{x},\bar{s}) \to (\bar{y},\bar{f})} \Big\} \, .
\end{equation}
%When $k = 1$, a staircase attaining the maximum is called a geodesic from $(x_1,s_1)$ to~$(y_1,f_1)$.
An $k$-tuple of staircases that attains this maximum may be called a multi-geodesic.

\subsection{Maximizer uniqueness}\label{s.maxunique}

Next is~\cite[Lemma~$A.1$]{Patch}. Since this result is drawn the final article in our four-paper study, it is worth mentioning that the result has a fairly short and self-contained proof in \cite[Appendix~$A$]{Patch}.
\begin{lemma}\label{l.severalpolyunique}
Let $k \in \N$ and let $(x_i,s_i)$ and $(y_i,f_i)$, $i \in \intint{k}$, be a collection of  pairs of points in $\R \times \N$
%that satisfy $(x_i,s_i) \preceq (y_i,f_i)$ for each
% $i \in \intint{k}$.
such that  $\upright^k_{(\bar{x},\bar{s}) \to (\bar{y},\bar{f})}$ is non-empty.
Then, 
except on a $\PP$-null set, there is a unique element of
$\upright^k_{(\bar{x},\bar{s}) \to (\bar{y},\bar{f})}$
whose energy attains $M^k_{(\bar{x},\bar{s}) \to (\bar{y},\bar{f})}$. 
\end{lemma}

When it exists, we denote the unique maximizer by
%whose almost sure uniqueness has just been asserted by
$$
 \Big( P^{k}_{(\bar{x},\bar{s}) \to (\bar{y},\bar{f});i} :  i \in \intint{k} \Big) \in  \upright^k_{(\bar{x},\bar{s}) \to (\bar{y},\bar{f})} \, .
$$

%Let $\preceq$ denote 
Consider the partial order on $\R^2$
such that $(x_1,x_2)$ is at most $(y_1,y_2)$
if and only if $x_1 \leq x_2$ and $y_1 \leq y_2$.

In the case that $k = 1$, and $(x,s)$ and $(y,f)$ are elements of $\R \times \N$ with $(y,f)$ exceeding $(x,s)$ in the partial order, the set 
$\upright^1_{(x,s) \to (y,f)}$ is non-empty. The maximizer in Lemma~\ref{l.severalpolyunique} is in this special case denoted by 
$P^{1}_{(x,s) \to (y,f)}$.

\subsection{Maximum weights of zigzag systems}\label{s.zigzagmax}

Each $n$-zigzag is the image under the scaling map~$R_n$ of a staircase.
Two such zigzags will be called horizontally separate if their staircase counterparts have this property. The weight of a collection of zigzags is the sum of the weights of those zigzags.

We wish to record notation for the maximum weight of a pairwise horizontally disjoint collection of $[t_1,t_2]$-lifetime zigzags with prescribed endpoints. 

To this end, let  $(n,t_1,t_2) \in \N \times \R^2_\leq$ be a compatible triple.
Taking $k \in \N$, a parameter that will denote the number of zigzags,  consider $k$ pairs $(x_i,y_i) \in \R^2$, $i \in \intint{k}$, where we suppose that $y_i \geq x_i - 2^{-1} n^{1/3} \tot$
for each $i \in \intint{k}$. Further suppose that the vectors $\overline{x} = \big( x_1,\cdots,x_k \big)$ and $\overline{y} = (y_1,\cdots,y_k)$ are non-decreasing. We will write $\weight_{n,k;(\bar{x},t_1)}^{(\bar{y},t_2)}$ for the maximum weight associated to $k$ horizontally separate zigzags moving consecutively between $(x_i,t_1)$ to $(y_i,t_2)$ for $i \in \intint{k}$.
Formally, we define~\hfff{collectiveweight}
\begin{equation}\label{e.multiweight}
\weight_{n,k;(\bar{x},t_1)}^{(\bar{y},t_2)} \,     =  \,   2^{-1/2} n^{-1/3} \Big(  M^k_{(n t_1 \bar{\bf 1} + 2n^{2/3}\bar{x},n t_1 \bar{\bf 1}) \to (n t_2 \bar{\bf 1} + 2n^{2/3}\bar{y},n t_2 \bar{\bf 1}  )} - 2 k n \tot -  2n^{2/3}\sum_{i=1}^k(y_i-x_i) \Big) \, ,
\end{equation}
where $\bar{\bf 1}\in \R^k$ denotes the vector whose components equal $1$. 

\subsection{Multi-polymers and their weights}\label{s.multipolymer}

Retaining these parameters, we now may offer a geometric view of this multi-weight. The maximum energy~$M^k$ seen on the right-hand side of~(\ref{e.multiweight})
is attained by a certain multi-geodesic which, after the scaling map is applied, may be viewed as a multi-polymer; the weight specified in~(\ref{e.multiweight}) is the sum of the weights of the zigzags that constitute this multi-polymer. 
To set up notation in this regard, note that
Lemma~\ref{l.severalpolyunique} ensures the almost sure existence of the $k$-tuple 
\begin{equation}\label{e.ptuple}
 \Big( P^{k}_{\big(n t_1 \bar{\bf 1} + 2 n^{2/3} \bar{x}, n t_1 \bar{\bf 1}) \to (n t_2  \bar{\bf 1}  + 2n^{2/3} \bar{y} , n t_2    \bar{\bf 1} \big)     ;i} :  i \in \intint{k} \Big) \in  \upright^k_{\big(n t_1 \bar{\bf 1} + 2 n^{2/3} \bar{x}, n t_1 \bar{\bf 1}) \to (n t_2  \bar{\bf 1}  + 2n^{2/3} \bar{y} , n t_2    \bar{\bf 1} \big)} \, .
\end{equation}
Indeed, the right-hand set is clearly non-empty under our hypotheses, permitting the use of the above lemma.
We define the multi-polymer
$$
\Big( \rho_{n,k,i;(\bar{x},t_1)}^{(\bar{y},t_2)} : i \in \intint{k} \Big)
$$
to be
the $k$-tuple 
of $n$-zigzags whose elements are the respective images under $R_n$
of the elements in the $k$-tuple in~(\ref{e.ptuple}). The $k$-tuple itself~\hfff{multipolymer} will be denoted by  $\rho_{n,k;(\bar{x},t_1)}^{(\bar{y},t_2)}$.
In this way, the $i\textsuperscript{th}$ component\hfff{multipolymercomp} $\rho_{n,k,i;(\bar{x},t_1)}^{(\bar{y},t_2)}$ is an $n$-zigzag with starting and ending points $(x_i, \btone)$ and $(y_i, \bttwo)$. This zigzag is not necessarily a polymer.

Each of the $k$ zigzag components of  $\rho_{n,k;(\bar{x},t_1)}^{(\bar{y},t_2)}$ has a weight via~(\ref{e.weightzigzag}): the  $i\textsuperscript{th}$ weight is 
  $\weight(Z) =  2^{-1/2} n^{-1/3} \big( E(S) - 2 n \tot  - 2n^{2/3}(y_i-x_i) \big)$,
where $S$ is the pre-image staircase in question. 
This $i\textsuperscript{th}$ weight will be denoted by $\weight_{n,k,i;(\bar{x},t_1)}^{(\bar{y},t_2)}$.
The weight of the multi-polymer  $\rho_{n,k;(\bar{x},t_1)}^{(\bar{y},t_2)}$
 equals the component sum $\sum_{i=1}^k \weight_{n,k,i;(\bar{x},t_1)}^{(\bar{y},t_2)}$, and is seen to be equal to~$\weight_{n,k;(\bar{x},t_1)}^{(\bar{y},t_2)}$.

An important special case arises when $\bar{x}$ and $\bar{y}$ are both constant vectors, of the form $\overline{x} = x \bar{\bf 1} \in \R^k$
and  $\overline{y} = y \bar{\bf 1} \in \R^k$. 
The multi-polymer may be called a {\em multi-polymer watermelon} in this case. The cases when merely one of $\bar{x}$ and $\bar{y}$ is a constant vector are also significant: these are the forward and backward bouquets that have been discussed in the road map.

\section{Multi-weight encoding line ensembles: definitions and key properties}\label{s.multiweight}

As the road map has described, a key input drawing in part on integrable probability, Corollary~\ref{c.neargeod.t}, will assert that   $\nearpoly_{n,k;(x,0)}^{(y,1)}(\eta)$ has probability $\eta^{k^2 - 1 + o(1)}$
uniformly in high $n$. This assertion about the weight of a multi-polymer watermelon
% maximum possible weight of a system of zigzags, with shared endpoints but otherwise disjoint, 
falls within the realm of integrable probability because this weight is the  maximum possible attained by a system of~$k$ zigzags, with given shared endpoints but otherwise disjoint, 
and, as such, it may be expressed as the value of the sum of the~$k$ uppermost curves in a line ensemble that is naturally associated to Brownian LPP (by means of the RSK correspondence). In this section, we set up the notation for these ensembles of random curves. 

\subsection{Line ensembles that encode polymer weights}\label{s.encode}

Let $(n,t_1,t_2) \in \N \times \R^2_\leq$ be a compatible triple,  and let $x \in \R$. 
We define the scaled  forward line ensemble~\hfff{scaledforwardle}
$$
\mc{L}_{n;(x,t_1)}^{\uparrow;t_2}: \intint{n \tot + 1} \times \big[ x - 2^{-1} n^{1/3} \tot  , \infty\big) \to \R
$$
rooted at $(x,t_1)$ with duration $\tot$
by declaring that, for each $k \in  \intint{n \tot + 1}$
and $y \geq x - 2^{-1} n^{1/3} \tot$,
\begin{equation}\label{e.scaledweight}
 \sum_{i=1}^k \mc{L}_{n;(x,t_1)}^{\uparrow;t_2}(i,y) \, =  \,
\weight_{n,k;(x \bar{\bf 1},t_1)}^{(y \bar{\bf 1},t_2)} \, .
\end{equation}
The ensemble is called `forward', and the notation is adorned with the symbol $\uparrow$, because it is the spatial location $y$ attached to the more advanced time $t_2$ that is treated as the variable.
We stand at $(x,t_1)$ and look forward in time to $t_2$ to witness behaviour as a function of location $y$. When $k=1$, the lowest indexed ensemble curve $y \to  \mc{L}_{n;(x,t_1)}^{\uparrow;t_2}(1,y)$
records the polymer weight profile from $(x,t_1)$ to $(y,t_2)$. The sum of the $k$ lowest indexed curves records the multi-polymer watermelon weight between these locations.

Equally, we may stand at $(y,t_2)$ and look backward in time to a variable location $x \in \R$ at time $t_1$. Indeed,  fixing $y \in \R$,
we define the scaled backward line ensemble\hfff{scaledbackwardle}
$$
\mc{L}_{n;t_1}^{\downarrow;(y,t_2)}: \intint{n \tot + 1} \times \big( - \infty ,  y + 2^{-1} n^{1/3} \tot \big] \to \R
$$
 rooted at $(y,t_2)$ with duration $\tot$
by declaring that, for each $k \in  \intint{n \tot + 1}$
and $x \leq y + 2^{-1} n^{1/3} \tot$,
$$
\sum_{i=1}^k  \mc{L}_{n;t_1}^{\downarrow;(y,t_2)}(i,x) \, = \, 
\weight_{n,k;(x \bar{\bf 1},t_1)}^{(y \bar{\bf 1},t_2)} \, .
$$
 
It is valuable to hold a vivid picture of the two ensembles. Each is an ordered system of random continuous curves, with the lowest indexed curve uppermost. The ensemble curve of any given index locally resembles Brownian motion but globally follows the contour of the parabola $- 2^{-1/2} (y-x)^2 \tot^{-4/3}$, as a function of $y$ or $x$ in the forward or backward case. Each ensemble adopts values of order $\tot^{1/3}$
when $x$ and $y$ differ by an order of $\tot^{2/3}$. More negative values, dictated by parabolic curvature, are witnessed outside this region. This description holds sway in a region that expands from the origin as the parameter $n$ rises.

\begin{figure}[ht]
\begin{center}
\includegraphics[height=10cm]{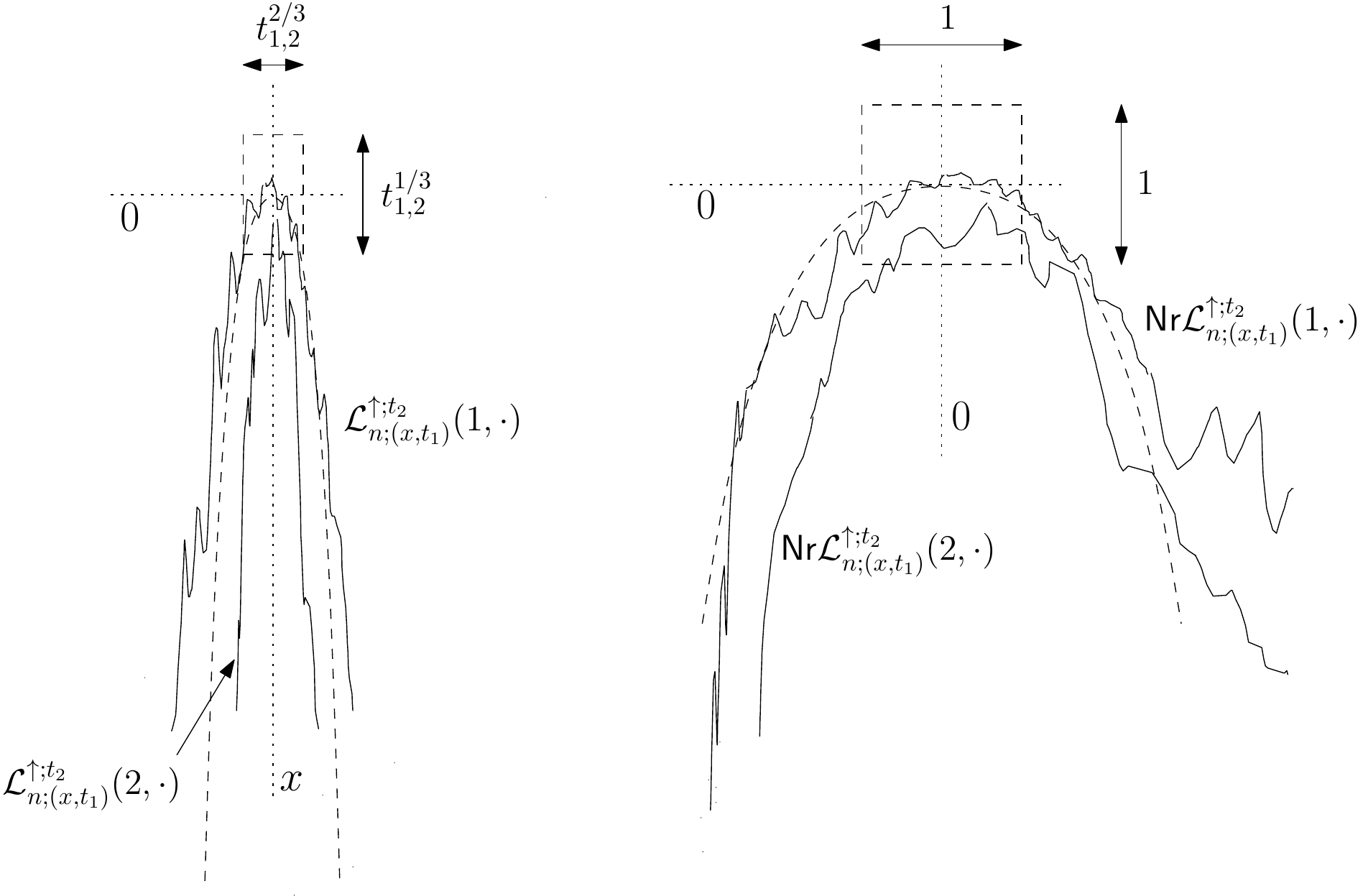}
\caption{Let $(n,t_1,t_2)$ be a compatible triple, with $\tot$ some small fixed value, and $n$ large. For $x \in \R$ given, the two sketches show the highest two curves in the scaled, and normalized, forward line ensembles rooted at $(x,t_1)$ with duration $\tot$.}
\label{f.lineensembles}
\end{center}
\end{figure}

Clearly, then, our ensembles have fundamental differences according to the value of $\tot$: sharply peaked ensemble curves when $\tot$ is small, and much flatter curves when $\tot$ is large.
 A simple further parabolic transformation will serve to put the ensembles on a much more equal footing. 
 Since the ensembles are already scaled objects, we will use the term `normalized' to allude to the newly transformed counterparts. A scaled forward ensemble and its normalized counterpart are depicted in Figure~\ref{f.lineensembles}.

That is, 
we define the {\em normalized} forward ensemble\hfff{normalizedle}  
\begin{equation}\label{e.forward}
 \scaledle_{n;(x,t_1)}^{\uparrow;t_2}: \intint{n \tot + 1} \times \big[- 2^{-1} (n \tot)^{1/3}  , \infty \big) \to \R \, ,
\end{equation}
%with base $(x,t_2)$ and duration $t_2 - t_1 > 0$.
setting 
\begin{equation}\label{e.scaledln}
\scaledle_{n;(x,t_1)}^{\uparrow;t_2}\big( k, z \big)
 = \tot^{-1/3} \mc{L}_{n;(x,t_1)}^{\uparrow;t_2}\big( k, x + \tot
^{2/3} z \big)
\end{equation}
for each $k \in \intint{n \tot + 1}$ and $z \geq - 2^{-1} n^{1/3} \tot^{-2/3}$.
Similarly, the normalized backward ensemble  
\begin{equation}\label{e.backward}
 \scaledle_{n;t_1}^{\downarrow;(y,t_2)}: \intint{n \tot + 1} \times \big(-\infty, 2^{-1} (n \tot)^{1/3}   \big] \to \R \, ,
\end{equation}
%with base $(x,t_2)$ and duration $t_2 - t_1 > 0$.
is specified by setting 
$$
%\begin{equation}\label{e.scaledln.backward}
\scaledle_{n;t_1}^{\downarrow;(y,t_2)}\big( k, z \big)
 = \tot^{-1/3} \mc{L}_{n;t_1}^{\downarrow;(y,t_2)}\big( k, y + \tot
^{2/3} z \big)
$$
for each $k \in \intint{n \tot + 1}$ and $z \leq  2^{-1} n^{1/3} \tot^{-2/3}$.

The curves in the new ensembles locally resemble Brownian motion as before, but they have been centred and squeezed so that now the parabola that dictates their overall shape is $- 2^{-1/2} z^2$.
This picture is accurate in a region that expands as the parameter $n \tot$ rises.

\subsection{Brownian Gibbs line ensembles}

The notion of a Brownian Gibbs line ensemble was introduced in~\cite{AiryLE} to capture a system of ordered curves that arise by conditioning Brownian motions or bridges on mutual avoidance. The precise definition is not logically needed in this article, but we offer an informal summary next.

\subsubsection{An overview}  
Let $n \in \N$ and 
let $I$ be a closed interval in the real line.
A $\intint{n}$-indexed line ensemble defined on $I$ is a random collection of continuous curves  $\mc{L}:\intint{n} \times I \to \R$ specified under a probability measure $\PP$. The $i\textsuperscript{th}$ curve is thus $\mc{L}(i,\cdot): I \to \R$. (The adjective `line' has been applied to these systems perhaps because of their origin in such models as Poissonian LPP, where the counterpart object has piecewise constant curves. We will omit it henceforth.)
An ensemble is called {\em ordered} if $\mc{L}(i,x) > \mc{L}(i+1,x)$ whenever $i \in \intint{n-1}$ and $x$ lies in the interior of $I$.
The curves may thus assume a common value at any finite endpoint of $I$.
We will consider ordered ensembles that satisfy a condition called the Brownian Gibbs property.
Colloquially, we may say that an ordered ensemble is called Brownian Gibbs if it arises from a system of Brownian bridges or Brownian motions defined on $I$ by conditioning on the mutual avoidance of the curves at all times in $I$.
 
\subsubsection{Defining $(c,C)$-regular ensembles} 
We are interested in ensembles that are not merely Brownian Gibbs but that hew to the shape of a parabola and have one-point distributions for the uppermost curve that enjoy tightness properties. We will employ the next definition, which specifies a  $(\bar\phi,\rsc,\rsC)$-regular ensemble from~\cite[Definition~$2.4$]{BrownianReg},
%d.regularsequence  
in the special case where the vector $\bar\phi$  equals  $(1/3,1/9,1/3)$.

\begin{definition}\label{d.regularsequence} 
%Let $n \in \N$ and $\xnmac \in [0,\infty)$. 
%
Consider a Brownian Gibbs ensemble of  the form 
$$
\mc{L}: \intint{\nmac} \times \big[ - \xnmac , \infty \big) \to \R  \,  ,
$$
and which is defined on a probability space under the law~$\PP$.
The number $\nmac = \nmac(\mathcal{L})$ of ensemble curves and the absolute value $\xnmac$ of the finite endpoint may take any values in $\N$ and $[0,\infty)$.

Let \hfff{parabola} $\para:\R \to \R$ denote the parabola $\para(x) = 2^{-1/2} x^2$.
 
Let $\rsC$ and $\rsc$ be two positive constants. The ensemble $\mc{L}$
is said to be $(\rsc,\rsC)$-regular~\hfff{regular} if the following conditions are satisfied.
\begin{enumerate}
\item {\bf Endpoint escape.} $\xnmac \geq  \rsc N^{1/3}$.
\item {\bf One-point lower tail.} If $z \in [ -\xnmac, \infty)$ satisfies $\vert z \vert \leq \rsc \nmac^{1/9}$, then
$$
\PP \Big( \mc{L} \big( 1,z\big) + \para(z) \leq - s \Big) \leq \rsC \exp \big\{ - \rsc s^{3/2} \big\}
$$
for all $s \in \big[1, \nmac^{1/3} \big]$.
\item {\bf One-point upper tail.}  If $z \in [ -\xnmac, \infty)$ satisfies $\vert z \vert \leq \rsc \nmac^{1/9}$, then
$$
\PP \Big( \mc{L} \big( 1,z\big) +  \para(z) \geq  s \Big) \leq \rsC \exp \big\{ - \rsc s^{3/2} \big\}
$$
for all $s \in [1, \infty)$.
\end{enumerate}
A Brownian Gibbs ensemble of the form 
$$
\mc{L}: \intint{\nmac} \times \big( -\infty , \xnmac  \big] \to \R
$$
is also said to be $(\rsc,\rsC)$-regular if the reflected ensemble $\mc{L}( \cdot, - \cdot)$ is. This is equivalent to the above conditions when instances of $[ - \xnmac, \infty)$
are replaced by $(-\infty, \xnmac]$.
\end{definition}

We will refer to these three regular ensemble conditions as $\rmreg(1)$,  $\rmreg(2)$ and $\rmreg(3)$.

\subsubsection{The normalized forward and backward ensembles are $(c,C)$-regular}

Our reason for invoking the theory of regular Brownian Gibbs ensembles is that the normalized Brownian LPP ensembles verify the definition.
This assertion is made by the next result; as we explain shortly, in Section~\ref{s.organization}, 
its proof has in essence been provided in~\cite{BrownianReg} and will formally be given in Appendix~\ref{s.normal}.
%submit discrepancy
\begin{proposition}\label{p.scaledreg}
Let $(n,t_1,t_2) \in \N \times \R^2_<$ be a compatible triple.
\begin{enumerate}
\item  Let $x \in \R$.  The ensemble $\mc{L}$
given by
$$
 \scaledle_{n;(x,t_1)}^{\uparrow;t_2}: \intint{n \tot + 1} \times \big[- 2^{-1} (n \tot)^{1/3}  , \infty \big) \to \R \, ,
$$
is Brownian Gibbs, where $\nmac(\mc{L}) = n \tot  + 1$ and $\xnmac = 2^{-1} (n \tot)^{1/3}$.
\item 
 There exist positive constants $\rsC$ and $\rsc$, which may be chosen independently of all such choices of the parameters $t_1$, $t_2$, $x$ and $n$,
such that the ensemble~$\mc{L}$ is $(\rsc,\rsC)$-regular. 
\item
If, in place of $x$, we consider $y \in \R$, and the ensemble $\mc{L}$
given by
$$
\scaledle_{n;t_1}^{\downarrow;(y,t_2)}(1,\cdot)  : \intint{n \tot + 1} \times  \big( - \infty ,   2^{-1} (n \tot)^{1/3}  \big] \to \R \, ,
$$ 
then the two preceding assertions hold, the second now  independently of the parameters $t_1$, $t_2$, $y$ and $n$.
\end{enumerate}
\end{proposition}
\subsection{Some useful properties of regular ensembles}\label{s.useful}

We are about to state Proposition~\ref{p.mega}, whose four parts assert the various properties, beyond the regular sequence conditions $\rmreg(1)$, $\rmreg(2)$ and $\rmreg(3)$, of the normalized ensembles that we will be employing.

We fix henceforth the values of the two positive constants~$\rsC$ and $\rsc$, specifying them by Proposition~\ref{p.scaledreg}. 
Since bounding the constants would render hypotheses of our results to be explicit, we mention that they are determined in~\cite[Appendix~$A.1$]{BrownianReg}
via Ledoux~\cite[(5.16)]{Ledoux} and Aubrun's~\cite[Proposition~$1$]{Aubrun} bounds on the lower and upper tail of the maximum eigenvalue of a matrix in the Gaussian unitary ensemble. 

We now specify two sequences $\big\{ C_k : k \geq 1 \big\}$ and $\big\{ c_k: k \geq 1 \big\}$, their values expressed in terms of $\rsC$ and $\rsc$. This usage for $C_k$ and $c_k$ is retained throughout the paper. We set, for each $k \geq 2$, 
  \begin{equation}\label{e.formere}
 \formerE_k = \max \Big\{  10 \cdot 20^{k-1} 5^{k/2} \Big( \tfrac{10}{3 - 2^{3/2}} \Big)^{k(k-1)/2} C \, , \, e^{c/2} \Big\} 
 \end{equation}
as well as $C_1 = 140 C$; and
\begin{equation}\label{e.littlec}
 c_k =   \big( (3 - 2^{3/2})^{3/2} 2^{-1} 5^{-3/2} \big)^{k-1} c_1 \, ,
\end{equation}
with  $c_1 = 2^{-5/2} c \wedge 1/8$. 
Note that 
 $\limsup \formerE_k^{1/k^2} < \infty$ and $\liminf c_k^{1/k} > 0$.

Recall from Definition~\ref{d.regularsequence} that $\para:\R \to \R$ denotes the parabola $\para(x) = 2^{-1/2} x^2$.

%The next result is the first three parts of [ModCon, Proposition *].

\begin{proposition}\label{p.mega}
Suppose that  $\mc{L} = \mc{L}_N$, mapping either $\intint{N} \times [-\xnmac,\infty)$ or  $\intint{N} \times (-\infty,\xnmac]$, to $\R$,
is a $(\rsc,\rsC)$-regular ensemble, where $N \in \N$ and $\xnmac \geq 0$.
\begin{enumerate}
\item (Pointwise curve lower bound)
Let $k \in \N$ and $z,s \in \R$.
Suppose that  $N \geq k 
  \vee  (c/3)^{-18} \vee  6^{36}$,
 $\vert z \vert \leq 2^{-1} \rsc N^{1/18}$ and $s \in \big[0, 2 N^{1/18} \big]$. Then
$$
\PP \Big( \, \mc{L}_N\big( k,z\big) + \para(z) \leq - s \, \Big) \leq \formerE_k \exp \big\{ - c_k s^{3/2} \big\} \, .
$$
\item (Uniform curve lower bound) 
For $k \in \N$,
 let  $\Cstrong_k = 20^{k-1} 2^{k(k-1)/2} \Cstrong_1$  where $\Cstrong_1 = 10C$. 
Set 
$$
\rcon = 5 (3 - 2^{3/2})^{-1} \, , \, r_1 = 2^{3/2} \, , 
\,  \textrm{ and  $r_k = \max \{ 5^3 , \rcon r_{k-1} \big\}$ for $k \geq 2$} \, .
$$
 Whenever  $k \in \N$ and $(t,r,y) \in \R$ satisfy $N \geq k 
  \vee  (c/3)^{-18} \vee  6^{36}$, 
  $t \in \big[ 0 , N^{1/18} \big]$, $r \in \big[ r_k \, , \, 2N^{1/18} \big]$
  and $y \in \rsc/2 \cdot [-  N^{1/18},  N^{1/18}]$,
$$
\PP \Big( \inf_{x \in [y-t,y+t]} \big( \mc{L}_N(k,x) + \para(x) \big) \leq - r \Big) \, \leq \, \Big( t \vee  5 \vee (3 - 2^{3/2})^{-1/2} r_{k-1}^{1/2}   \Big)^k \cdot \Cstrong_k \exp \big\{ - c_k r^{3/2} \big\} \, .
$$
\item (No Big Max) 
For  $\vert y \vert \leq 2^{-1} c  N^{1/9}$, $r \in \big[0,4^{-1} \rsc  N^{1/9}\big]$, $t \in \big[ 2^{7/2} , 2 N^{1/3} \big]$ and $N \geq c^{-18}$,
$$
\PP \Big( \sup_{x \in [y-r,y+r]} \big( \mc{L}_N ( 1,x ) + 2^{-1/2}x^2 \big) \geq t \Big) \leq  (r + 1) \cdot  6  \rsC \exp \big\{ - 2^{-11/2} \rsc  t^{3/2} \big\} \, . 
$$
\item (Collapse near infinity) 
For $\eta \in (0,\rsc]$, let
$\ell = \ell_\eta:\R \to \R$ denote the even function which is affine on $[0,\infty)$ and has gradient $ - 5 \cdot 2^{-3/2} \eta \nmac^{1/9}$ on this interval, and which satisfies
$\ell(\eta \nmac^{1/9}) = \big( - 2^{-1/2} + 2^{-5/2} \big) \eta^2 \nmac^{2/9}$. If  $\nmac \geq 2^{45/4} \rsc^{-9}$, then  
\begin{eqnarray*}
 & & \PP \Big( \mc{L}_N \big(1,z\big) > \ell(z) \, \, \textrm{for some} \, \,  z \in  D \setminus \big[ - \eta  \nmac^{1/9} , \eta \nmac^{1/9} \big] \Big) \\
  & \leq &  
6C \exp \Big\{ - c \eta^3  2^{-15/4}   \nmac^{1/3} \Big\} \, .
\end{eqnarray*}
The set $D$ is the spatial domain of $\mc{L}$, either $[ - \xnmac , \infty )$  or  $(-\infty,\xnmac]$.  
\end{enumerate}
\end{proposition}

These four assertions are proved in~\cite{BrownianReg}. Respectively, they appear as  the following results in that article:
Proposition~$2.7$,
%p.othercurves
Proposition~$A.2$,
%p.strongothercurves
Proposition~$2.28$,
%p.nobigmax.gen
and
Proposition~$2.30$.
%p.collapsenearinfinity

A few words about the meaning of the four parts of this proposition: 
the first part asserts a lower bound for a curve in a regular ensemble of any given index.
This holds by definition when the index is $k=1$ via $\rmreg(2)$ but is non-trivial in the other cases. The second part strengthens this conclusion to speak of the minimum value of such a curve on a compact interval. The third similarly strengthens the one-point upper tail $\rmreg(3)$. In regard to the fourth, note that $\rmreg(2)$ and $\rmreg(3)$ do not assert that curves hew to the parabola $-2^{-1/2}z^2$
globally, but only in an expanding region about the origin, of width $2cN^{1/9}$ centred at the origin, where $N$ is the ensemble curve cardinality. Proposition~\ref{p.mega}(4) offers a substitute control on curves far from the origin, showing them to decay at a rapid but nonetheless linear rate in the region beyond scale $N^{1/9}$.
  
These four assertions from~\cite{BrownianReg} are all consequences of the theory of Brownian Gibbs resampling introduced in~\cite{AiryLE} and developed in~\cite{BrownianReg}. They are all rather simple consequences of this theory, with short proofs in~\cite{BrownianReg}. It is in fact in the upcoming Theorem~\ref{t.neargeod} that we cite a result from~\cite{BrownianReg} that harnesses the Brownian Gibbs theory from~\cite{BrownianReg} in a substantial way. 
  
 \section{Some further key inputs}\label{s.input}

\subsection{The scaling principle}\label{s.scalingprinciple}

Fundamental to the theory of the KPZ fixed point is the triple $(1/3,2/3,1)$ of exponents that reflects the scaling laws for weight, polymer geometry, and polymer lifetime. The triple manifests itself in the context of our use of scaled coordinates. It is a simple consequence of the definition of the scaling transformation~$R_n$ and of the energy-weight relationship~(\ref{e.weightzigzag}) that the following useful fact holds true.
 
 \noindent{\em The scaling principle.}\hfff{scaling} 
Let $(n,t_1,t_2) \in \N \times \R^2_<$ be a compatible triple.
 Any statement concerning the system of $n$-zigzags, including weight information, is equivalent to the corresponding statement concerning the system of $n\tot$-zigzags, provided that the following changes are made:
 \begin{itemize}
 \item the index $n$ is replaced by $n\tot$;
 \item any time is multiplied by $\tot^{-1}$;
 \item any weight is multiplied by $\tot^{1/3}$;
 \item and any horizontal distance is multiplied by $\tot^{-2/3}$.
 \end{itemize}
 A little more explanation of the scaling principle appears in~\cite[Section~$2.3$]{ModCon}.

\subsection{Several disjoint near polymers with shared endpoints}

Our road map offering a route to Theorem~\ref{t.disjtpoly.pop} involves specifying and bounding the probability of an event $\nearpoly$. This work has been carried out in~\cite{BrownianReg}, and we here recall the needed result. 

Define 
$$
\neargeod_{n,k;(x,t_1)}^{(y,t_2)}(\eta) = \Big\{ \, \weight_{n,k;(x\bar{\bf 1},t_1)}^{(y\bar{\bf 1},t_2)}    \geq k \cdot \weight_{n;(x,t_1)}^{(y,t_2)} -   \tot^{1/3} \eta \, \Big\} \, .
$$
The presence of the $\tot^{1/3}$ factor that multiplies $\eta$ is consistent with the scaling principle. It means that $\eta$ has the role of a measure of discrepancy from the maximum weight attainable in principle, independently of the value of $\tot$. This $\neargeod$ may also be expressed using a scaled forward 
ensemble: namely,
$$
\neargeod_{n,k;(x,t_1)}^{(y,t_2)}(\eta) \, = \, \bigg\{ \,
\sum_{i=1}^k \mc{L}_{n;(x,t_1)}^{\uparrow;t_2}(i,y)   \geq k \cdot \mc{L}_{n;(x,t_1)}^{\uparrow;t_2}(1,y) \,  -  \,  \tot^{1/3}  \eta \, \bigg\} \, .
$$

The next result is perhaps the most consequential of this article's inputs. It is~\cite[Theorem~$1.12$]{BrownianReg}    up to a relabelling of one parameter.
%Conditions such as $K_1 > 0$ are not stated in BrownianReg.

\begin{theorem}\label{t.neargeod}
There exist constants $K_0 \geq 1$, $K_1 > 0$, $a_0 \in (0,1)$ and $\eta_0 > 0$, and a positive sequence $\{ \beta_k: k \in \N \}$ with $\limsup \beta_k^{1/k} < \infty$,
such that,
for $n,k \in \N$,  $x \in \R$ and $\eta \in \big(0,(\eta_0)^{k^2}\big)$ satisfying $k \geq 2$, $n \geq k \vee \, (K_0)^{k^2} \big( \log \eta^{-1} \big)^{K_0}$ and   $\vert x \vert \leq a_0 n^{1/9}$,
$$
 \eta^{k^2 - 1}  \cdot \exp \big\{ - e^{K_1 k} \big\} \, \leq \,  \PP \Big(  
\neargeod_{n,k;(0,0)}^{(x,1)} \big(   \eta \big) \Big) \, \leq \, \eta^{k^2 - 1} \cdot 
\exp \Big\{ \beta_k \big( \log \eta^{-1} \big)^{5/6} \Big\} \, .
$$
\end{theorem}

The actual consequence that will be used when we implement the road map is the upper bound in the result now stated.

\begin{corollary}\label{c.neargeod.t}
There exist  constants $K_0 \geq 1$, $K_1 > 0$, $a_0 \in (0,1)$ and $\eta_0 > 0$
 and a positive sequence $\{ \beta_k: k \in \N \}$ with $\limsup \beta_k^{1/k} < \infty$
such that,
for $n,k \in \N$, $(t_1,t_2) \in \R^2_<$,  $x,y \in \R$ and $\eta \in \big(0,(\eta_0)^{k^2}\big)$ satisfying $k \geq 2$, $\tot n \geq k \vee \, (K_0)^{k^2} \big( \log \eta^{-1} \big)^{K_0}$ and   $\tot^{-2/3} \vert y - x \vert \leq a_0 n^{1/9}$,
$$
 \eta^{k^2 - 1} \cdot \exp \big\{ - e^{K_1 k} \big\}  \, \leq  \, \PP \Big(  
\neargeod_{n,k;(x,t_1)}^{(y,t_2)} \big(   \eta \big) \Big) \, \leq \, \eta^{k^2 - 1} \cdot 
\exp \Big\{ \beta_k \big( \log \eta^{-1} \big)^{5/6} \Big\} \, .
$$
\end{corollary}
\noindent{\bf Proof.} This follows from Theorem~\ref{t.neargeod} by the scaling principle. \qed

\subsection{A useful tool: tail behaviour for polymer weight suprema and infima}\label{s.usefultool}

Short scale rewiring of zigzags is a central part of the plan in the road map. The next definition and result will be used in order to prove that such rewiring occurs at a manageable cost.

For $x,y \in \R$, $w_1,w_2 \geq 0$ and  $r > 0$, let $\maxmin_{n;([x,x+w_1],t_1)}^{([y,y+w_2],t_2)}( r )$~\hfff{pwr} denote the {\em polymer weight regularity} event that, for all $(u,v) \in [0,w_1] \times [0,w_2]$,
$$
\Big\vert \, \tot^{-1/3} \weight_{n;(x+u,t_1)}^{(y+v,t_2)} +  2^{-1/2} \tot^{-4/3}  \big(y+v-x - u \big)^2 \, \Big\vert \, \leq \, r \, .
$$
When $w_1 = 0$, so that the interval $[x,x+w_1]$ is a singleton, we write $x$ in place of $\{ x \}$ or $[x,x]$ in using this notation.   

The next result, \cite[Corollary~$2.1$]{ModCon}, expresses the one-third power law for polymer weight that was discussed in Section~\ref{s.roadmap}.
\begin{corollary}\label{c.maxminweight}
Let $(n,t_1,t_2) \in \N \times \R^2_<$ be a compatible triple for which 
$n\tot \in \N$ is at least $10^{29} \vee 2(c/3)^{-18}$.
 Let $x,y \in \R$, and let $a,b \in \N$, 
 % be positive, 
 be such that
$\big\vert x - y  \big\vert \tot^{-2/3} + \max\{ a,b\} - 1 \leq   6^{-1}  \rsc  (n\tot)^{1/18}$.
 Let $r \in \big[  34 \, , \, 4 (n \tot)^{1/18} \big]$.
Then
$$
 \PP \Big( \neg \, \maxmin_{n;([x,x+a \tot^{2/3}],t_1)}^{([y,y+b\tot^{2/3}],t_2)}(r) \Big) \leq  
 ab \cdot  400 C \exp \big\{ - c_1 2^{-10} r^{3/2} \big\} \, .
$$
\end{corollary}

\subsection{Local weight regularity}\label{s.polyweightreg}

Recall from the road map that a power of one-half dictates the H\"older continuity of polymer weights as the endpoints are varied horizontally. A rigorous version of this assertion has been presented in~\cite[Theorem~$1.1$]{ModCon}, and we recall it now, again using the notation $Q:\R \to \R$, $Q(u) = 2^{-1/2} u^2$, to denote the parabola that dictates the global shape of the weight profile.

\begin{theorem}\label{t.differenceweight}
Let $n \in \N$ and   $x,y \in \R$  satisfy 
 $n \geq 10^{32} c^{-18}$ and   $\big\vert x - y  \big\vert \leq 2^{-2} 3^{-1} \rsc  n^{1/18}$.
Let $\e \in (0,2^{-4}]$ and 
 $R \in \big[10^4 \, , \,   10^3 n^{1/18} \big]$.
Then
$$
\PP \left( \sup_{\begin{subarray}{c} u_1,u_2 \in [x,x+\e] \\
    v_1,v_2 \in [y,y+\e]  \end{subarray}} \Big\vert \weight_{n;(u_2,0)}^{(v_2,1)} + Q(v_2 - u_2) - \weight_{n;(u_1,0)}^{(v_1,1)} - Q(v_1 - u_1) \Big\vert  \, \geq \, \e^{1/2}
  R  \right)
$$
  is at most  $10032 \, C  \exp \big\{ - c_1 2^{-21}   R^{3/2}   \big\}$.
\end{theorem}

It is also convenient to record a version of this result in which the parabolic curvature term is absent.
Let $I$ and $J$ denote two closed intervals in the real line, each of length $\e$.
Define the {\em local weight regularity} event~\hfff{lwr}  
$$
\lwr_{n;(I,0)}^{(J,1)}\big( \e , r  \big) \, = \,
 \left\{   \sup_{\begin{subarray}{c} x_1,x_2 \in I \\
    y_1, y_2 \in J  \end{subarray}} \Big\vert \weight_{n;(x_2,0)}^{(y_2,1)}  - \weight_{n;(x_1,0)}^{(y_1,1)}  \Big\vert \, \leq \, r \e^{1/2} \right\} \, . 
$$

%From ModCon: The relevant control is offered by Theorem~\ref{t.differenceweight}, except that this theorem addresses parabolically adjusted weight.
%The next result is the one we will apply in proving 
%Theorem~\ref{t.wlp.one}: the main new hypothesis, $\vert x - y \vert \leq \e^{-1/2}$, limits parabolic curvature.  
%The complement of an event $A$ is denoted~$\neg \, A$.

Next is~\cite[Corollary~$6.3$]{ModCon}. The new hypothesis $\vert x - y \vert \leq \e^{-1/2}$ controls parabolic curvature and permits comparison of polymer weights without parabolic adjustment.
\begin{corollary}\label{c.ordweight}
Let 
$n \in \N$ and  $x,y \in \R$ satisfy 
$n \geq 10^{32} c^{-18}$ and   $\big\vert x - y  \big\vert \leq \e^{-1/2} \wedge 2^{-2} 3^{-1} \rsc  n^{1/18}$.
Let  $\e \in (0,2^{-4}]$ and
 $R \in \big[2 \cdot 10^4 \, , \,   10^3 n^{1/18} \big]$.
Then
\begin{equation}\label{e.ordweight}
 \PP \, \Big( \, \neg \, \lwr_{n;([x,x+\e],0)}^{([y,y+\e],1)}\big( \e , R  \big)   \, \Big) \, \leq \, 
 10032 \, C  \exp \big\{ - c_1 2^{-22 - 1/2} R^{3/2}     \big\}
 \, .
\end{equation}
\end{corollary}
%: this result expresses the one-half power law for weight differences.
%\begin{theorem}\label{t.differenceweight}
%Let $\e \in (0,2^{-4}]$. 
%Let $n \in \N$ be an even integer that satisfies
%$n \geq 10^{29} c^{-18}$ and let $x,y \in \R$ satisfy  $\big\vert x - y  \big\vert \leq 2^{-5/3} 3^{-1} \rsc  n^{1/18}$.
%Let 
% $R \in \big[10^6 \, , \,   10^4 n^{1/18} \big]$.
%Then
%$$
% \PP \Big( \neg \, \lwr_{n;([x,x+\e],0)}^{([y,y+\e],1)}\big( \e , R  \big)    \Big) \leq 
% 25064 \, C  \exp \big\{ - c_1 2^{-31} R^{3/2}     \big\}
% \, .
%$$
%\end{theorem}

\subsection{Polymer ordering}\label{s.polyorder}

An important challenge that lies ahead is bouquet construction, and a key difficulty here will be to ensure the disjointness, except at the shared endpoint, of the concerned zigzags. Some natural monotonicity properties of multi-polymers will be needed in the proof. Lemma~\ref{l.tworelations} is our result in this regard. It will be proved in Appendix~\ref{s.mgo}.
%submit discrepancy

Let $(n,t_1,t_2)$ be a compatible triple.
We introduce two ordering relations, $\prec$ and $\preceq$, on the space of $n$-zigzags with lifetime $[t_1,t_2]$. To define the relations, let $(x_1,x_2),(y_1,y_2) \in \R^2$ and  consider
a zigzag  $Z_1$ from $(x_1,t_1)$ to $(y_1,t_2)$ and another
 $Z_2$ from $(x_2,t_1)$ to $(y_2,t_2)$.

\noindent{\em The $\prec$ relation.}
 We declare that $Z_1 \prec Z_2$ if $x_1 \leq y_1$, $x_2 \leq y_2$,
 and the two polymers are horizontally separate.

\noindent{\em The $\preceq$ relation.} Consider again $Z_1$ and $Z_2$.
We declare that $Z_1 \preceq Z_2$ if `$Z_2$ lies on or to the right of $Z_1$': formally, if $Z_2$ is contained in the union of the closed horizontal planar line segments whose left endpoints lie in $Z_1$.

%The next result, [Guide, Lemma], provides two key pieces of information about these relations.

\begin{lemma}\label{l.tworelations}
Let $(n,t_1,t_2)$ be a compatible triple.
\begin{enumerate}
\item Let $(x_1,x_2),(y_1,y_2) \in \R_\leq^2$ and $k \in \N$. Then the multi-polymer components satisfy 
$$
 \PP \Big( \, \rho_{n,k,i;(x_1 \bar{\bf 1)},t_1}^{(y_1 \bar{\bf 1},t_2)} \preceq  \rho_{n,k,i;(x_2 \bar{\bf 1)},t_1}^{(y_2 \bar{\bf 1},t_2)} \, \, \forall \, \, i \in \intint{k} \, \Big) \, = \, 1 \, .
$$
\item Let $Z_1$, $Z_2$ and $Z_3$  be $n$-zigzags of lifetime~$[t_1,t_2]$
that verify $Z_1 \prec Z_2$ and $Z_2 \preceq Z_3$. Then $Z_1 \prec Z_3$.
\end{enumerate}
\end{lemma}

A rather simple sandwiching fact about polymers will also be needed.
\begin{lemma}\label{l.sandwich}
Let $(n,t_1,t_2)$ be a compatible triple, and let  $(x_1,x_2),(y_1,y_2) \in \R_\leq^2$. Suppose that there is a unique $n$-polymer from $(x_i,t_1)$ to $(y_i,t_2)$, both when $i=1$ and $i=2$. (This circumstance occurs almost surely, and the resulting polymers have been labelled $\rho_{n;(x_1,t_1)}^{(y_1,t_2)}$ and  $\rho_{n;(x_2,t_1)}^{(y_2,t_2)}$.)
Now let $\rho$ denote any $n$-polymer that begins in $[x_1,x_2] \times \{ t_1\}$
and ends in  $[y_1,y_2] \times \{ t_2 \}$. Then $\rho_{n;(x_1,t_1)}^{(y_1,t_2)} \preceq \rho \preceq \rho_{n;(x_2,t_1)}^{(y_2,t_2)}$.
\end{lemma}
\noindent{\bf Proof.}
% of Lemma~\ref{l.sandwich}
When two $n$-zigzags $Z_1$ and $Z_2$ share their starting and ending heights, it is easy enough to define associated minimum and maximum zigzags $Z_1 \wedge Z_2$
and $Z_1 \vee Z_2$. The minimum begins and ends at the leftmost of the starting and ending points of $Z_1$ and $Z_2$; during its lifetime, it follows one or other of the two zigzags, always keeping as far to the left as possible.
It is a similar story for the maximum. A formal definition is made for staircases (but this hardly changes for zigzags) in the opening paragraphs of Appendix~\ref{s.mgo}.
%[Guide, Section~{s.mgo}].
 The polymer $\rho$ begins and ends its journey to the right of   $\rho_{n;(x_1,t_1)}^{(y_1,t_2)}$. If the condition  $\rho_{n;(x_1,t_1)}^{(y_1,t_2)} \preceq \rho$ is to be violated, $\rho$ must pass a sojourn in its lifetime to the left of the other polymer. This would make the minimum $\rho \wedge \rho_{n;(x_1,t_1)}^{(y_1,t_2)}$ distinct from $\rho_{n;(x_1,t_1)}^{(y_1,t_2)}$. A moment's thought shows, however, that this minimum is itself a polymer between $(x_1,t_1)$ and $(y_1,t_2)$. This is a contradiction to polymer uniqueness for these endpoints, which is known as a direct consequence of Lemma~\ref{l.severalpolyunique} with~$\ell = 1$.
Thus, $\rho_{n;(x_1,t_1)}^{(y_1,t_2)} \preceq \rho$. The second claimed ordering has a similar derivation. \qed

\subsection{Operations on polymers: splitting and concatenation}\label{s.split}
Two very natural operations are now discussed.

A polymer may be split into two pieces. 
Let  $(n,t_1,t_2) \in \N \times \R^2_\leq$ is a compatible triple, and let $(x,y) \in \R^2$
satisfy  $y \geq x - 2^{-1} n^{1/3} \tot$. Let $t \in (t_1,t_2) \cap n^{-1} \Z$.
Suppose that the almost sure event that
 $\rho_{n;(x,t_1)}^{(y,t_2)}$ is well defined occurs.  Select any element $(z,t) \in \rho_{n;(x,t_1)}^{(y,t_2)}$.
The removal of $(z,t)$ from $\rho_{n;(x,t_1)}^{(y,t_2)}$ creates two connected components. 
Taking the closure of each of these amounts to adding the point $(z,t)$ to each of them. The resulting sets are $n$-zigzags from $(x,t_1)$ to $(z,t)$, and from $(z,t)$ to $(y,t_2)$; indeed,
it is straightforward to see that these are the unique $n$-polymers given their endpoints. We use a concatenation notation~$\circ$ to represent this splitting. In summary,   $\rho_{n;(x,t_1)}^{(y,t_2)} = \rho_{n;(x,t_1)}^{(z,t)} \circ \rho_{n;(z,t)}^{(y,t_2)}$. Naturally, we also have  $\weight_{n;(x,t_1)}^{(y,t_2)} = \weight_{n;(x,t_1)}^{(z,t)} + \weight_{n;(z,t)}^{(y,t_2)}$.
Indeed, the concatenation operation may be applied to any two $n$-zigzags
for which the ending point of the first equals the starting point of the second. Since the zigzags are subsets of $\R^2$, it is nothing other than the operation of union. The weight is additive under the operation.

\subsection{A closure property for a space of several disjoint polymers}\label{s.closure}

A rather technical point about stability under closure of systems of disjoint polymers is now addressed.

Let $(n,t_1,t_2)$ be a compatible triple, 
%Let $x,y \in \R$. An $n$-zigzag $\rho$ with starting and ending points $(x,t_1)$ and $(y,t_2)$ is called {\em proper}
%if $\rho$'s intersections with the vertical coordinates $t_1$ and $t_2$
%is limited to the singleton sets $\{ x \}$ and $\{ y \}$. That is, $\rho$
%begins and ends life sloping, rather than in horizontal motion. 
let $I$ and $J$ be two closed real intervals, and 
let $k \in \N$. We now define two random subsets of $I^k \times J^k \subset \R^{2k}$.  Set $\disjtindex_{n,k;(I,t_1)}^{(J,t_2)}$ to be the collection of vectors
 $\big( x_1,x_2,\cdots,x_k,y_1,y_2,\cdots,y_k \big) \in I^k \times J^k$
 such that there exists a collection of $k$ pairwise disjoint  $n$-polymers which consecutively move from $(x_i,t_1)$ to $(y_i,t_2)$.
 This set is non-empty precisely when $\maxpoly_{n;(I,t_1)}^{(J,t_2)} \geq k$.
 We further define 
 $\horsepindex_{n,k;(I,t_1)}^{(J,t_2)}$, by replacing the condition of pairwise disjoint by that of pairwise horizontally separate in the above.
 
% Note that 
% $\disjtindex_{n,k;(I,t_1)}^{(J,t_2)}$ is non-empty when  $\maxpoly_{n;(I,t_1)}^{(J,t_2)} \geq k$ occurs, because any polymer may be made proper by removing its initial and final horizontal segments. We also clearly have that 
% $\disjtindex_{n,k;(I,t_1)}^{(J,t_2)} \subseteq 
% \horsepindex_{n,k;(I,t_1)}^{(J,t_2)}$. A stronger property holds.
 
 \begin{lemma}\label{l.disjthorsep}
The closure of $\disjtindex_{n,k;(I,t_1)}^{(J,t_2)}$, when  
viewed as a subset of $\R^{2k}$,  is contained in $\horsepindex_{n,k;(I,t_1)}^{(J,t_2)}$.
 \end{lemma}
 \noindent{\bf Proof.} Let $\big( \bar{x}^i, \bar{y}^i \big)$, $i \in \N$,
 be elements of $\disjtindex_{n,k;(I,t_1)}^{(J,t_2)}$ that converge 
 to $\big(\bar{x}, \bar{y} \big) \in \R^k \times \R^k$. Associated to index~$i$
 is a system of disjoint polymers from $\big(\bar{x},t_1\big)$ to $\big(\bar{y},t_2\big)$. It is a simple matter to extract a subsequence of these indices such that all the endpoint locations of the horizontal intervals in each of the $k$ polymers converge pointwise. By considering the $k$ $n$-zigzags with these locations given by the limiting values, we may note that we are dealing with a collection of polymers, due to the continuity of the underlying Brownian motions $B(\ell,\cdot)$, $\ell \in \Z$, in Brownian LPP. Moreover, these $k$ $n$-zigzags are pairwise horizontally separate: indeed, we are equipping the space of $n$-zigzags of lifetime $[t_1,t_2]$ with the Hausdorff topology, so that it is enough to note that  the property of a $k$-tuple of zigzags having a  pair that is not  horizontally separate is an open condition for the $k$-wise product topology on such tuples of zigzags. \qed

\noindent{\em Remark.} Theorems~\ref{t.disjtpoly.pop} and~\ref{t.maxpoly.pop}
make assertions about the random variable  $\maxpoly_{n;(I,t_1)}^{(J,t_2)}$.
The value of this random variable cannot decrease if in its definition in Subsection~\ref{s.maxpoly} we replace the italicized {\em pairwise disjoint} by {\em pairwise horizontally separate}. If we were to redefine the random variable with this change, then the new versions of the results would imply the old ones. In fact, all proofs in this article are valid for the altered definition, and we take the liberty of adopting it henceforth.

\subsection{Organization of the remainder of the paper}\label{s.organization}

Our principal results are the disjoint polymer estimates Theorems~\ref{t.disjtpoly.pop} and~\ref{t.maxpoly.pop}, and the polymer fluctuation bound Theorem~\ref{t.polyfluc}.
A road map has been outlined for the proof of Theorem~\ref{t.disjtpoly.pop}, and, in order to communicate this conceptually central idea with suitable emphasis, the next section, Section~\ref{s.rarity}, is devoted to rigorously formulating it. We mention here that, in fact, this proof will invoke Theorem~\ref{t.polyfluc}: this aspect of the road map has yet to be explained, but the resolution of the third challenge in the road map will involve polymer fluctuation bounds.
Theorem~\ref{t.maxpoly.pop} is a rather straightforward consequence of Theorem~\ref{t.disjtpoly.pop}, and its proof appears at the end of Section~\ref{s.rarity}.

The proof of Theorem~\ref{t.polyfluc} has a rather different flavour. We present the proof next, in Section~\ref{s.polyfluc}. An outline of the proof is offered at the beginning of the section.

The paper ends with four appendices. Appendix~\ref{s.glossary} recalls in a list some of the article's principal notation. Appendices~\ref{s.mgo} and~\ref{s.normal} 
%Sections~\ref{s.mgo} and~\ref{s.normal} 
present proofs for two remaining tools that we have described.

The multi-polymer ordering Lemma~\ref{l.tworelations} is proved in Appendix~\ref{s.mgo}. 
In our presentation of results, and in the overall perspective of this article, we have been eager to focus attention on scaled coordinates. But of course  Lemma~\ref{l.tworelations} is really simply the scaled counterpart to a result that concerns multi-geodesics. This Lemma~\ref{l.severalpolyorder} 
%submit discrepancy
appears late in the paper, in accordance with our wish to encourage the reader to focus on the scaled coordinate interpretation of the subject. It is, however, a basic result that may have an independent interest.

Proposition~\ref{p.scaledreg}, which asserts that our normalized ensembles~(\ref{e.forward}) and~(\ref{e.backward}) are regular, has in essence been proved already in~\cite{BrownianReg}, on the basis of a fundamental result of O'Connell and Yor
describing a counterpart ensemble of curves in unscaled coordinates. Appendix~\ref{s.normal} explains this connection and provides the proof of this proposition, which involves a mundane argument to reconcile a slight difference in ensemble notation between~\cite{BrownianReg} and the present work.

We close this section by discussing two conventions that will govern the presentation of upcoming proofs. %arxiv/submit The second of these comments also explains what is contained in Appendix~\ref{s.calcder}.

\subsubsection{Boldface notation for quoted results}

During the upcoming proofs, we will naturally be making use of the various tools that we have recalled: the $\rmreg$ conditions and the regular ensemble properties Proposition~\ref{p.mega},  
and results presented earlier in this section. 
The statements of such results involve several parameters, in several cases including $(n,t_1,t_2)$, spatial locations $x$ and $y$, and positive real parameters such as $r$.
We will employ a  device that will permit us to disregard notational conflict between the use of such parameters in the contexts of the ongoing proof in question and  the statements of quoted results. When specifying the parameter settings of a particular application, we will allude to the parameters of the result being applied in boldface, and thus permit the reuse of the concerned symbols.

\subsubsection{The role of hypotheses invoked during proofs}\label{s.rolehyp}

When we quote results in order to apply them, we will take care, in addition to specifying the parameters according to the just described convention, to indicate explicitly what the conditions on these parameters are that will permit the quoted result in question to be applied. Of course, it is necessary that the hypotheses of the result being proved imply all such conditions. 
%arXiv/submit The task of verifying that the hypotheses of a given result are adequate for the purpose of obtaining all conditions needed to invoke the various results used during its proof 
%may be called the calculational derivation of the result in question. This derivation
%is necessary, but also in some cases lengthy and unenlightening: a succession of trivial steps. We have chosen to separate the calculational derivations of most of our results from their proofs. 
%These derivations may be found in Appendix~\ref{s.calcder} at the end of the paper. This appendix also contains the proofs of four lemmas that are invoked in the proof of Theorem~\ref{t.polyfluc}.
%These proofs are exclusively calculational matters.
 The task of verifying that the hypotheses of a given result are adequate for the purpose of obtaining all conditions needed to invoke the various results used during its proof 
may be called the calculational derivation of the result in question. This derivation
is necessary, but also in some cases lengthy and unenlightening: a succession of trivial steps. We have chosen to separate the calculational derivations of most of our results from their proofs. 
These derivations may be found in Appendix~$D$
%\ref{s.calcder} 
at the end of the version~\url{http://math.berkeley.edu/~alanmh/papers/NonIntPolymer.pdf} of the paper on the author's webpage; 
%arxiv not submit
the latex source code for this version is an ancillary file to the present arXiv submission.
%~\url{https://arxiv.org/abs/1709.04110}
 Appendix~$D$ also contains the proofs of four lemmas that are invoked in the proof of Theorem~\ref{t.polyfluc}.
These proofs are exclusively calculational matters.

\section{Rarity of many disjoint polymers: proofs of Theorems~\ref{t.disjtpoly.pop} and~\ref{t.maxpoly.pop}}\label{s.rarity}

Theorems~\ref{t.disjtpoly.pop} and~\ref{t.maxpoly.pop}
are slightly simplified versions of the next two results, in which hypotheses on parameters are stated in a more explicit form.
The results make reference to
$\eta_0$, $K_0$ and $a_0$. These are positive constants that are fixed by Theorem~\ref{t.neargeod}.
 The sequence of positive constants $\big\{ \beta_i: i \in \N \big\}$, which verifies $\limsup_{i \in \N} \beta_i^{1/i} < \infty$, is also supplied by this theorem.
 \begin{theorem}\label{t.disjtpoly}
Let
$(n,t_1,t_2) \in \N \times \R^2_<$ be a compatible triple.
Let $\kay \in \N$,  $\e > 0$ and $x,y \in \R$ satisfy the conditions that $\kay \geq 2$,
\begin{equation}\label{e.epsilonbound}
 \e \leq \min \Big\{ \, (\eta_0)^{4\kay^2} \, , \, 10^{-616}  c_\kay^{22}  \kay^{-115} \,  , \, \exp \big\{ - C^{3/8} \big\} \, \Big\} \, ,
\end{equation}
\begin{equation}\label{e.nlowerbound}
    n \tot \geq \max \bigg\{ \, 2(K_0)^{\kay^2} \big( \log \e^{-1} \big)^{K_0}
    \, , \, 10^{606}   c_\kay^{-48}  \kay^{240}  \rsc^{-36} \e^{-222}  \max \big\{   1  \, , \,   \vert x - y \vert^{36} \tot^{-24}  \big\} \, , \, 
  a_0^{-9} \vert y - x \vert^9 \tot^{-6}  \,  \Bigg\} \, ,
\end{equation}
as well as $\vert y - x \vert \tot^{-2/3} \leq \e^{-1/2} \big( \log \e^{-1} \big)^{-2/3} \cdot 10^{-8} c_\kay^{2/3} \kay^{-10/3}$.
Then
\begin{eqnarray*}
 & & \PP \bigg( \maxpoly_{n;\big([x-\tot^{2/3}\e,x+\tot^{2/3}\e],t_1\big)}^{\big([y-\tot^{2/3}\e,y+\tot^{2/3}\e],t_2\big)}  \geq \kay \bigg) \\
 & \leq & 
  \e^{(\kay^2 - 1)/2}  \cdot   10^{32\kay^2}  \kay^{15\kay^2}  c_\kay^{-3\kay^2} C_\kay   \big( \log \e^{-1} \big)^{4\kay^2}   \exp \big\{ \beta_\kay \big( \log \e^{-1} \big)^{5/6} \big\} \, .
\end{eqnarray*}
\end{theorem}

\begin{theorem}\label{t.maxpoly}
There exists a positive constant $\emm_0$, and  a sequence of positive constants $\big\{ \conseqmac_i: i \in \N \big\}$ for which $\sup_{i \in \N} \conseqmac_i \exp \big\{ - 2 (\log i)^{11/12} \big\}$ is finite, such that the following holds.
Let
$(n,t_1,t_2) \in \N \times \R^2_<$ be a compatible triple. Let $x,y \in \R$, $a,b \in \N$ and $\emm \in \N$. Write $h = a \vee b$. Suppose that 
 $$
 \emm \geq   \emm_0 \vee \big( \vert x - y \vert \tot^{-2/3} + 2h \big)^3  
 $$ 
 and
\begin{eqnarray}
 n \tot  & \geq  & \max \bigg\{ \,    2(K_0)^{(12)^{-2} (\log \log \emm)^2} \big( \log \emm \big)^{K_0}
    \, , \,  a_0^{-9} \big( \vert y - x \vert  \tot^{-2/3} + 2h \big)^9 \, ,
    \label{e.nmaxpoly} \\
  & & \qquad \qquad    
     10^{325}       \rsc^{-36} \emm^{465}  \max \big\{   1  \, , \,   \big(\vert y - x \vert \tot^{-2/3} + 2h \big)^{36}   \big\}  
   \,  \bigg\} \, . \nonumber
\end{eqnarray}

Then
$$
 \PP \bigg(  \maxpoly_{n;\big([x,x+a\tot^{2/3}],t_1\big)}^{\big([y,y+b \tot^{2/3}],t_2\big)} \geq \emm \bigg)  \, \leq \,  \emm^{-   (145)^{-1} ( \log \beta)^{-2} (0 \vee \log \log \emm)^2}  \cdot h^{(\log \beta)^{-2} (0 \vee \log \log \emm)^2/{288}  + 3/2} \conseqmac_\emm  \, .
$$
Here, $\beta$ is specified to be $e \vee \limsup_{i \in \N} \beta_i^{1/i}$.
\end{theorem}

\noindent{\bf Proof of Theorems~\ref{t.disjtpoly.pop} and~\ref{t.maxpoly.pop}.}
These results are direct consequences of their just stated counterparts. \qed

In this section, then, our job is to prove Theorems~\ref{t.disjtpoly} and~\ref{t.maxpoly}. Our principal task is to prove Theorem~\ref{t.disjtpoly}.
This we do first, by implementing the road map rigorously. Theorem~\ref{t.maxpoly}
then emerges, at the end of Section~\ref{s.rarity}, as a fairly straightforward consequence: after all, many disjoint polymers running between unit intervals entail several disjoint polymers running between short intervals. 

\subsection{Bouquet construction}\label{s.bouquet}

As we noted at the end of Section~\ref{s.roadmap}, there are three main challenges involved in implementing the road map, beyond the $\nearpoly$ probability estimate Corollary~\ref{c.neargeod.t}.  The technical input concerning the first of the three challenges, polymer weight similarity, 
has been cited already, as Corollary~\ref{c.maxminweight}.

We now address the second challenge: the construction of forward and backward bouquets of suitable weight, dictated by the exponent of one-third. 
 
In Section~\ref{s.zigzagmax}, we introduced the maximum weight 
$\weight_{n,k;(\bar{x},t_1)}^{(\bar{y},t_2)}$ of a system of $k$ pairwise horizontally separate zigzags that make the journey from $(x_i,t_1)$ to $(y_i,t_2)$ for $i \in \intint{k}$. In forward bouquet construction, we specialise to $\bar{x} = x \bar{\bf 1}$ for some $x \in \R$; in backward bouquet construction, we instead choose  
$\bar{y} = y \bar{\bf 1}$ for a given $y \in \R$. 

The aim of bouquet construction will be achieved in the guise of Corollary~\ref{c.bouquetreg}, which shows that the cumulative weight of the maximum weight bouquet, when measured in the suitable one-third power units, is tight, uniformly in the scaling parameter $n$. This result is a direct consequence of Proposition~\ref{p.sumweight}, which asserts an upper and a lower bound on the maximum bouquet weight. The first of these bounds is easy, because polymer weights offer a control from above on bouquet weights. The lower bound is more delicate: horizontal separateness of zigzags must be ensured without damaging weight significantly. Our technique will be a diagonal argument, with the ordering Lemma~\ref{l.tworelations} playing an important role in establishing the necessary form of disjointness.

\begin{proposition}\label{p.sumweight}
Let $(n,t_1,t_2) \in \N \times \R^2_<$ be a compatible triple. Let  $k \in \N$,  $\bar{u} \in \R^k_{\leq}$ and $r > 0$.
\begin{enumerate}
\item Suppose that 
$n \tot$ is at least  $\tot^{-6} \big\vert u_i - x \big\vert^9  \rsc^{-9}$ for $i \in \intint{k}$, and that $r \geq k$. Then
$$
\PP \bigg( \, \tot^{-1/3} \cdot \weight_{n,k;(x \bar{\bf 1},t_1)}^{(\bar{u},t_2)} \geq - \, 2^{-1/2} \sum_{i=1}^k  \tot^{-4/3}  \big( u_i - x \big)^2  \, \, + \, r \, \bigg) \, \leq \,  k \cdot C \exp \Big\{ - c  k^{-3/2} r^{3/2} \Big\} \, .
$$
\item Suppose that
$n \tot$ is at least $\max \big\{  k  \, , \,  3^{18} c^{-18}  \, ,  \, 6^{36}  \big\}$ and is also bounded below for each $i \in \intint{k}$ by $2^{18} \rsc^{-18} \tot^{-12} \big\vert u_i - x \big\vert^{18}$. 
Suppose further that 
$r \in \big[  4k^2 ,  9 k^2 (n\tot)^{1/3} \big]$.
Then
$$
\PP \bigg( \tot^{-1/3} \cdot \weight_{n,k;(x \bar{\bf 1},t_1)}^{(\bar{u},t_2)} \leq - \, 2^{-1/2} \sum_{i=1}^k  \tot^{-4/3}  \big( u_i - x \big)^2  \, \, - \, r  \bigg) \, \leq 
 \,  
3 k^2 
C_k \exp \Big\{ - 2^{-3}  k^{-3} c_k r^{3/2}     \Big\}   \, . 
$$
\end{enumerate}
\end{proposition}

As we prepare to state our bouquet construction tool, Corollary~\ref{c.bouquetreg},
we mention a further aspect of this construction. It would be natural to suppose, on the basis of the ideas presented in Section~\ref{s.roadmap}, that, on the event $\maxpoly_{(n;[x-\e,x+\e],0)}^{([y-\e,y+\e],1)} \geq k$, if the $k$ disjoint polymers that move between $[x-\e,x+\e] \times \{ 0 \}$ and   $[y-\e,y+\e] \times \{ 1 \}$ have endpoints $(u_i,0)$ and $(v_i,1)$,  for $i \in \intint{k}$, then the forward and backward bouquets will be chosen to be the multi-polymers whose weights are $\weight_{n,k;(x \bar{\bf 1},-\e^{3/2})}^{(\bar{u},0)}$ and $\weight_{n,k;(\bar{v},1)}^{(y \bar{\bf 1},1 + \e^{3/2})}$. A microscopic detail must be addressed, however. The first of these multi-polymers is measurable with respect to the randomness in the region $\R \times (-\infty,0]$ and the second to the randomness in  $\R \times [1,\infty)$. Our purpose in the upcoming surgery will be better served were a choice of the bouquet pair to be made for which this assertion is valid with the two {\em open} regions  $\R \times (-\infty,0)$ and  $\R \times (1,\infty)$ instead. 
What is needed is a modified definition in which no use is made of the randomness in the horizontal line with the height of the fixed endpoint: height zero or height one for the forward and backward bouquets.

In the modification, the weight function $\weight$ will be replaced by a proper weight function $\properweight$, in a sense we now specify.
Recall that if $(n,t_1,t_2) \in \N \times \R^2_<$ is a compatible triple and $Z$ is an $n$-zigzag with starting point $(x,t_1)$ and ending point $(y,t_2)$, then $Z$ begins with a planar line segment abutting $(x,t_1)$, which is either horizontal or sloping, and likewise ends with such a segment abutting $(y,t_2)$. We will call $Z$ {\em backward proper} if the line segment that abuts $(x,t_1)$ is sloping, and {\em forward proper} if the line segments that abuts $(y,t_2)$ is sloping.   

We now define the {\em proper} weight $\properweight_{n,k;(x \bar{\bf 1},t_1)}^{(\bar{u},t_2)}$  to equal the maximum weight associated to~$k$ horizontally separate {\em forward proper} zigzags moving consecutively between $(x,t_1)$ to $(u_i,t_2)$ for $i \in \intint{k}$. 
The $k$-tuple maximizer will be denoted by 
$\rho_{n,k;(x\bar{\bf 1},t_1)}^{{\rm prop};(\bar{u},t_2)}$. (Almost sure uniqueness of the maximizer follows in essence from Lemma~\ref{l.severalpolyunique} but is not needed for our purpose.) 
The proper weight definition modifies $\weight_{n,k;(x \bar{\bf 1},t_1)}^{(\bar{u},t_2)}$ by introducing the insistence that the concerned zigzags be forward proper.
Similarly, we define $\properweight_{n,k;(\bar{v},t_1)}^{(y \bar{\bf 1},t_2)}$ to equal   the maximum weight associated to~$k$ horizontally separate {\em backward proper} zigzags moving consecutively between $(v_i,t_1)$ to $(y,t_2)$ for $i \in \intint{k}$.
The $k$-tuple maximizer is
$\rho_{n,k;(\bar{v},t_1)}^{{\rm prop};(y \bar{\bf 1},t_2)}$.

Define the {\em forward bouquet regularity} event~\hfff{fbr} $\fbr_{n,k;(x,t_1)}^{(\bar{u},t_2)}(  r )$ to equal 
$$
 \bigg\{ \, \bigg\vert  \, (\tot - n^{-1})^{-1/3} \cdot \properweight_{n,k;(x \bar{\bf 1},t_1)}^{(\bar{u},t_2)}  \, + \, 2^{-1/2} \sum_{i=1}^k  (\tot - n^{-1})^{-4/3}  \big( u_i +  2^{-1} n^{-2/3}  - x  \big)^2  \, \bigg\vert \, \leq \, r \bigg\}
$$
and the  {\em backward bouquet regularity} event~\hfff{bbr} $\bbr_{n,k;(\bar{v},t_1)}^{(y,t_2)}(  r )$ to equal 
$$
 \bigg\{ \, \bigg\vert  \, (\tot - n^{-1})^{-1/3} \cdot \properweight_{n,k;(\bar{v},t_1)}^{(y \bar{\bf 1},t_2)}  \, + \, 2^{-1/2} \sum_{i=1}^k   (\tot - n^{-1})^{-4/3}  \big( v_i - 2^{-1} n^{-2/3} - y \big)^2 \, \bigg\vert \, \leq \, r \bigg\} \, .
$$
\begin{corollary}\label{c.bouquetreg}
Let $(n,t_1,t_2) \in \N \times \R^2_<$ be a compatible triple. 
Let $k \in \N$ satisfy $k \leq n$, and let $x,y \in \R$, $\bar{u},\bar{v} \in \R^k_{\leq}$ and $r > 0$.
Suppose that
\begin{eqnarray*}
n \tot - 1 & \geq & \max \Big\{  k  \, , \,  3^{18} c^{-18}  \, ,  \, 6^{36} \, , \, (\tot - n^{-1})^{-6} \big( \vert u_i + 2^{-1} n^{-2/3}  - x \vert \vee  \vert v_i - 2^{-1} n^{-2/3}   - y \vert \big)^9  \rsc^{-9} \, ,   \\
 & & \qquad \qquad \qquad \qquad   2^{18}  \rsc^{-18} (\tot - n^{-1})^{-12}  \big( \vert u_i + 2^{-1} n^{-2/3} - x \vert \vee  \vert v_i -  2^{-1} n^{-2/3} - y \vert \big)^{18} \Big\}
\end{eqnarray*}
where here $i$ varies over $\intint{k}$.
Suppose also that $r \in \big[ 4k^2  ,  9 k^2 (\tot - n^{-1})^{1/3} n^{1/3} \big]$. 
Then
$$
\PP \Big( \neg \,  \fbr_{n,k;(x,t_1)}^{(\bar{u},t_2)} \, \Big) \vee \PP \Big(  \, \neg \, \bbr_{n,k;(\bar{v},t_1)}^{(y,t_2)} \Big) \leq   
4 k^2 
C_k \exp \Big\{ - 2^{-3}  k^{-3} c_k r^{3/2}     \Big\}   \, .
$$
\end{corollary}
The form of the corollary is a little cluttered, not least because the microscopic detail we have introduced forces the bouquets to have lifetime $\tot - n^{-1}$ rather than $\tot$.
Setting $n = \infty$ may serve to focus attention on what is essential here, even if it obscures this detail.

\noindent{\bf Proof of Corollary~\ref{c.bouquetreg}.} 
We begin by claiming that 
\begin{equation}\label{e.properweightforward}
\properweight_{n,k;(x \bar{\bf 1},t_1)}^{(\bar{u},t_2)} = \weight_{n,k;(x \bar{\bf 1},t_1)}^{(\bar{u} + 2^{-1}n^{-2/3}{\bf 1} , t_2 - n^{-1})} \, .
\end{equation}
To understand why this is true,  we consider the multi-polymer that realizes the left-hand side and apply the inverse scaling map $R_n^{-1}$, so that it is depicted in unscaled coordinates.
The result is a  $k$-tuple of staircases, each element of which ends in a vertical unit displacement. Omitting these displacements, we obtain a $k$-tuple of staircases which is a multi-geodesic: the sum of the energies of the staircases is maximal given the set of endpoints. Applying the scaling map $R_n$ to reconsider this multi-geodesic in scaled coordinates, we see that its weight is the above right-hand side. 
Similar considerations yield the formula
\begin{equation}\label{e.properweightbackward}
\properweight_{n,k;(\bar{v},t_1)}^{(y \bar{\bf 1},t_2)}  =  \weight_{n,k;(\bar{v} - 2^{-1}n^{-2/3} \bar{\bf 1}  ,  t_1 + n^{-1})}^{(y \bar{\bf 1},t_2)} \, .
\end{equation}

The bound on the first probability in the corollary follows from~(\ref{e.properweightforward})
and an application of the two parts of Proposition~\ref{p.sumweight} 
with parameter settings ${\bf t_1} = t_1 + n^{-1}$, ${\bf t_2} = t_2$, ${\bf x} = x$ and ${\bf \bar{u}} = \bar{u} + 2^{-1}n^{-2/3}\bar{\bf 1}$.
The bounds $C_k \geq C$ and $c_k \leq c$ are also used.
The estimate on the second probability is reduced to this proposition via~(\ref{e.properweightbackward})
 and by noting the half-circle rotational symmetry of Brownian LPP that arises by reindexing its constituent curves $B(k,z)$ in the form $B(-k,-z)$. \qed

\noindent{\bf Proof of Proposition~\ref{p.sumweight}.}
The proposition's first assertion is the simpler of the two to derive.  To verify it, note that, since 
$\rho_{n,k,i ; (x \bar{\bf 1},t_1)}^{(\bar{u},t_2)}$ has starting and ending points $(x,\btone)$ and $(u_i,\bttwo)$, the weight of this $n$-zigzag is at most 
$\weight_{n;(x,t_1)}^{(u_i,t_2)}  $. Summing over $i \in \intint{k}$, we find that
\begin{equation}\label{e.swub}
 \weight_{n,k;(x \bar{\bf 1},t_1)}^{(\bar{u},t_2)} 
\leq \sum_{i=1}^k 
\weight_{n;(x,t_1)}^{(u_i,t_2)}  \, .
\end{equation}
Set $w_i = \tot^{-2/3} \big( u_i - x \big)$ and note that 
$$
 \weight_{n;(x,t_1)}^{(u_i,t_2)} =  \mc{L}^{\uparrow,t_2}_{n;(x,t_1)}\big( 1 , u_i \big) =  \tot^{1/3} \cdot \tot^{-1/3} \mc{L}^{\uparrow,t_2}_{n;(x,t_1)} \big( 1, x + \tot^{2/3} w_i  \big) =  \tot^{1/3} \, \scaledle^{\uparrow,t_2}_{n;(x,t_1)}(1,w_i) \, . 
$$
Applying {\em one-point upper tail} \rmreg(3) with parameter settings ${\bf z} = w_i$ and ${\bf s} = r k^{-1}$ to
the $(n\tot + 1)$-curve ensemble $\scaledle^{\uparrow,t_2}_{n;(x,t_1)}$, we see that the event
$$
 \Big\{   \tot^{-1/3} \cdot \weight_{n;(x,t_1)}^{(u_i,t_2)}  \geq  - 2^{-1/2}  \tot^{-4/3} \big( u_i -  x \big)^2  + r/k  \Big\} = \Big\{ \scaledle^{\uparrow,t_2}_{n;(x,t_1)}(1,w_i) + 2^{-1/2} w_i^2 \geq r/k \Big\} 
$$
has $\PP$-probability at most $C \exp \big\{ - c (r/k)^{3/2} \big\}$, provided that 
the hypotheses
$$
 \tot^{-2/3} \big\vert u_i - x \big\vert  \leq \rsc (n \tot)^{1/9} \,  \, \textrm{and} \, \, r \geq k 
$$
are satisfied. 
 Thus, (\ref{e.swub}) implies Proposition~\ref{p.sumweight}(1).

To derive the lower-tail bound Proposition~\ref{p.sumweight}(2), we will find~$k$ pairwise horizontally separate $n$-zigzags  of suitable weight, each of which begins at $(x,\btone)$ and which end successively at $(u_i,\bttwo)$, $i \in \intint{k}$.

\begin{figure}[ht]
\begin{center}
\includegraphics[height=9cm]{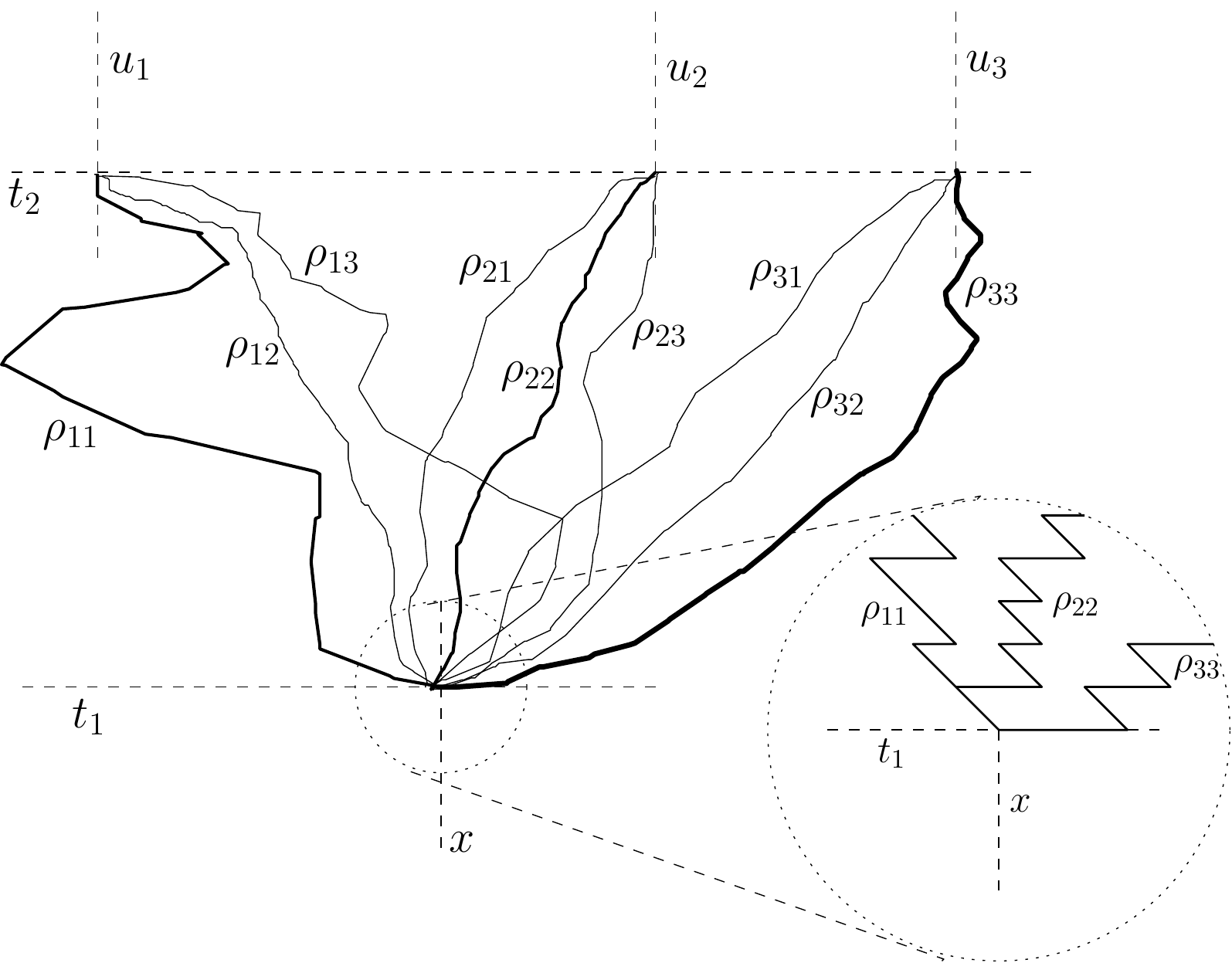}
\caption{The diagonal vector  $\big( \rho_{1,1}, \cdots, \rho_{k,k} \big)$ is illustrated in bold for an example with $k=3$. Under the microscope, we see how it is that these three zigzags may remain horizontally separate despite their sharing the birthplace $(x,t_1)$.}
\label{f.diag}
\end{center}
\end{figure}

Our technique for finding the zigzags is a diagonal argument: see Figure~\ref{f.diag}. For each $i \in \intint{k}$, consider the $k$-tuple multi-polymer watermelon
\begin{equation}\label{e.maximizer}
 \bigg( \rho_{n,k,i;(x  \bar{\bf{1}} ,t_1)}^{\big(u_i  \bar{\bf{1}}  ,  t_2 \big)}:  i \in \intint{k} \bigg) \, ,
 %\in \mc{P}_{n,k;(x,t_1)}^{\big(u_i {\bf 1} , t_2\big)} \, ,
\end{equation}
where here recall that $z  \bar{\bf{1}} \in \R^k$ denotes the vector each of whose components equals $z \in \R$. 
Our multi-polymer notation used here is specified in Section~\ref{s.multipolymer}. 
The object in question is the multi-polymer whose $k$ components are tethered at both endpoints, to   $(x,t_1)$ and $(u_i,t_2)$. As we explained when the notation was introduced, it is Lemma~\ref{l.severalpolyunique} that ensures the almost sure existence and uniqueness of the maximizer~(\ref{e.maximizer}), which for brevity we will denote by $\big( \rho_{i,1}, \cdots, \rho_{i,k} \big)$. 

We now consider the diagonal vector $\big( \rho_{i,i}: i \in \intint{k} \big)$. We first verify that this vector's components are pairwise horizontally separate zigzags that consecutively run from $(x,t_1)$ to $(u_i,t_2)$. The endpoint locations are not in doubt, so it is our task to verify the separateness condition between consecutively indexed zigzags.
% belongs to  $\mc{P}_{n,k;(x,t_1)}^{(\bar{u},t_2)}$. 
 In the language of Section~\ref{s.polyorder}, we want to check  that  $\rho_{i,i} \prec \rho_{i+1,i+1}$ for $i \in \intint{k-1}$. 
 We have that $\rho_{i,i} \prec \rho_{i,i+1}$,
because these players are components of a given multi-polymer, and we also have  
 $\rho_{i,i+1} \preceq \rho_{i+1,i+1}$ almost surely, by Lemma~\ref{l.tworelations}(1). Thus,  Lemma~\ref{l.tworelations}(2) implies the desired fact.

The confirmed property of the diagonal vector implies that
\begin{equation}\label{e.swlb}
\weight_{n,k;(x  \bar{\bf{1}} ,t_1)}^{(\bar{u},t_2)} \geq 
\sum_{i=1}^k \weight \big( \rho_{i,i} \big) \, ,
\end{equation}
where note that $\weight\big(\rho_{i,i}\big)$ is equal to
$\weight_{n,k,i;(x \bar{\bf{1}} ,t_1)}^{\big(u_i \bar{\bf{1}}   ,t_2\big)}$.

We now present a lemma regarding the values of these weights.
\begin{lemma}\label{l.kweight}
Let $(n,t_1,t_2) \in \N \times \R^2_<$ be a compatible triple. 
For  $(x,y) \in \R^2$, $y \geq x - 2^{-1} n^{1/3} \tot$, $k \in \N$ and $i \in \intint{k}$,
\begin{equation}\label{e.kweight}
 \mathcal{L}^{\uparrow;t_2}_{n;(x,t_1)}(1,y) - \sum_{j=2}^k \Big(\mathcal{L}^{\uparrow;t_2}_{n;(x,t_1)}(1,y) -  \mathcal{L}^{\uparrow;t_2}_{n;(x,t_1)}(j,y) \Big)  \leq  \weight_{n,k,i;(x\bar{\bf{1}},t_1)}^{\big(y\bar{\bf{1}},t_2\big)}
  \leq \mathcal{L}^{\uparrow;t_2}_{n;(x,t_1)}(1,y) \, .
\end{equation}
%where $y\bar{\bf{1}} \in \R^k$ denotes the vector each of whose components is $y$.
\end{lemma}
\noindent{\bf Proof.}
First note that the condition $y \geq x - 2^{-1} n^{1/3} \tot$ is assumed simply in order to assure that the ensemble  $\mathcal{L}^{\uparrow;t_2}_{n;(x,t_1)}$ is well defined at location $y$.
 Note that
$\mathcal{L}^{\uparrow;t_2}_{n;(x,t_1)}(1,y) = \weight_{n;(x,t_1)}^{(y,t_2)}$.
Since $\rho_{n,k,i;(x \bar{\bf{1}} ,t_1)}^{\big(u_i  {\bf{1}},t_2 \big)}$ is an $n$-zigzag from $(x,t_1)$ to $(u_i,t_2)$  for $i \in \intint{k}$, we have that
$\weight_{n,k,i;(x \bar{\bf{1}} ,t_1)}^{\big(y\bar{\bf{1}},t_2\big)} \leq  \weight_{n;(x,t_1)}^{(y,t_2)}$ for such $i$. Thus, we obtain the second bound in~(\ref{e.kweight}). As for the first, note that, by~(\ref{e.scaledweight}),
$$
 \sum_{j=1}^k  \mathcal{L}^{\uparrow;t_2}_{n;(x,t_1)}(j,y) \, = \,
 \sum_{j=1}^k  \weight_{n,k,j;(x \bar{\bf{1}}  ,t_1)}^{\big(y\bar{\bf{1}},t_2\big)}
  \, ,
$$
an equality that we may rewrite 
%using~(\ref{e.weightopcurve}) as 
\begin{equation}\label{e.wequal}
k \cdot \mathcal{L}^{\uparrow;t_2}_{n;(x,t_1)}(1,y)
  \, \, +  \, \, \sum_{j=2}^k \Big(\mathcal{L}^{\uparrow;t_2}_{n;(x,t_1)}(j,y) -  \mathcal{L}^{\uparrow;t_2}_{n;(x,t_1)}(1,y) \Big)  \, = \,   \sum_{j=1}^k \weight_{n,k,j;(x  \bar{\bf{1}} ,t_1)}^{\big(y\bar{\bf{1}},t_2\big)}   \  \, .
\end{equation}
Set $W_j = \weight_{n,k,j;(x  \bar{\bf{1}} ,t_1)}^{\big(y\bar{\bf{1}},t_2\big)}$ for $j \in \intint{k}$. 
Our remaining task is to derive the lower bound on $W_i$ given in the first inequality of~(\ref{e.kweight}), where $i \in \intint{k}$ is given. 
The right-hand side in the equality~(\ref{e.wequal}) takes the form $\sum_{j=1}^k W_j$.
The summand $W_j$ is at most $\weight_{n,(x,t_1)}^{(y,t_2)} = \mathcal{L}^{\uparrow;t_2}_{n;(x,t_1)}(1,y)$ for any $j \in \intint{k}$.  Suppose that we push
each $W_j$ up to its maximum possible value  $\mathcal{L}^{\uparrow;t_2}_{n;(x,t_1)}(1,y)$ 
for every  $j \in \intint{k}$ with $j \not= i$ -- and note that, in so doing, we push the value of $W_i$ down, because the equality~(\ref{e.wequal}) must be satisfied. That is, the value of the variable~$W_i$ in this scenario is determined by~ the satisfaction of~(\ref{e.wequal}); and this value offers a lower bound on the actual value of $W_i$. Since the resulting inequality is the first bound in~(\ref{e.kweight}), the proof of Lemma~\ref{l.kweight} is completed.  \qed

\medskip

In light of~(\ref{e.swlb}),  Proposition~\ref{p.sumweight}(2)  will follow once we show that, under the result's hypotheses,
\begin{equation}\label{e.sumbound}
\PP \bigg( 
 \tot^{-1/3} \, \sum_{i=1}^k 
\weight_{n,k,i;(x \bar{\bf{1}} ,t_1)}^{\big(u_i\bar{\bf{1}},t_2\big)}  \leq \, - \,  2^{-1/2} \sum_{i=1}^k \tot^{-4/3}(u_i - x)^2  \, \, - \, r \bigg) \leq  
3 k^2 
C_k \exp \Big\{ - 2^{-3}  k^{-3} c_k r^{3/2}     \Big\}  \, .
\end{equation}

To begin showing this, note that the occurrence of the left-hand event 
entails that the bound
$$ \tot^{-1/3} \,  
\weight_{n,k,i;(x \bar{\bf{1}} ,t_1)}^{\big(u_i\bar{\bf{1}},t_2\big)}  \leq \, - \,  2^{-1/2} \tot^{-4/3}(u_i - x)^2  \, \, - \,  r/k
$$
is satisfied for at least one index $i \in \intint{k}$. In view of Lemma~\ref{l.kweight}'s left-hand bound applied with ${\bf y} = u_i$ (for given $i \in \intint{k}$), the last inequality implies that either 
\begin{equation}\label{e.tsl.one}
\tot^{-1/3} \,  \mathcal{L}^{\uparrow;t_2}_{n;(x,t_1)}(1,u_i) \leq  \, - \,  2^{-1/2} \tot^{-4/3}(u_i - x)^2  \, \, - \, \frac{r}{2k} 
\end{equation}
or at least one among the inequalities
\begin{equation}\label{e.tsl.two} 
\tot^{-1/3} \Big(  \mathcal{L}^{\uparrow;t_2}_{n;(x,t_1)}\big(1,u_i\big)  \, - \,  \mathcal{L}^{\uparrow;t_2}_{n;(x,t_1)}\big(j,u_i\big) \Big) \geq  
  \frac{r}{2k(k-1)} \, ,
\end{equation}
indexed by $j \in \llbracket 2,k \rrbracket$,
is satisfied.

Condition~(\ref{e.tsl.one}) asserts that the highest curve of the normalized forward line ensemble, rooted at $(x,t_1)$  and of duration $\tot$,
$$
 \scaledle^{\uparrow;t_2}_{n;(x,t_1)}\big(1,  w_i \big)  = \tot^{-1/3} \,  \mathcal{L}^{\uparrow;t_2}_{n;(x,t_1)}\big(1, x + \tot^{2/3} w_i \big) \, ,
$$
evaluated at $w_i = \tot^{-2/3} \big( u_i - x \big)$, is at most $- \,  2^{-1/2} \tot^{-4/3}(u_i - x)^2   \, - \,  2^{-1} r k^{-1}$. 
Applying one-point lower tail $\rmreg(2)$ to the $(n \tot + 1)$-curve ensemble $\scaledle^{\uparrow;t_2}_{n;(x,t_1)}$,
%, $n \in \N$, (where Proposition~\ref{p.scaledreg} assures that this sequence is regular), 
with parameter settings ${\bf z} = w_i$ and ${\bf s} = 2^{-1} rk^{-1}$,
  we find that this eventuality has probability at most $C \exp \big\{ - c \, 2^{-3/2} k^{-3/2} r^{3/2} \big\}$.
  This application of  $\rmreg(2)$ may be carried out if the hypotheses 
 $$
 \tot^{-2/3} \vert u_i - x \vert \leq \rsc (n \tot)^{1/9}    \, \,  \textrm{and} \, \, r \in \big[2k, 2k (n \tot)^{1/3} \big]
 $$ 
  are satisfied.

When the bound~(\ref{e.tsl.two}) is rewritten in normalized coordinates, it asserts that
\begin{eqnarray*}
 & & \Big( \scaledle^{\uparrow;t_2}_{n;(x,t_1)}\big(1,  w_i \big) + 2^{-1/2} w_i^2 \Big)  - \Big(  \scaledle^{\uparrow;t_2}_{n;(x,t_1)}\big(j,  w_i \big) + 2^{-1/2} w_i^2 \Big)   \\
 & = & \tot^{-1/3} \,  \mathcal{L}^{\uparrow;t_2}_{n;(x,t_1)}\big(1, x + \tot^{2/3} w_i \big) \, - \, 
\tot^{-1/3} \,  \mathcal{L}^{\uparrow;t_2}_{n;(x,t_1)}\big(j, x + \tot^{2/3} w_i \big) \, ,
\end{eqnarray*}
when 
evaluated at $w_i = \tot^{-2/3} \big( u_i - x \big)$, is at least $\frac{r}{2k(k-1)}$. The bound forces at least one of two cases: in the above difference of two terms,
either the first term is at least $\frac{r}{4k(k-1)}$ or the second is at most  $-\frac{r}{4k(k-1)}$.

The first case is handled by applying one-point upper tail $\rmreg(3)$
to $\scaledle^{\uparrow;t_2}_{n;(x,t_1)}$, taking ${\bf z} = w_i$ and ${\bf s} = \frac{r}{4k(k-1)}$, with an upper bound of
$$
C \exp \Big\{ - c \big(\tfrac{r}{4k(k-1)}\big)^{3/2} \Big\} 
$$
being found on the probability of the event in question. This use of $\rmreg(3)$ may be made when
$$
 \tot^{-2/3} \big\vert u_i - x \big\vert \leq \rsc (n \tot)^{1/9} \, \, \textrm{and} \, \,  r \geq 4k(k-1) \, .
$$
 The second case is treated by applying pointwise lower tail Proposition~\ref{p.mega}(1)
%  Proposition~\ref{p.othercurves}
 to the $(n\tot+1)$-curve ensemble $\mc{L}_n =  \scaledle^{\uparrow;t_2}_{n;(x,t_1)}$,  with parameter settings ${\bf k} = j$, ${\bf z} = w_i$ and   ${\bf s} = \frac{r}{4k(k-1)}$. The event in question is thus found to have $\PP$-probability at most   
 $$
  C_j \exp \Big\{ - c_j \big(\tfrac{r}{4k(k-1)}\big)^{3/2} \Big\} \, .
 $$
 These applications of  Proposition~\ref{p.mega}(1)
   are made for each $j \in \llbracket 2, k \rrbracket$, and the $k-1$ applications may be made provided that 
  $$
  n \tot  \geq k 
  \vee  (c/3)^{-18} \vee  6^{36} \, , \, 
  \tot^{-2/3} \big\vert u_i - x \big\vert  \leq  2^{-1} \rsc  (n \tot)^{1/18} \, \,  \textrm{and}  \, \,  r \leq 8 k (k-1) (n \tot)^{1/18} \, .
 $$ 
  We find then that  
 the probability that~(\ref{e.tsl.two}) is satisfied for a given pair $(i,j) \in \intint{k}  \times \llbracket 2, k \rrbracket$ is at most 
$$
C  \exp \Big\{ - c \big(\tfrac{r}{4k(k-1)}\big)^{3/2} \Big\} \, + \, C_j \exp \Big\{ - c_j \big(\tfrac{r}{4k(k-1)}\big)^{3/2} \Big\} \, .
$$

Gathering together these inferences by means of a union bound over this set of index pairs, we find that the left-hand side of~(\ref{e.sumbound}) is at most 
$$
C k \exp \big\{ - c 2^{-3/2} k^{-3/2} r^{3/2} \big\} \, + \, 
C k^2 \exp \Big\{ - c \big(\tfrac{r}{4k(k-1)}\big)^{3/2} \Big\} \, + \, k \sum_{j=2}^k 
C_j \exp \Big\{ - c_j \big(\tfrac{r}{4k(k-1)}\big)^{3/2} \Big\}  \, .
$$
Noting that 
$C_k$ is increasing,  $c_k$ is decreasing, $C_k \geq C$ and $c_k \leq c$, we verify~(\ref{e.sumbound}).
This completes the proof of Proposition~\ref{p.sumweight}(2). \qed

\subsection{The road map's third challenge: the solution in overview}\label{s.thirdchallenge}
We have assembled the elements needed to implement the road map and thus to prove Theorem~\ref{t.disjtpoly}, with the exception of addressing the third challenge, which we labelled {\em final polymer comparison} in Section~\ref{s.roadmap}.

Recall that the road map advocates the construction of a system of $\kay$ near polymers, each running from $(x,-\e^{3/2})$ to $(y,1+\e^{3/2})$, and pairwise disjoint otherwise.
These will be zigzags, to be called $\rho_1$, $\rho_2$, $\cdots$, $\rho_\kay$, that during the time interval $[0,1]$
follow the course of the $\kay$ polymers, to be called $\phi_1$, $\phi_2$, $\cdots$, $\phi_\kay$, whose existence is ensured by 
the occurrence of the event  $\maxpoly_{n;([x-\e,x+\e],0))}^{([y-\e,y+\e],1)} \geq \kay$. They begin by following elements in a forward bouquet of lifetime $[-\e^{3/2},0]$ with shared endpoint $(x,-\e^{3/2})$
and end by following elements in a backward bouquet of lifetime $[1,1+\e^{3/2}]$
with shared endpoint $(y,1+\e^{3/2})$.

Now the third challenge involves arguing that each of the $\rho_i$ is indeed a near polymer, with a shortfall in weight from the maximum $\weight_{n;(x,-\e^{3/2})}^{(y,1+\e^{3/2})}$ of order $\e^{1/2}$. As we pointed out in Section~\ref{s.roadmap}, the third challenge seems merely a restatement of the overall problem. A little more detail is needed to explain even in outline how we propose to solve the third challenge, and we now offer such an outline.

Our task is to establish that, for each $i \in \intint{\kay}$, the quantity $\weight \big( \rho_{n;(x,-\e^{3/2})}^{(y,1+\e^{3/2})} \big) - \weight(\rho_i)$, which is necessarily non-negative, is in fact at most $O(\e^{1/2})$.
Now the value $\weight(\rho_i)$ is naturally written as a sum of three terms: the weight of a forward bouquet element, the weight of the polymer $\phi_i$, which moves from $[x-\e,x+\e] \times \{0\}$ to $[y-\e,y+\e] \times \{ 1 \}$, and the weight of a backward bouquet element. The first and third weights are known to be of order $\e^{1/2}$ typically by Corollary~\ref{c.bouquetreg}. The weight
$\weight \big( \rho_{n;(x,-\e^{3/2})}^{(y,1+\e^{3/2})} \big)$ may correspondingly be split into three terms, by splitting the polymer in question at the pair of times $(0,1)$. To resolve the third challenge, we want to argue that:
\begin{itemize}
\item The first and third pieces of the polymer  $\rho_{n;(x,-\e^{3/2})}^{(y,1+\e^{3/2})}$, with lifetimes $[-\e^{3/2},0]$
and  $[1,1+\e^{3/2}]$, have weight of order $\e^{1/2}$.
\item The second piece, the grand middle section with lifetime $[0,1]$, has a weight that differs from any of the weights $\weight(\phi_i)$ for $i \in \intint{\kay}$, by an order of $\e^{1/2}$.
\end{itemize}

\begin{figure}[ht]
\begin{center}
\includegraphics[height=9cm]{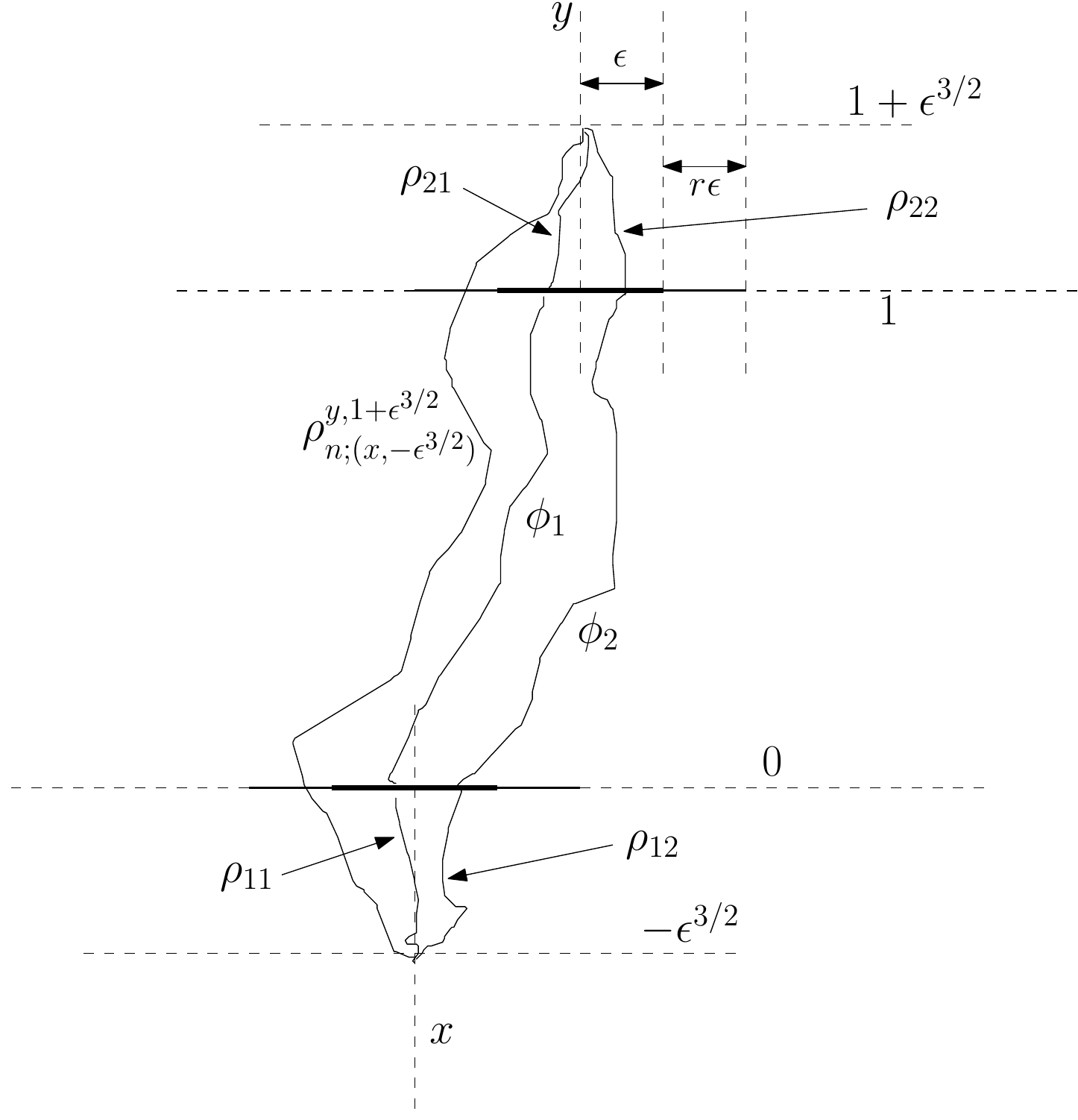}
\caption{Solving the third challenge of the road map. The planar intervals $I \times \{ 0 \}$ and $J \times \{ 1 \}$, with $I = [x-\e,x+\e]$ and $J = [y-\e,y+\e]$, are drawn with thick bold lines.
Extended intervals 
$I^+ \times \{ 0 \}$ and $J^+ \times \{ 1 \}$ are also drawn in solid lines. Here, $I^+ = [x-(r+1)\e,x+(r+1)\e]$ and $J^+ = [y-(r+1)\e,y+(r+1)\e]$ with $r > 0$ given. 
The notation will be used in the upcoming proof of Theorem~\ref{t.disjtpoly}. The polymer  $\rho_{n;(x,-\e^{3/2})}^{(y,1+\e^{3/2})}$ typically visits these extended intervals at times zero and one, and, in this event, the weight of its $[0,1]$-duration subpath may be closely compared to that of the polymers $\phi_1$ and $\phi_2$. The $\rho_{ij}$ notation will be used in the proof to indicate paths in the forward and backward bouquets.}
\label{f.surgery}
\end{center}
\end{figure}

Now to reach these conclusions, we will in fact seek control of the geometry of the polymer $\rho_{n;(x,-\e^{3/2})}^{(y,1+\e^{3/2})}$. This we will do by invoking one of the main results of the present article, Theorem~\ref{t.polyfluc}. We will learn that, typically, this polymer at time zero is located at distance of order $\e$
from $x$, and at time one, at distance of order $\e$ from $y$.
With this understanding, we see that the polymer's first and third pieces do not suffer significant lateral shifts, so that Corollary~\ref{c.maxminweight}
will ensure that these pieces have the desired order $\e^{1/2}$ weight.
The same understanding means that the middle section of $\rho_{n;(x,-\e^{3/2})}^{(y,1+\e^{3/2})}$ and the middle section $\phi_i$ of any $\rho_i$
have starting and ending points whose locations differ by an order of $\e$:
close to $x$ or to $y$ respectively. This is crucial information, because it permits us to invoke Theorem~\ref{t.differenceweight}
to conclude that the middle section weights indeed differ by an order of~$\e^{1/2}$. These ideas are illustrated by Figure~\ref{f.surgery}.

Armed with this elaboration of the third challenge in the road map, we are ready to derive Theorem~\ref{t.disjtpoly}.

\subsection{Proof of Theorem~\ref{t.disjtpoly}.} \label{s.proof}
By the scaling principle,
it suffices to prove the result when $(t_1,t_2) = (0,1)$, and this we now do. 
Set $I = [x-\e,x+\e]$ and $J = [y - \e,y+\e]$.
%, and let $\bar{u} \in I^m_<$ and $\bar{v} \in J^m_<$ be given. 
For a given parameter $r > 0$, we also define the extended intervals  $I^+ = [x-(r+1)\e,x+(r+1)\e]$ and $J^+ = [y - (r+1)\e,y+(r+1)\e]$ that respectively contain $I$ and $J$.

On the event that $\maxpoly_{n;(I,0))}^{(J,1)} \geq \kay$,
there exist collections of $\kay$ disjoint polymers that make the journey from $I \times \{ 0 \}$ to $J \times \{ 1 \}$. We now seek to choose one of them. 
%The concerned polymers may begin and end with horizontal segments and, if they do so, they infringe, albeit microscopically, on randomness in the before- and after-life regions $\R \times (-\infty,0]$
%and $\R \times [1,\infty)$. In order not to trespass into these regions, we recall the notion of a proper polymer from Section~\ref{s.closure}. In fact, using the theory from that section,
%we are able 
%on the event that $\maxpoly_{n;(I,0))}^{(J,1)} \geq m$ to specify (and henceforth fix) an element $\big(\bar{U},\bar{V} \big) \in I^m_< \times J^m_<$ such that there exist pairwise horizontally separate proper $n$-polymers with consecutive starting and ending points $(U_i,0)$ and $(V_i,1)$.
Indeed, using the language of Section~\ref{s.closure}, we may choose  $\big(\bar{U},\bar{V} \big)$ to be the lexicographically minimal element in the closure of the set  $\disjtindex_{n,\kay;(I,0)}^{(J,1)}$.
This definition makes sense because we are dealing with a closed set. Invoking 
Lemma~\ref{l.disjthorsep}, we see that $\big(\bar{U},\bar{V} \big)$ is an element of  $\horsepindex_{n,\kay;(I,0)}^{(J,1)}$.
That is, we have explicitly selected an $\kay$-tuple of polymers that, while not necessarily pairwise disjoint, is pairwise horizontally separate; this is enough for the upcoming surgery.
% as we wished to confirm. 

It is helpful to recall from the road map that it is our aim to show that the occurrence of the event $\maxpoly_{n;(I,0))}^{(J,1)} \geq \kay$ {\em typically}
entails 
$\neargeod_{n,\kay;(x,-\e^{3/2})}^{(y,1 + \e^{3/2})}( \eta )$, where the parameter~$\eta$ will soon be selected explicitly to have order $\e^{1/2}$. 
What do we mean by typically? 
Typical behaviour means a collection of  provably standard circumstances needed to undertake our rewiring surgery successfully. It has four types:
\begin{itemize}
\item $\fbr$/$\bbr$:
The $\e^{1/2}$-order weight of bouquets.
\item  $\pwr$: The $\e^{1/2}$-order weight of any $[-\e^{3/2},0]$- or $[1,1+\e^{3/2}]$-lifetime polymer respectively near $x$ and $y$.
\item $\pdr$: Control on the geometry of $\rho_{n;(x,-\e^{3/2})}^{(y,1 + \e^{3/2})}$ at times zero and one.
\item $\lwr$: The one-half power weight difference for middle section polymers. 
\end{itemize}
We will shortly define a {\em favourable surgical conditions} event, which specifies exactly what we mean by `typical', as an intersection of events to which
the labels just listed correspond.   

Before offering the definition, we introduce a variation of the $\pdr$ notation specified in  Subsection~\ref{s.polyflucintro}. When the event 
$\pdr_{n;(x,t_1)}^{(y,t_2)}\big(a,r\big)$ is specified in~(\ref{e.pdr}), the parameter $a$ refers to the proportion of the polymer lifetime $t_{1,2}$
which has elapsed at the intermediate time $(1-a)t_1 + at_2$. We now wish to allude to the intermediate time directly. To do so,
 we will employ a square bracket notation. Taking $t \in [t_1,t_2] \cap n^{-1}\Z$,
we will write
$$
\pdr_{n;(x,t_1)}^{(y,t_2)}\big[t,r\big]  \, = \, \bigg\{  \,  \Big\vert \, \rho_{n;(x,t_1)}^{(y,t_2)} (  t ) - \ell_{(x,t_1)}^{(y,t_2)}  (  t ) \, \Big\vert \, \leq \, r \big(  (t - t_1) \wedge (t_2 - t) \big)^{2/3} \bigg\} \, ,
$$
so that 
$\pdr_{n;(x,t_1)}^{(y,t_2)}\big[t,r\big]$  equals 
$\pdr_{n;(x,t_1)}^{(y,t_2)}\big( (t-t_1) \tot^{-1}  \, , r \big)$. We also extend this notation so that the first square bracket argument $t$ is replaced by a pair of such times. This notation refers to the intersection of the two single-time $\pdr$
events.  

We now may define the favourable surgical conditions\hfff{fsc} event  
 $\fsc_{n;(I,0)}^{(J,1)}\big(\kay; \bar{U},\bar{V} ; \e,r \big)$. 
 The quantities $\e > 0$ and $r > 0$ will retain their roles as parameters in this event for the remainder of the proof of Theorem~\ref{t.disjtpoly}.
 The event will only be considered when the event $\maxpoly_{n;(I,0))}^{(J,1)} \geq \kay$ occurs, so that the 
 random vectors~$\bar{U}$ and $\bar{V}$ that are used as parameters in the event's definition will always make sense. 
 The new event is defined to equal
\begin{eqnarray*}
 &  & \lwr_{n;(I^+,0)}^{(J^+,1)}\big(2(r+1)\e,r\big) \cap \pdr_{n;(x,-\e^{3/2})}^{(y,1 + \e^{3/2})} \big[ \{0,1\} , r \big] \\
  & & \quad \cap \, \fbr_{n,\kay;(x,-\e^{3/2})}^{(\bar{u},0)}(r) \cap \bbr_{n,\kay;(\bar{v},1)}^{(y,1+\e^{3/2})}(r) \\
  & & \quad \cap \, \pwr_{n;(x,-\e^{3/2})}^{(I^+,0)}(r^2) \cap \pwr_{n;(J^+,1)}^{(y,1+\e^{3/2})}(r^2) \, .
\end{eqnarray*}
The $\lwr$, $\fbr$/$\bbr$ and $\pwr$ events
have been defined in Sections~\ref{s.polyweightreg},~\ref{s.bouquet} and~\ref{s.usefultool}.
In the use of the new square bracket notation for $\pdr$, the polymer is being controlled at times zero and one. The reason of the choice of parameter $r^2$ in the two $\pwr$ events will become clearer, and will be discussed, in due course.

Note that the quantity $\e^{3/2}$ plays the role of the duration of certain zigzags concerned in this definition, these zigzags having lifetime $[-\e^{3/2},0]$ and $[1,1+\e^{3/2}]$.
In order that these zigzags begin and end at vertical coordinates in the $n^{-1}$-mesh, we will insist throughout that $\e > 0$ satisfies $\e^{3/2} \in n^{-1} \Z$. It would seem then that this condition should enter as a hypothesis of Theorem~\ref{t.disjtpoly}, though it does not; we omit it because, in view of the hypotheses~(\ref{e.epsilonbound}) and~(\ref{e.nlowerbound}), the condition can be forced by a multiplicative adjustment in $\e$ that differs from one by  order say $10^{-100}$, 
and this causes a tiny similar adjustment in the theorem's conclusion, which the reader may readily confirm is easily absorbed by tightening estimates during the proof.  
In summary, though  $\e > 0$ is $n$-dependent, it may in practice  be considered to be fixed at a given small positive value, the upper bound specified by several upcoming demands.  We may later omit mention of the mesh membership condition  $\e^{3/2} \in n^{-1} \Z$. 

The next two results assert that favourable surgical conditions are indeed typical, and that, in the presence of these conditions, the occurrence of $\maxpoly_{n;(I,0))}^{(J,1)} \geq \kay$ indeed forces the occurrence of
$\neargeod_{n,\kay;(x,-\e^{3/2})}^{(y,1 + \e^{3/2})}( \eta )$ with $\eta = O(\e^{1/2})$. 
%Note that the role of $m$ is played by $k$ in each result.
\begin{lemma}\label{l.fsc}
Suppose that $n \in \N$, $k \in \N$, $x,y \in \R$, $\e > 0$ and $r > 0$ satisfy 
$$
n \geq  10^{40}  k  c^{-18} \e^{-75/2} \max \big\{ 1  \, , \,   \vert x - y \vert^{36} \big\} \, \, \, \textrm{and} \, \, \, \vert x - y \vert \leq \big( 2 (r+1) \e \big)^{-1/2} \, .
$$
 Suppose also that 
$$
 r \geq \max \Big\{  \, 10^9 c_1^{-4/5}  \, , \, 15 C^{1/2} \, , \, 4 k^2  \, , \, 70 \e^{1/2} \vert x - y \vert   \, \Big\} \, , \,  \e  \leq 2^{-5} (r+1)^{-1}         \, \, \textrm{ and } \, \, r  \leq  2^{-6} c \, \e^{25/6} n^{1/36} \, .
 $$
 Then 
\begin{eqnarray*}
 & &  \PP \Big(  \maxpoly_{n;(I,0))}^{(J,1)} \geq k  \, , \, \neg \,   \fsc_{n;(I,0)}^{(J,1)}\big(k; \bar{U},\bar{V} ; \e,r  \big) \Big) \\
 & \leq &     14062 \, k^2 r C_k \exp \big\{ -  10^{-11}   k^{-3} c_k  r^{3/4} \big\}   \, .
\end{eqnarray*}
\end{lemma}

\begin{proposition}\label{p.doubletie}
Let $n,k \in \N$, $n \geq k \geq 1$. Let $\e > 0$ and $x,y \in \R$
satisfy $n \e^{3/2} \geq 10^2$, $\vert x - y \vert \leq \e^{-1/2}$, and let $r \geq 1$. Then, whenever 
$\big\{ \maxpoly_{n;([x-\e,x+\e],0))}^{([y-\e,y+\e],1)} \geq k  \big\} \cap \fsc_{n;(I,0)}^{(J,1)}\big(k; \bar{U},\bar{V} ; \e,r   \big)$ occurs,
\begin{equation}\label{e.doubletie} 
\weight_{n,k;(x \bar{\bf 1},  - \e^{3/2})}^{(y \bar{\bf 1},1   + \e^{3/2})} \geq k \cdot \weight_{n;(x,  - \e^{3/2})}^{(y,1   + \e^{3/2})} \, - \,
 15k r^2  \e^{1/2}
\, ,
\end{equation}
which is to say, the event 
$\neargeod_{n,k;(x,-\e^{3/2})}^{(y,1 + \e^{3/2})}( \eta )$ occurs where $\big(1+2\e^{3/2} \big)^{1/3} \eta = 15k r^2  \e^{1/2}$.
\end{proposition}

We now apply these two results to close out the proof of Theorem~\ref{t.disjtpoly}
and then prove them in turn. 

Note that, for the value of $\eta$ specified in Proposition~\ref{p.doubletie},  $\neargeod_{n,k;(x,-\e^{3/2})}^{(y,1 + \e^{3/2})}( \eta )$ is a subset of $\neargeod_{n,k;(x,-\e^{3/2})}^{(y,1 + \e^{3/2})}\big( 15k r^2  \e^{1/2} \big)$.
The proposition 
%with ${\bf k} = m$ 
thus implies that 
\begin{eqnarray*}
& & \PP \Big( \maxpoly_{n;(I,0))}^{(J,1)} \geq \kay  \, , \, \fsc_{n;(I,0)}^{(J,1)}\big( \kay ; \bar{U},\bar{V}  ; \e,r  \big)  \Big) \\
& \leq & \PP \Big( \neargeod_{n,\kay;(x,-\e^{3/2})}^{(y,1 + \e^{3/2})}\big( 15 \kay r^2  \e^{1/2}  \big) \Big) \, ,
\end{eqnarray*}
provided that $\vert x - y \vert \leq \e^{-1/2}$ and $r \geq 1$.  
 
We now apply
Corollary~\ref{c.neargeod.t} 
to bound above the right-hand probability.  
Parameter settings are ${\bf t_1} = -\e^{3/2}$, ${\bf t_2} = 1 + \e^{3/2}$,  ${\bf k} = \kay$, ${\bf x} = x$,  ${\bf y} = y$ and   ${\bm \eta} =  15\kay r^2  \e^{1/2}$.
Note that, since ${\bm \eta} \geq \e$, the corollary's hypotheses are satisfied when 
$$
15\kay r^2  \e^{1/2}  < (\eta_0)^{\kay^2} \, , \,
 \kay \geq 2 \, ,  \,  (1 + 2\e^{3/2}) n/2 \geq \kay \vee \, (K_0)^{\kay^2} \big( \log \e^{-1} \big)^{K_0}
 $$ 
 and   $(1 + 2\e^{3/2})^{-2/3} \vert y - x \vert \leq a_0 n^{1/9}$. (It is also necessary that $\big(n,1-\e^{3/2},1+\e^{3/2}\big) \in \N \times \R^2_<$ be a compatible triple, which amounts to the already imposed condition that 
 $n \e^{3/2} \in \N$.) When these conditions are met, the conclusion of the corollary implies that 
\begin{eqnarray*}
&  & \PP \Big( \neargeod_{n,\kay;(x,-\e^{3/2})}^{(y,1 + \e^{3/2})}\big( 15 \kay r^2  \e^{1/2}  \big) \Big) \\
& \leq & 15^{\kay^2} \kay^{\kay^2} r^{2\kay^2} \e^{(\kay^2 - 1)/2}    \exp \big\{ \beta_\kay \big( \log \e^{-1} \big)^{5/6} \big\}
 \, ,
\end{eqnarray*}
where we again used $\eta \geq \e$ to write the final term in the product.

Applying Lemma~\ref{l.fsc},
% with ${\bf k} = m$, 
we find that
\begin{eqnarray}
  \PP \Big( \maxpoly_{n;(I,0))}^{(J,1)} \geq \kay \Big) 
  & \leq &  15^{\kay^2} \kay^{\kay^2} r^{2\kay^2} \e^{(\kay^2 - 1)/2}     \exp \big\{ \beta_\kay \big( \log \e^{-1} \big)^{5/6} \big\} \nonumber \\ 
  & & +  \,  \,     14062 \, \kay^2 r C_\kay \exp \big\{ -  10^{-11}   \kay^{-3} c_\kay  r^{3/4} \big\}  \label{e.explast} \, .
\end{eqnarray}
This application requires that $n \geq  10^{40}  \kay  c^{-18} \e^{-75/2} \max \big\{ 1  \, , \,   \vert x - y \vert^{36} \big\}$, $\vert x - y \vert \leq \big( 2 (r+1) \e \big)^{-1/2}$, 
$$
 r \geq \max \Big\{  \, 10^9 c_1^{-4/5}  \, , \, 15 C^{1/2} \, , \, 4 \kay^2  \, , \, 70 \e^{1/2} \vert x - y \vert   \, \Big\} \, , \,  \e  \leq 2^{-5} (r+1)^{-1}         \, \, \textrm{ and } \, \, r  \leq  2^{-6} c \, \e^{25/6} n^{1/36} \, .
 $$

We now set  the parameter $r > 0$ so that   $10^{-11}   \kay^{-3} c_\kay  r^{3/4} =  2^{-1} (\kay^2 - 1) \log \e^{-1}$.  
The choice is made in order that   the exponential term in~(\ref{e.explast}) equal $\e^{(\kay^2 - 1)/2}$.
Thus, 
\begin{equation}\label{e.rchoice}
 r = 10^{44/3} \kay^4 c_\kay^{-4/3}  2^{-4/3} (\kay^2 - 1)^{4/3} \big( \log \e^{-1} \big)^{4/3} \, .
\end{equation}
 Note that 
$$
 r \leq  10^{15}  c_\kay^{-4/3}  \kay^{20/3} \big( \log \e^{-1} \big)^{4/3} \, .
$$
We then obtain 
\begin{eqnarray*}
 & & \PP \Big( \maxpoly_{n;(I,0))}^{(J,1)} \geq \kay \Big) \\
  & \leq &  15^{\kay^2}  10^{30\kay^2}  \kay^{43\kay^2/3}  c_\kay^{-8\kay^2/3}   \big( \log \e^{-1} \big)^{4\kay^2} \e^{(\kay^2 - 1)/2}    \exp \big\{ \beta_\kay \big( \log \e^{-1} \big)^{5/6} \big\} \\ 
  & & +  \,   10^{20}  \kay^{26/3}      c_\kay^{-4/3}   C_\kay  \big( \log \e^{-1} \big)^2 \cdot \e^{(\kay^2 - 1)/2}  \\
  & \leq &  10^{32\kay^2}  \kay^{15\kay^2}  c_\kay^{-3\kay^2} C_\kay   \big( \log \e^{-1} \big)^{4\kay^2} \e^{(\kay^2 - 1)/2}    \exp \big\{ \beta_\kay \big( \log \e^{-1} \big)^{5/6} \big\} \, ,
\end{eqnarray*}
where we used  $c_\kay \leq c_1$ and $\e \leq e^{-1}$ in the first inequality and  $\kay \geq 1$, $C_\kay \geq 1$, $c_\kay \leq 1$ and $\e \leq e^{-1}$ in the second.
This completes the proof of Theorem~\ref{t.disjtpoly} in the case that $t_1=0$ and $t_1=1$; as we noted at the outset, this special case implies the general one. \qed

\medskip

\noindent{\bf Proof of Lemma~\ref{l.fsc}.} The $\fsc$ event is an intersection of six events. Each of these events is known to be typical, by several results proved or cited in this article. Thus, the present proof need merely gather together the estimates, invoked with suitable parameter settings.

 By Corollary~\ref{c.ordweight}
 % was Theorem~\ref{t.differenceweight} 
 with ${\bm \e} = 2(r+1)\e$, ${\bf R} = r$, ${\bf x} = x - (r+1)\e$ and ${\bf y} = y - (r+1)\e$,  
$$
 \PP \Big( \neg \, \lwr_{n;(I^+,0)}^{(J^+,1)}\big(2(r+1)\e,r\big) \Big) \leq   
 10032 \, C  \exp \big\{ - c_1 2^{-23} r^{3/2}     \big\} 
% was 25064 \, C  \exp \big\{ - c_1 2^{-31} r^{3/2}     \big\} 
$$
when $n  \geq 10^{32} c^{-18}$  and
$$
2(r+1)\e \in (0,2^{-4}] \, , \, \big\vert x -  y  \big\vert \leq \big( 2(r+1) \e \big)^{-1/2} \wedge 2^{-2} 3^{-1} \rsc  n^{1/18} \, \, \textrm{and} \, \, 
 r \in \big[ 2 \cdot 10^4 \, , \,   10^3 n^{1/18} \big] \, .
 $$

We apply Proposition~\ref{p.polyfluc}, and the remark following this proposition, in order to find that  
$$
 \PP \Big( \neg \,  \pdr_{n;(x,-\e^{3/2})}^{(y,1 + \e^{3/2})} \big[ \{0,1\} , r \big] \Big) \leq 22 Cr    \exp \big\{ - 10^{-11} c_1  r^{3/4} \big\}  \, .
$$
When the proposition is applied, it is with parameter settings
 ${\bf x} = x$, ${\bf y} = y$, ${\bf t_1} = -\e^{3/2}$, ${\bf t_2} = 1 + \e^{3/2}$, ${\bf r} = r$
and with ${\bf a}$ chosen equal to  $1 - \e^{3/2}(1 + 2\e^{3/2})^{-1}$.
When the remark is applied, we instead take ${\bf a} = \e^{3/2}(1 + 2\e^{3/2})^{-1}$.
Note that ${\bf a} \wedge (1 - {\bf a} ) \geq \e^{3/2}/2$ since $\e \leq 2^{-2/3}$.
Using  ${\bf a}^{-1} \vee (1 - {\bf a} )^{-1} \leq 2 \e^{-3/2}$ and ${\bf \tot} \in [1,2]$, we see that the hypotheses of this application are met provided that
$$
 n    \geq  \max \bigg\{ 
 10^{32} \big(\e^{3/2}/2 \big)^{-25} c^{-18} \, \, , \, \, 
 10^{24} c^{-18}  \big(\e^{3/2}/2 \big)^{-25} \vert x - y \vert^{36} (1 + 2\e^{3/2})^{-24} 
 \bigg\} \, ,
$$
$$ 
 r  \geq 
 \max \bigg\{  10^9 c_1^{-4/5} \, \, , \, \,  15 C^{1/2} \, \, , \, \,  87  \big(\e^{3/2}/2 \big)^{1/3}  (1 + 2\e^{3/2})^{-2/3}  \vert x - y \vert   \bigg\}
 $$
and, since $\e^{3/2}(1 + 2\e^{3/2})^{-1}  \geq 2^{-1} \e^{3/2}$,
$r \leq 3   \big(   2^{-1} \e^{3/2}  \big)^{25/9} n^{1/36} (1 + 2\e^{3/2})^{1/36}$.

Corollary~\ref{c.bouquetreg} is applied to find that
\begin{equation}\label{e.bouquetapp}
 \PP \Big( \neg \, \fbr_{n,k;(x,-\e^{3/2})}^{(\bar{U},0)}(r) \Big) \, + \, \PP \Big(  \neg \, \bbr_{n,k;(\bar{V},1)}^{(y,1+\e^{3/2})}(r) \Big) \leq
8 k^2 
C_k \exp \big\{ - 2^{-3}  k^{-3} c_k r^{3/2} \big\}   \, .
\end{equation}
Here, parameter settings are  ${\bf k} = k$,
${\bf t_1} = - \e^{3/2}$, ${\bf t_2} = 0$, ${\bf x} = x$, ${\bf y} = y$, ${\bf \bar{u}} = \bar{U}$, ${\bf \bar{v}} = \bar{V}$
and ${\bf r}  = r$.  Since $\vert U_i - x \vert \vee  \vert V_i - y \vert \leq \e$ for each $i \in \intint{k}$,
this application may be made provided that
\begin{eqnarray*}
n \e^{3/2} - 1  & \geq & \max \Big\{  k  \, , \,  3^{18} c^{-18}  \, ,  \, 6^{36} \, , \,   \rsc^{-9}(1 - n^{-1}\e^{-3/2})^{-6}(1 + 2^{-1} n^{-2/3} \e^{-1} )^9 \, , \\
& & \qquad \qquad    \qquad \qquad   \qquad    2^{18}  \rsc^{-18} (1 - n^{-1}\e^{-3/2})^{-12}(1 + 2^{-1} n^{-2/3} \e^{-1} )^{18} \Big\} 
\end{eqnarray*}
and  $r \in \big[ 4k^2  ,  9 k^2 \e^{1/2} (1 - n^{-1}\e^{3/2})^{1/3} n^{1/3} \big]$.

There is a subtlety in this application of 
Corollary~\ref{c.bouquetreg}  that deserves mention. In fact, it is in this application that the microscopic detail discussed after Proposition~\ref{p.sumweight} is implicated. Our parameter choice $\big( \bar{u}, \bar{v} \big) = \big( \bar{U}, \bar{V} \big)$ is a random one. Even if, as is the case, our random variable $\big( \bar{U}, \bar{V} \big)$
almost surely verifies the necessary bounds for this use, it is not formally admissible in this role. However, recall that $\big( \bar{U}, \bar{V} \big)$ has been specified to be measurable with respect to randomness in the  region $\R \times [0,1]$. This random variable is thus independent of the randomness in the open region
$\R \times \big((-\infty,0) \cup (1,\infty) \big)$ that specifies the two events in~(\ref{e.bouquetapp}) with which the application of the corollary is concerned.
It is the use of proper weights in the definition of the events $\fbr$ and $\bbr$ which ensures that these events are measurable with respect to the randomness in $\R \times \big((-\infty,0) \cup (1,\infty) \big)$.  It is this independence which renders the application admissible.

Note that
$$
\pwr_{n;(x,-\e^{3/2})}^{([x - (r+1)\e,x+(r+1)\e],0)}(r^2) \subseteq \maxmin_{n;([x,x+\e],-\e^{3/2})}^{([y,y+d\e],0)}(r^2) 
$$ 
where $y = x - (r+1)\e$  and $d = \lceil 2(r+1) \rceil$, and that
$$
 \pwr_{n;([y - (r+1)\e,y+(r+1)\e],1)}^{(y,1+\e^{3/2})}(r^2)  \subseteq  \maxmin_{n;([z,z+d\e],-\e^{3/2})}^{([y,y+\e],0)}(r^2)  
$$
where $z = y - (r+1)\e$.
We may apply Corollary~\ref{c.maxminweight} in order to bound the $\PP$-probability of the two right-hand events. In the first case, we take ${\bf x} = x$, ${\bf y} =  x - (r+1)\e$, ${\bf a} =1$, ${\bf b} =  \lceil 2(r+1) \rceil$, ${\bf t_1} = -\e^{3/2}$, ${\bf t_2} = 0$ and ${\bf r} = r^2$ to find that
$$
 \PP \Big( \neg \, 
\maxmin_{n;([x,x+\e],-\e^{3/2})}^{([y,y+d\e],0)}(r^2)  \Big) \leq  
 (2r+3) \cdot  400 C \exp \big\{ - c_1 2^{-10} r^3 \big\} \, ,
$$
the application being valid provided that 
$$
n \e^{3/2} \geq   10^{29} \vee 2(c/3)^{-18}  \, , \,     3(r + 1)    \leq   6^{-1}  \rsc  n^{1/18} \e^{1/12} \, \textrm{ and } \, r \in \big[  34^{1/2} \, , \, 2 \e^{1/24} n^{1/36} \big] \, .
$$
In fact, we also need to use that $n \e^{3/2}$ is an integer. It was for this reason that we earlier imposed this requirement on $\e > 0$.

In the second case,  we take ${\bf x} = y - (r+1)\e$, ${\bf y} =  y$, ${\bf a} =\lceil 2(r+1) \rceil$, ${\bf b} =  1$,  ${\bf t_1} = 1$, ${\bf t_2} = 1+ \e^{3/2}$ and ${\bf r} = r^2$ to find that
$$
 \PP \Big( \neg \, 
\maxmin_{n;([z,z+d\e],-\e^{3/2})}^{([y,y+\e],0)}(r^2)   \Big) \leq  
 (2r+3) \cdot  400 C \exp \big\{ - c_1 2^{-10} r^3 \big\} \, ,
$$
the conclusion valid under the same hypotheses as in the first case.

We find then that
$$
 \PP \Big( \neg \,
\pwr_{n;(x,-\e^{3/2})}^{([x - (r+1)\e,x+(r+1)\e],0)}(r^2) 
  \Big)  \, + \, \PP \Big( \neg \, \pwr_{n;([y - (r+1)\e,y+(r+1)\e],1)}^{(y,1+\e^{3/2})}(r^2) \Big)
$$
is at most  $2 (2r+3) \cdot  400 C \exp \big\{ - c_1 2^{-10} r^3 \big\}$.
Thus,
\begin{eqnarray*}
  & & \PP \Big(   \maxpoly_{n;(I,0))}^{(J,1)} \geq k  \, , \, \neg \,   \fsc_{n;(I,0)}^{(J,1)}\big(m; \bar{U},\bar{V} ; \e,r  \big)  \Big)  \\
   & \leq & 
   10032 \, C  \exp \big\{ - c_1 2^{-23} r^{3/2}     \big\} 
   \, + \,  22 Cr    \exp \big\{ - 10^{-11} c_1  r^{3/4} \big\}   \\
    &  & \qquad + \, 
8 k^2 
C_k \exp \big\{ - 2^{-3}  k^{-3} c_k r^{3/2} \big\}  \, + \, 2 (2r+3) \cdot  400 C \exp \big\{ - c_1 2^{-10} r^3 \big\}
   \\
    & \leq &   22
 Cr    \exp \big\{ - 10^{-11} c_k  r^{3/4} \big\}  \, + \,   (10032 + 8 + 4000 ) k^2 r C_k \exp \big\{ - 2^{-23}  k^{-3} c_k r^{3/2} \big\}  
   \\
    & \leq &     14062 \, k^2 r C_k \exp \big\{ -  10^{-11}   k^{-3} c_k  r^{3/4} \big\}  \, .
\end{eqnarray*}
where the second inequality used that  the sequence $c_i$ is decreasing,
 as well as $C_k \geq C$, $k \geq 1$ and  $r \geq 1$.
The third used $r \geq 1$, $k \geq 1$, and $C_k \geq C$.
This completes the proof of Lemma~\ref{l.fsc}. \qed

\medskip

\noindent{\bf Proof of Proposition~\ref{p.doubletie}.} 
In this argument, we rigorously implement the resolution of the third challenge of the road map, expressed in outline in the explanation in  Section~\ref{s.thirdchallenge} that led to the two bullet point comments.

Set $I = [x-\e,x+\e]$ and $J = [y - \e,y+\e]$.
When $\maxpoly_{n;(I,0))}^{(J,1)} \geq k$, recall that $\bar{U} \in I^k_\leq$ and  $\bar{V} \in J^k_\leq$
are such that   the collection $\big\{ \rho_{n;(U_i,0)}^{(V_i,1)}: i \in \llbracket 1,k \rrbracket  \big\}$ of $n$-polymers is pairwise horizontally separate.
Set $\phi_i =  \rho_{n;(U_i,0)}^{(V_i,1)}$ for $i \in \llbracket 1,k \rrbracket$.
It is the task of surgery to tie together the $k$ multi-polymer component starting points  $\bar{U}  \times \{ 0 \}$ to what we may call the `lower knot', $\big(x, - \e^{3/2}\big)$, and the $k$ ending points  $\bar{v} \times \{ 1 \}$ to the upper knot $\big(y,1 + \e^{3/2}\big)$.

  Recall the discussion of bouquets and proper weights that followed Proposition~\ref{p.sumweight}.  
In surgery, 
the maximizer 
$\rho_{n,k;(x\bar{\bf 1},- \e^{3/2})}^{{\rm prop};(\bar{U},0)}$ is selected.
(This object is unique, because $\bar{U}$ is measurable with respect to randomness in the region $\R \times [0,\infty)$, 
 so that Lemma~\ref{l.severalpolyunique} applies in view of~(\ref{e.properweightforward}) and~(\ref{e.properweightbackward}).)
 The maximizer is the lower bouquet and we denote it by $\big(\rho_{11},\cdots,\rho_{1k}\big)$.  Forward bouquet regularity $\fbr$ ensures that the bouquet's weight $\properweight_{n,k;(x \bar{\bf 1},- \e^{3/2})}^{(\bar{U},0)}$ satisfies
$$
\bigg\vert \,  \big(\e^{3/2} - n^{-1} \big)^{-1/3}  \properweight_{n,k;(x \bar{\bf 1},- \e^{3/2})}^{(\bar{U},0)}  \, + \,  2^{-1/2}  \big(\e^{3/2} - n^{-1} \big)^{-4/3} \sum_{i=1}^k \big( U_i  - x + 2^{-1} n^{-2/3} \big)^2 \, \bigg\vert \,  \leq \,  r \, .
$$
Since $\vert U_i - x \vert \leq \e$ and $n\e^{3/2} \geq 10^2$, this leads to the simpler
\begin{equation}\label{e.lowersimple}
\e^{-1/2} \, \Big\vert \, \properweight_{n,k;(x \bar{\bf 1},- \e^{3/2})}^{(\bar{U},0)} \, \Big\vert 
 \leq   2^{1/2} (r + k)     \, .
 \end{equation}
The story of the upper bouquet's construction is no different. This bouquet is the maximizer $\rho_{n,k;(\bar{V},1)}^{{\rm prop};(y \bar{\bf 1} ,1 + \e^{3/2})}$, is denoted by $\big(\rho_{21},\cdots,\rho_{2k} \big)$, and has weight that satisfies
\begin{equation}\label{e.uppersimple}
\e^{-1/2} \, \Big\vert \, \properweight_{n,k;(\bar{V},1)}^{(y \bar{\bf 1},1 + \e^{3/2})} \, \Big\vert \leq   2^{1/2} (r + k)  
 \, .
 \end{equation}

Surgery is completed by the construction of a $k$-tuple $\big(\rho_1,\cdots,\rho_k\big)$  of $n$-zigzags each running between $(x,- \e^{3/2})$ and $(y,1 + \e^{3/2})$, where 
$$
\rho_i = \rho_{1i} \circ \phi_i \circ \rho_{2i} \, \, \textrm{for $i \in \llbracket 1,k \rrbracket$} \, .
$$
Given that this $k$-tuple has this set of endpoints, 
the sum of the weights of its elements offers a lower bound on 
$\weight_{n,k;(x \bar{\bf 1},  - \e^{3/2})}^{(y \bar{\bf 1},1   + \e^{3/2})}$. To show, as Proposition~(\ref{p.doubletie}) asserts, that the latter weight 
 is at least  $k \cdot \weight_{n;(x,  - \e^{3/2})}^{(y,1   + \e^{3/2})} \, - \,
 15k r^2  \e^{1/2}$, we will now argue that this is true of $\sum_{i=1}^k \weight(\rho_i)$.

 The polymer weight $\weight_{n;(x,- \e^{3/2})}^{(y,1 + \e^{3/2})}$  may be written as a sum of three terms, $\omega_1$, $\omega_2$ and $\omega_3$, these being the weights of the three $n$-polymers formed by intersecting $\rho_{n;(x,-\e^{3/2})}^{(y,1 + \e^{3/2})}$ 
with the strips $\R \times [ -\e^{3/2} , 0]$, $\R \times  [0,1]$ and $\R \times  [1,1 + \e^{3/2}]$.

Note then that 
\begin{eqnarray}
 & & \bigg\vert \, \sum_{i=1}^k \weight(\rho_i) \, \, - \, k \cdot \weight_{n;(x,-\e^{3/2})}^{(y,1 + \e^{3/2})} \, \bigg\vert \label{e.returnbound} \\
  & \leq & \Big\vert \,   \weight_{n,k;(x \bar{\bf 1},- \e^{3/2})}^{(\bar{u},0)} \, \Big\vert \,  + \,  \Big\vert \,    \weight_{n,k;(\bar{v},1)}^{(y \bar{\bf 1},1 + \e^{3/2})} \, \Big\vert \, + \, \Big\vert \sum_{i=1}^k \big( \weight(\phi_i)  -  \omega_2 \big) \Big\vert  \, +    k \big( \vert \omega_1 \vert + \vert \omega_3 \vert \big) \, , \nonumber
\end{eqnarray}
since for example $\sum_{i=1}^k \weight(\rho_{1i}) =  \weight_{n,k;(x \bar{\bf 1},- \e^{3/2})}^{(\bar{u},0)}$ and  $\sum_{i=1}^k \weight(\rho_{2i}) =   \weight_{n,k;(\bar{v},1)}^{(y \bar{\bf 1},1+ \e^{3/2})}$.
The occurrence of $\pdr_{n;(x,-\e^{3/2})}^{y,1+\e^{3/2}}\big[ \{ 0,1  \}   , r \big]$  implies that 
$$
 \bigg\vert \, \rho_{n;(x,-\e^{3/2})}^{(y,1 + \e^{3/2})}(0) - \Big(  \big( 1 - \tfrac{\e^{3/2}}{1 + 2\e^{3/2}} \big) x  + \tfrac{\e^{3/2}}{1 + 2\e^{3/2}}  y \Big) \, \bigg\vert \, \leq \, r \e(1 + 2\e^{3/2})^{-2/3}
$$
and 
$$
 \bigg\vert \, \rho_{n;(x,-\e^{3/2})}^{(y,1 + \e^{3/2})}(1) - \Big(   \tfrac{\e^{3/2}}{1 + 2\e^{3/2}}  x  +   \big( 1 - \tfrac{\e^{3/2}}{1 + 2\e^{3/2}} \big) y \Big) \, \bigg\vert \, \leq \, r \e(1 + 2\e^{3/2})^{-2/3}
 \, .
$$
Since $\vert y - x \vert \leq \e^{-1/2}$, these imply that 
\begin{equation}\label{e.pairloc}
 \Big\vert \, \rho_{n;(x,-\e^{3/2})}^{(y,1 + \e^{3/2})}(0) -  x \, \Big\vert  \vee  \Big\vert \,  \rho_{n;(x,-\e^{3/2})}^{(y,1 + \e^{3/2})}(1)  -  y \, \Big\vert  \, \leq \, (r+1) \e \, .
\end{equation}

The occurrence of 
 $\pwr_{n;(x,-\e^{3/2})}^{([x - (r+1)\e,x+(r+1)\e],0)}(r^2) \cap \pwr_{n;([y - (r+1)\e,y+(r+1)\e],1)}^{(y,1+\e^{3/2})}(r^2)$, which is one of the favourable surgical conditions,
  then ensures that
\begin{equation}\label{e.omega13}
\max \big\{  \vert \omega_1 \vert , \vert \omega_3 \vert \big\} \leq   \big( r^2 +  2^{-1/2}(r+1)^2 \big) \e^{1/2}  \, .
\end{equation}
It is by considering this bound that we see the reason that the two concerned $\pwr$ events have been chosen to have the form `$(r^2)$': parabolic curvature introduces a term~$2^{-1/2}(r+1)^2 \big) \e^{1/2}$ of the form $\Theta(1)r^2 \e^{1/2}$ into the right-hand side of~(\ref{e.omega13}), and the value $r^2$ has an order which is the highest compatible with this right-hand side maintaining the form  $\Theta(1)r^2 \e^{1/2}$.

For any given $i \in \intint{k}$, the quantity $\weight(\phi_i)  -  \omega_2$ is the difference in weight of two $n$-polymers whose lifetime is  $[0,1]$.
The first polymer begins at a location in  $[x-\e,x+\e]$ and ends at one in $[y-\e,y+\e]$ while we see from~(\ref{e.pairloc}) that the second  begins in  $[x-(r+1)\e,x+(r+1)\e]$ and ends in  $[y-(r+1)\e,y+(r+1)\e]$. 
Thus,~$\lwr_{n;(I^+,0)}^{(J^+,1)} \big( 2(r+1)\e,r \big)$ entails that 
\begin{equation}\label{e.phiomegadifference}
 \vert \weight(\phi_i)  -  \omega_2 \vert\, \leq \, r \big(2(r+1) \e \big)^{1/2} 
\end{equation}
for each $i \in \llbracket 1,k \rrbracket$. 

By revisiting the bound~(\ref{e.returnbound}) equipped with our knowledge of~(\ref{e.omega13}) and~(\ref{e.phiomegadifference}) as well as~(\ref{e.lowersimple}) and~(\ref{e.uppersimple}),
we come to learn that
$$
\bigg\vert \, \sum_{i=1}^k \weight(\rho_i) \, - k \cdot  \weight_{n;(x,-\e^{3/2})}^{(y,1 + \e^{3/2})}  \, \bigg\vert  \leq   \Big( 
  2^{3/2} (r + k)    + k r  2^{1/2} (r+1)^{1/2}  + 2k   \big( r^2 +  2^{-1/2}(r+1)^2 \big) \Big)  \, \e^{1/2}  \, .
$$
Since $r \geq 1$ and $k \geq 1$, the right-hand side is at most
$$
 \big( 2^{3/2}r + 2^{3/2}k +    2kr^{3/2} +  2 ( 1 + 2^{3/2} )  k r^2 \big)  \e^{1/2}  \leq 15 k r^2  \e^{1/2}  \, . 
$$
%15 can read 13
This completes the proof of Proposition~\ref{p.doubletie}. \qed

\subsection{Proof of Theorem~\ref{t.maxpoly}.}
By the scaling principle,
it suffices to consider the case where $t_1 = 0$ and $t_2 = 1$ (so that $\tot = 1$). 
On the event that  $\maxpoly_{n;([x,x+a],0)}^{([y,y+b],1)} \geq \emm$, let   $\bar{U} \in [x,x+a]^\emm_\leq$ and $\bar{V} \in [y,y+b]^\emm_\leq$ be such that the polymer collection $\big\{ \rho_{n;(U_i,t_1)}^{(V_i,t_2)}: i \in \llbracket 1,\emm \rrbracket  \big\}$ is pairwise horizontally separate. At the start of Section~\ref{s.proof}
an explicit means for the selection of a measurable choice of $(\bar{U},\bar{V})$
is offered.

Recall that we denote $h = a \vee b$. Each sequence $\big\{ U_{i+1} - U_i: i \in \llbracket 1, \emm-1 \rrbracket  \big\}$ and $\big\{ V_{i+1} - V_i:  i \in \llbracket 1, \emm-1 \rrbracket  \big\}$ consists of positive terms that sum to at most $h$. Let $\kay \in \intint{\emm-1}$; the value of $\kay$ will be fixed later.
Call an index $i \in \intint{\emm - 1}$ {\em unsuitable}
if  $( U_{i+1} - U_i) \vee ( V_{i+1} - V_i) > 2h(\kay+1)(\emm-1)^{-1}$.
Note that the number of unsuitable indices is less than $(\emm-1)(\kay+1)^{-1}$,
and therefore at most $\lfloor (\emm-1)(\kay+1)^{-1} \rfloor$.
When $\kay \in \intint{\emm-1}$
is chosen so that $5 \kay^2 \leq \emm$, we claim that  there exist $\kay$ consecutive indices $i \in \intint{\emm-1}$ 
such that $( U_{i+1} - U_i) \vee ( V_{i+1} - V_i) \leq 2h(\kay+1)(\emm-1)^{-1}$.
Indeed, were there not such an interval, the number of unsuitable indices would be at least $\lfloor (\emm-1)\kay^{-1} \rfloor$, and this is incompatible with our upper bound, since $(\emm-1)\kay^{-1}$ exceeds $(\emm-1)(\kay+1)^{-1}$ by at least one when  $5 \kay^2 \leq \emm$.

Let  $\Theta \in \intint{\emm -  \kay}$ be chosen so that $\llbracket \Theta, \Theta + \kay - 1 \rrbracket$ is an interval of such not unsuitable indices.
%Suppose that $M \geq L$.
 The polymers $\rho_{n;(U_i,0)}^{(V_i,1)}$ with index $i \in \llbracket \Theta, \Theta + \kay - 1 \rrbracket$
start and end in the  planar intervals 
$$
\big[ U_\Theta, U_\Theta +   2h\kay(\kay+1)(\emm-1)^{-1} \big] \times \{ 0 \} \, \, \, \, 
\textrm{and} \, \, \, \,   \big[ V_\Theta, V_\Theta + 2h\kay(\kay+1)(\emm-1)^{-1}  \big] \times \{ 1 \} \, .
$$ 
These planar intervals are contained in 
$$
\big[ U_{\Theta}^-, U_{\Theta}^- +   (2h\kay+1)(\kay+1)(\emm-1)^{-1}  \big] \times \{ 0 \} \, \, \, \,
\textrm{and} \, \, \, \,  \big[ V^-_{\Theta}, V^-_{\Theta} +  (2h\kay+1)(\kay+1)(\emm-1)^{-1}   \big] \times \{  1 \} \, , 
$$
where we set 
$$
U^-_{\Theta} = (\kay+1)(\emm-1)^{-1} \big\lfloor  (\emm-1)(\kay+1)^{-1} U_\Theta
  \big\rfloor \, \, \,  \textrm{and} \, \, \, 
V^-_{\Theta} = (\kay+1)(\emm-1)^{-1} \big\lfloor  (\emm-1)(\kay+1)^{-1} V_\Theta
  \big\rfloor
$$ 
to be left-displacements onto a~$(\kay+1)(\emm-1)^{-1}$-mesh.
That is, 
$$
\Big\{ \maxpoly_{n;([x,x+a],0)}^{([y,y+b],1)} \geq \emm \Big\} \, \subseteq \, \bigg\{ \maxpoly_{n; 
\big( \big[ U^-_{\Theta}, U^-_{\Theta} +   (2h\kay+1)(\kay+1)(\emm-1)^{-1}  \big] , 0 \big)}^{\big( \big[ V^-_{\Theta}, V^-_{\Theta} +  (2h\kay+1)(\kay+1)(\emm-1)^{-1}   \big] , 1 \big)} \geq \kay \bigg\} \, .
$$
This is a useful moment to recall that, in a remark in Section~\ref{s.closure}, we permitted that the definition of the $\maxpoly$
random variable be respecified so that the concerned polymers are merely pairwise horizontally separate, rather than pairwise disjoint. We actually need to use this form of the definition to obtain the preceding inequality. Indeed, our explicit selection of $(\bar{U},\bar{V})$ means that we are forced to deal with polymers that may be merely horizontally separate.   

Let $X$ denote  the set formed by adding to the mesh points in $(x,x+a]$
the greatest mesh point that is at most $x$, and let $Y$ be the counterpart set where $(y,y+b]$ replaces $(x,x+a]$. The cardinality of $X$ is at most $a(\emm-1)(\kay+1)^{-1} + 2$ and thus, in view of $a \geq 1$ and $\emm \geq \kay^2$, at most $2a\emm \kay^{-1}$. Similarly, $\vert Y \vert \leq 2b\emm \kay^{-1}$.
Since $U^-_\Theta \in U$ and $V^-_\Theta \in V$, we find that
\begin{eqnarray}
 & & \PP \Big( \maxpoly_{n;([x,x+a],0)}^{([y,y+b],1)} \geq \emm  \Big) \nonumber \\
 & \leq &  \sum_{u \in X ,v \in Y}  \PP  \bigg(\maxpoly_{n; 
\big( \big[ u , u +   (2h\kay+1)(\kay+1)(\emm-1)^{-1}  \big] , 0 \big)}^{\big( \big[ v , v  +  (2h\kay+1)(\kay+1)(\emm-1)^{-1}   \big] , 1 \big)} \geq \kay  \bigg) \nonumber \\
  & \leq &   2a\emm \kay^{-1} \cdot 2b\emm \kay^{-1}  \cdot 
   \big( 15h \kay^2/4 \cdot \emm^{-1} \big)^{(\kay^2 - 1)/2} \nonumber \\
 & & \qquad    
    \cdot 
 10^{26\kay^2}  \kay^{43\kay^2/3}  c_\kay^{-8\kay^2/3} C_\kay   \big( \log \emm \kay^{-2} \big)^{4\kay^2}     \exp \big\{ \beta_\kay \big( \log \emm \kay^{-2} \big)^{5/6} \big\} \nonumber \\
& \leq & \emm^{5/2 -\kay^2/2} \cdot h^{(\kay^2 + 3)/2}
10^{27\kay^2}   \kay^{46\kay^2/3}  c_\kay^{-8\kay^2/3} C_\kay   \big( \log \emm \big)^{4\kay^2}     \exp \big\{ \beta_\kay \big( \log \emm \big)^{5/6} \big\} \, . \label{e.expfinalline}
\end{eqnarray}

The last inequality is nothing more than a simplifying of terms that uses $h = a \vee b$ and $\kay \geq 1$. In the inequality that precedes it, 
the summand in the second line is bounded above by using  Theorem~\ref{t.disjtpoly}
with parameter settings ${\bf t_1} = 0$, ${\bf t_2} = 1$, ${\bf x} = u +  \tfrac{1}{2}(2h\kay+1)(\kay+1)(\emm-1)^{-1}$, ${\bf y} = v +  \tfrac{1}{2}(2h\kay+1)(\kay+1)(\emm-1)^{-1}$, ${\bf m} = \kay$ and
$$
{\bm \e} = \tfrac{1}{2}(2h\kay+1)(\kay+1)(\emm-1)^{-1} \, .
$$

In this second inequality, we also make use of the bounds
\begin{equation}\label{e.eulbounds}
\kay^2  \emm^{-1}  \leq {\bm \e} \leq \tfrac{15}{4} h \kay^2  \emm^{-1} \, ,
\end{equation} 
where the latter bound is due to $\kay \geq 2$. 

Of course, the hypotheses of Theorem~\ref{t.disjtpoly}
must be validated with the above parameter settings for the above bound to hold. In fact, we will now set the value of the parameter $\kay \in \N$ in terms of $\emm$, and then justify that, for this choice, these hypotheses are indeed validated.

Recall from Theorem~\ref{t.maxpoly}
that the constant
$\beta$ is set equal to $e \vee \limsup_{i \in \N} \beta_i^{1/i}$ and that Corollary~\ref{c.neargeod.t} implies that $\beta < \infty$. Since $\beta > 1$, we may choose $\kay \in \N$ to be maximal so that $\beta^\kay \leq \big( \log \emm \big)^{1/12}$. 
That is, $\kay = \lfloor \tfrac{1}{12 \log \beta} \log \log \emm \rfloor$. 
Recalling that $\liminf c_i^{1/i}$ is positive, and also that $c_\kay \leq c$, we see that the condition~(\ref{e.epsilonbound}) on ${\bm \e}$ in Theorem~\ref{t.disjtpoly} is met provided that $\emm \geq \emm_0$, where  $\emm_0 \in \N$ is a certain constant; a suitable choice of $\emm_0$ also ensures that $\kay \geq 2$.

Also note that $\vert {\bf x} - {\bf y} \vert = \vert u - v \vert \leq \vert x - y \vert + 2h$ because
 $\vert u - v \vert
 \leq \vert x - y \vert + h + (\kay+1)(\emm-1)^{-1}$ while  $h \geq 1 \geq  (\kay+1)(\emm-1)^{-1}$ due to $\emm \geq \kay + 2$.

Since $\kay \leq (12)^{-1} \log \log \emm$ (due to $\beta \geq e$), $K_0 \geq 1$, ${\bm \e} \geq \emm^{-1}$ and  $\liminf c_i^{1/i} > 0$, the bound~(\ref{e.nlowerbound}) is verified when $n$ is at least  $2(K_0)^{(12)^{-2} (\log \log \emm)^2} \big( \log \emm \big)^{K_0}$ and 
$$
  \max \bigg\{ \,
     10^{584}     \big( (12)^{-1} \log \log \emm \big)^{240}  \rsc^{-36} \emm^{225}  \max \big\{   1  \, , \,   (\vert y - x \vert + 2h)^{36}   \big\} \, , \, 
  a_0^{-9} ( \vert y - x \vert  + 2h ) \,  \bigg\} \, ,
$$
where we used $\emm \geq \emm_0$ and adjusted the value of $\emm_0$ if need be.
Finally, the condition 
$$
\vert {\bf y} - {\bf x} \vert ({\bf \tot})^{-2/3} \leq {\bm \e}^{-1/2} \big( \log {\bm \e}^{-1} \big)^{-2/3} \cdot 10^{-8} c_\kay^{2/3} \kay^{-10/3}
$$
 is verified provided that
$\emm \geq  \emm_0 \vee \big( \vert x - y \vert + 2h \big)^3$ in light of  ${\bm \e} \leq \emm^{-1 + o(1)}$ (with an increase if necessary in the value of $\emm_0$).
This lower bound on $n$ is implied by Theorem~\ref{t.maxpoly}'s hypothesis~(\ref{e.nmaxpoly}) due to $\log \log \emm \leq \emm$.
%and $(12)^{240} = 1.008 \cdots \times 10^{259}$.

Given the choice of $\kay$, we see that the expression~(\ref{e.expfinalline}) is at most
$$
 \emm^{-   (145)^{-1} ( \log \beta)^{-2} (0 \vee \log \log \emm)^2} h^{(\log \beta)^{-2} (0 \vee \log \log \emm)^2/{288}  + 3/2} \conseqmac_\emm \, ,
$$
%0 \vee added in second instance
where $\big\{ \conseqmac_i: i \in \N \big\}$ is a sequence of positive constants such that $\sup_{i \in \N} \conseqmac_i \exp \big\{ - 2 (\log i)^{11/12} \big\}$ is finite.

This proves Theorem~\ref{t.maxpoly} when $t_0 = 0$ and $t_1 = 1$; as we noted at the outset, there is no loss of generality in considering this case.  \qed

\section{Polymer fluctuation: proving Theorem~\ref{t.polyfluc}.}\label{s.polyfluc}

In this section, we will derive the polymer fluctuation Theorem~\ref{t.polyfluc}. The result asserts that any of the polymers that cross between unit-order length intervals separated at unit-order times deviates by a distance $r$ from the line segment that interpolates its endpoints with probability at most $\exp \big\{ - O(1) r^{3/4} \big\}$, uniformly in high choices of the scaling parameter $n$. 
We may easily reduce to the case where these unit-order intervals are instead singleton sets, however, by a simple application of polymer ordering. In Proposition~\ref{p.polyfluc}, this reduced version of the theorem has been stated. After noting how this proposition implies Theorem~\ref{t.polyfluc}, we will turn to discuss the ideas of the proof of the proposition, and then give the proof itself.

\noindent{\bf Proof of Theorem~\ref{t.polyfluc}.}
%For $(t_1,t_2) \in \R^2_<$, $a \in (0,1)$ and 
For $u,v \in \R$, consider the random variable
$$
 X_{n;(u,t_1)}^{(v,t_2)}(a) =   t_{1,2}^{-2/3} \big(  a \wedge (1-a) \big)^{-2/3}  \Big( \, \rho_{n;(u,t_1)}^{(v,t_2)} \big(  (1-a) t_1 +  a t_2 \big) - \ell_{(u,t_1)}^{(v,t_2)}  \big(  (1-a) t_1 +  a t_2 \big) \, \Big) \, .
 $$
 Set $\Rmac =   t_{1,2}^{2/3} \big(  a \wedge (1-a) \big)^{2/3} r$,
and recall from the theorem's statement that we define  $I = [x,x+\Rmac]$ and $J = [y,y+\Rmac]$.
When the event $\neg \, \pdr_{n;(I,t_1)}^{(J,t_2)}\big(a,2r\big)$ 
occurs, we may choose  $(U,V) \in I \times J$ such that  $\big\vert X_{n;(U,t_1)}^{(V,t_2)}(a) \big\vert \geq 2\Rmac$. It might seem reassuring to know that this selection may be made measurably, and with a modicum of effort we might show this, but in fact this information is not needed.

Recall the polymer sandwich Lemma~\ref{l.sandwich}.
By applying this result with parameter settings ${\bf x_1} = x$, ${\bf x_2} = x + \Rmac$, ${\bf y_1} = y$, ${\bf y_2} = y + \Rmac$, we readily find that
at least one of  $\big\vert X_{n;(x,t_1)}^{(y,t_2)}(a) \big\vert$
and  $\big\vert  X_{n;(x+\Rmac,t_1)}^{(y+\Rmac,t_2)}(a) \big\vert$ is at least $\big\vert  X_{n;(U,t_1)}^{(V,t_2)}(a) \big\vert  - \Rmac$.
Thus, 
\begin{eqnarray*}
 & & \PP \Big( \neg \, \pdr_{n;(I,t_1)}^{(J,t_2)}\big(a,2r\big) \Big) \\
 & \leq &  
\PP \Big( \neg \, \pdr_{n;(x,t_1)}^{(y,t_2)}\big(a,r\big) \, \Big) \,  + \, 
\PP \Big( \neg \, \pdr_{n;(x+\Rmac,t_1)}^{(y+\Rmac,t_2)}\big(a,r\big) \, \Big) \, .
\end{eqnarray*}
Two applications of Proposition~\ref{p.polyfluc},
with $({\bf x},{\bf y})$ equal to $(x,y)$ and $(x+\Rmac,y+\Rmac)$, complete 
 the proof. \qed

We now overview the ideas of the proof of Proposition~\ref{p.polyfluc}. Recall from Section~\ref{s.roadmap} the role of the powers of one-half and one-third.
The one-half power law is articulated in Theorem~\ref{t.differenceweight}: polymer weight has a H\"older exponent of one-half in response to horizontal displacement of endpoints.
The one-third power law is articulated by the scaling principle alongside the regular sequence conditions $\rmreg(2)$ and $\rmreg(3)$: a polymer of lifetime $t$ whose endpoints differ by $r t^{2/3}$
has weight $t^{1/3}\big( U - 2^{-1/2} r^2 \big)$; here, $U$ is a unit-order, random, quantity. 

These principles come into conflict unless the two-thirds principle for polymer geometry (which we are trying to prove) also obtains. 

To explain this concept in more detail, it is useful to begin by fixing a little notation. We set $Z = \rho_{n;(x,t_1)}^{(y,t_2)} \big( (1-a)t_1 + a t_2 \big)$ and 
$z = \ell_{(x,t_1)}^{(y,t_2)}  \big(  (1-a) t_1 +  a t_2 \big)$. Recall that $a$ is supposed to be close to one, and that we seek to show that $Z$ is close  to $z$.
Using notation already seen in the proof of Theorem~\ref{t.polyfluc}, we write $\Rmac = t_{1,2}^{2/3} (1-a)^{2/3} r$. 
The factor of $t_{1,2}^{2/3} (1-a)^{2/3}$ is a suitable scale for judging polymer fluctuation between the pairs of times $(1-a)t_1 + a t_2$ and $t_2$,
and our task is to show that it is a rare event, with a decaying probability expressed in terms of $r$, that $Z \leq z - \Rmac$.

\begin{figure}[ht]
\begin{center}
\includegraphics[height=8cm]{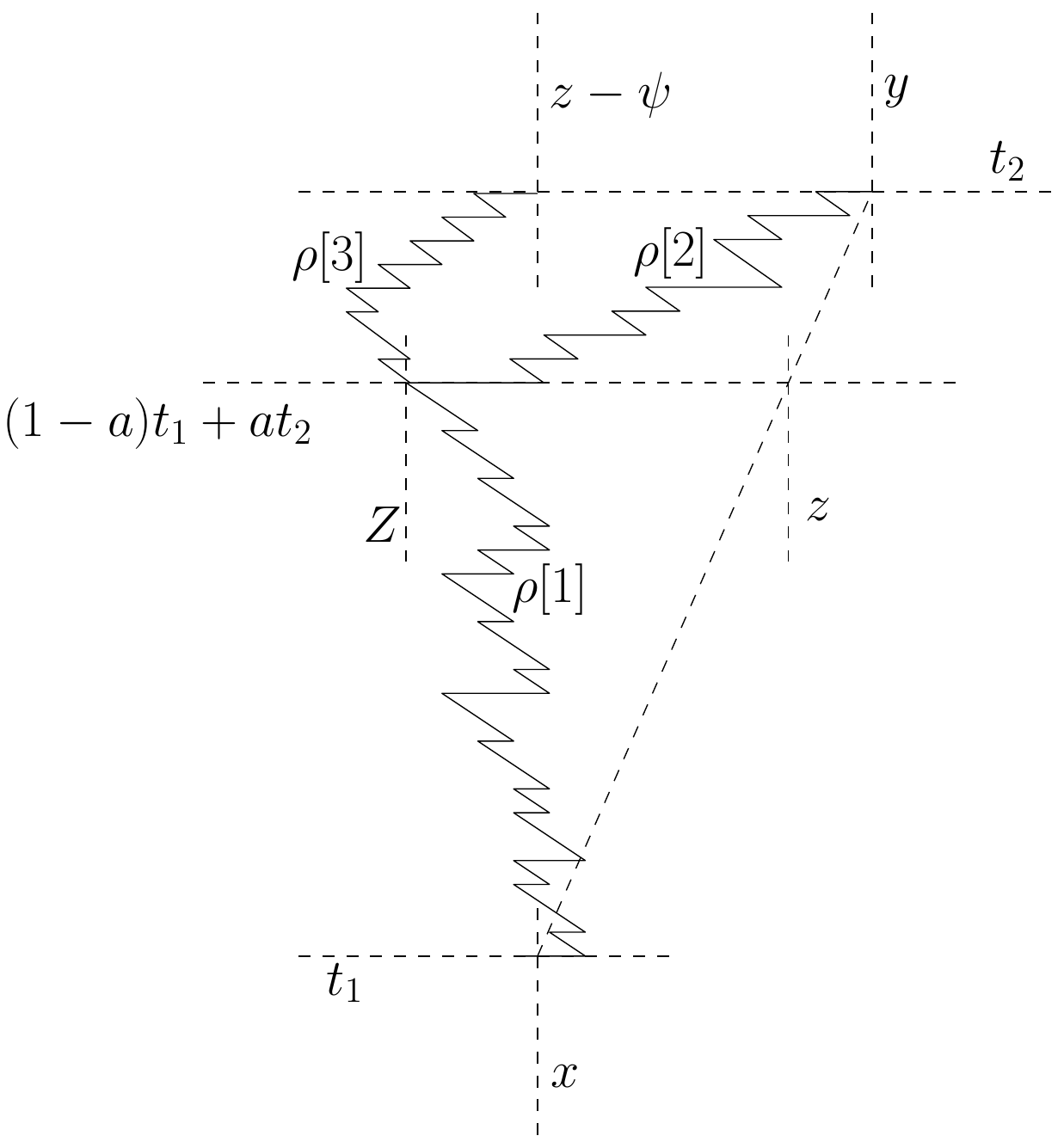}
\caption{The polymer $\rho = \rho_{n;(x,t_1)}^{(y,t_2)}$ is split at time $(1-a)t_1 + at_2$
and the resulting subpaths are labelled $\rho[1]$ and $\rho[2]$. The polymer $\rho[3]$ has a less hectic journey than does $\rho[2]$ during their shared lifetime. The concatenation $\rho[1] \circ \rho[3]$
will be labelled~$\hat\rho$ in the upcoming proof of Proposition~\ref{p.polyfluc} and shown in cases such as that depicted to have a significantly greater weight than $\rho$.}
\label{f.threerho}
\end{center}
\end{figure}

 Abbreviate $\rho = \rho_{n;(x,t_1)}^{(y,t_2)}$ and consult Figure~\ref{f.threerho}. The polymer $\rho$ may be split into two pieces, $\rho[1]$ and $\rho[2]$, by cutting it at the point $(Z,(1-a)t_1 + a t_2)$.
If the big fluctuation event $Z \leq z - \Rmac$ occurs, then $\rho[2]$ makes a rather sudden deviation, and its weight is dictated by parabolic curvature to be of order $-r^2 (1-a)^{1/3} \tot^{1/3}$.
We may consider the polymer $\rho[3]$ whose endpoints are the cut location $(Z,(1-a)t_1 + a t_2)$ and $(z - \Rmac,t_2)$. The latter endpoint being  $(z - \Rmac,t_2)$ means that $\rho[3]$
is not making the sudden deviation that $\rho[2]$ does, so that its weight is dictated by local randomness to be a random unit-order multiple of  $(1-a)^{1/3} \tot^{1/3}$.
Consider now the polymer weight profile $\weight_{n;(x,t_1)}^{(v,t_2)}$ as a function of $v$. Between the locations $v = y$ and $v = z - \Rmac$, which is an order of $\Rmac$ to the left of~$y$, the weight profile has risen by an order of least $-r^2 (1-a)^{1/3} \tot^{1/3}$, because the weight profile at $v = z - \Rmac$ is at least the sum of the weights of $\rho[1]$ and $\rho[3]$, while at $v = z$, it equals the sum of the weights of $\rho[1]$ and $\rho[2]$.

Recalling that $\Rmac = t_{1,2}^{2/3} (1-a)^{2/3} r$, we see that the weight profile is experiencing a change where a factor $r^{3/2}$ multiplies the square-root $\Theta\big(\Rmac^{1/2}\big)$ of the horizontal displacement in endpoint location. When $r$ is large, this is incompatible with the one-half power law principle recalled above.

This idea is important, but it is not adequate to prove Proposition~\ref{p.polyfluc}. To understand why the idea is not the end of the story, examine the form of
Corollary~\ref{c.ordweight},
% Theorem~\ref{t.differenceweight}, 
which articulates the one-half power law for the polymer weight profile. The quantity $\e$, whose role is to measure horizontal displacement, must be at most $2^{-4}$.
Indeed, the one-half power law does not govern the weight profile beyond a unit scale. This means that our idea will work only in a case that we may call {\em near}: when the quantity $r$ is large, but not so large that the concerned endpoint displacement goes above a unit scale.

When $r$ is large enough that the one-half power law estimate cannot be applied, we will note that the polymer weight profile change is unlikely for other reasons. Indeed, since $r$ is so  large, 
the observed change in the weight profile is very large. In fact, this change is so large that it either forces the profile at value $v = y$ to be abnormally low or the profile at $v = z - \Rmac$ to be abnormally high.
The improbability of these outcomes may then be gauged using the regular ensemble conditions $\rmreg(2)$ and $\rmreg(3)$.

In fact, even this is not the end of the story. These two regular ensemble conditions capture the parabolic curvature of the weight profile. But this curvature breaks down on an extremely large scale, as witnessed by the upper bounds on $\vert z \vert$ needed in  $\rmreg(2)$ and $\rmreg(3)$. This second argument thus works when $r$ is very large, but not extremely large (where the latter term involves an $n$ dependence). For this reason, this second argument will be said to apply in the {\em middle-distant} case.

When $r$ is extremely high, the linear (but not parabolic) collapse near infinity of the polymer weight profile expressed by Proposition~\ref{p.mega}(4) will replace the use of $\rmreg(3)$.
This third argument applies in what we will call the {\em far} case.

This then will be the structure of the argument, with three different estimates treating the near, middle-distant and far cases, in increasing order of the value of $r$.

\medskip

\noindent{\bf Proof of Proposition~\ref{p.polyfluc}.}
Our task is to understand the fluctuation  of the polymer $\rho_{n;(x,t_1)}^{(y,t_2)}$ at time $(1-a)t_1 + a t_2$.
Note first that, since the endpoints are given, we may harmlessly suppose that this polymer is unique, because, as we discussed in Section~\ref{s.brlpp}, the complementary event has zero probability.
Note also that we are working in a case where $a$ is close to one; specifically, $a \geq 1/2$. (The remark following the proposition treats the opposing case.)

For $K$ a real interval, we write $\fluc_{n;(x,t_1)}^{(y,t_2)}\big(a;K \big)$ for the event~\hfff{fluc} that $X_{n;(x,t_1)}^{(y,t_2)}(a)  \in K$. We omit the parameters that are fixed from this notation, writing for example
$\fluc[s,t] =  \fluc_{n;(x,t_1)}^{(y,t_2)}\big(a;[s,t]\big)$ for $(s,t) \in \R^2_<$.

 We will argue that, under the proposition's hypotheses,
\begin{equation}\label{e.fminusinfinityr}
 \PP \Big( \fluc (-\infty,-r] \Big) \, \leq \,  
22  C r    \exp \big\{ - 10^{-11} c_1  r^{3/4} \big\}      \,  ,
\end{equation}
and that the same bound holds on 
 $\PP \big( \fluc [r,\infty) \big)$. Indeed, since  the event  $\neg \, \pdr_{n;(x,t_1)}^{(y,t_2)}\big(a,r\big)$ equals the union of $\fluc (-\infty,-r)$ and $\fluc (r,\infty)$,
 these bounds suffice to prove Proposition~\ref{p.polyfluc}.
 
To derive the two bounds,  it suffices to prove~(\ref{e.fminusinfinityr}). 
The two parts of the next proposition and the following lemma are the three principal tools: they treat the near, middle-distant and far cases. The reader should bear in mind that, in these assertions and later ones, it is supposed, as it was in  Proposition~\ref{p.polyfluc}, that $(n,t_1,t_2)$ is  a compatible triple, and that $a\tot \in n^{-1} \Z$. 

\begin{proposition}\label{p.feventfirstbound}
 Set $\hat{n} = n (1-a) \tot$.
Suppose that
\begin{equation}\label{e.firstbound.h1}
\hat{n}   \geq   10^{32} c^{-18} \, \, , \, \, 
 r \in \Big[ 8 \vee 87(1-a)^{1/3} \tot^{-2/3} \vert x - y \vert  \, , \,  3 \hat{n}^{1/36}  \Big] 
 \end{equation}
 and 
\begin{equation}\label{e.firstbound.h2}
   \big\vert x - y  \big\vert  \tot^{-2/3} \leq 2^{-3} 3^{-1} \rsc  (n \tot)^{1/18} \, .
 \end{equation}
\begin{enumerate}
\item %Set $n' = n\tot$. 
Suppose in addition  that    
$$
r \in \big[ 2800  \, , \,  \tfrac{1}{18} (1-a)^{-2/3}  \big]  
$$
and $a \geq 1 - 10^{-3}$. Then
$$
  \PP \Big( \fluc \big[-(r+1),-r\big] \Big) \, \leq \,  
   (10178 + 3r) \, C   \exp \big\{ - c_1 2^{-26} 3^{-3} 5^{3/2}  r^{9/4} \big\} 
  % was (25210 + 3r) \, C   \exp \big\{ - c_1 2^{-34} 3^{-3} 5^{3/2}  r^{9/4} \big\}  \, .
$$
\item Now suppose in addition to~(\ref{e.firstbound.h1}) and~(\ref{e.firstbound.h2})  that
$r  \geq 6(1-a)^{-1/2}$ and $a \geq 1 - 10^{-2}$.
Then
$$
 \PP \Big(  \fluc \big[-(r+1),-r\big] \Big)
  \, \leq \,    
 \big( 148 C +  9r/4 \big)  \exp \big\{ - c_1  (1-a)^{1/2}   3^{-3/2} 5^{3/2} 2^{-43/4} r^3 \big\} \, .
$$
\end{enumerate}
\end{proposition}
\begin{lemma}\label{l.feventhigh}
Suppose that $a \in (1/2,1)$, and that  
 $n \tot$ is at least the maximum of 
$$
(1-a)^{-25} \tot^{-24} \vert x - y \vert^{36} \, \, , \, \,  
c^{-9} \tot^{-6} \vert x - y \vert^9 \, \, \textrm{ and } \, \,   10^5 (1-a)^{-17} c^{-12}
 \, . 
$$
Then
$$
   \PP \Big(  \fluc \big(-\infty,-2\hat{n}^{1/36} \big] \Big) \, \leq  \, 13C \exp \Big\{ - 2^{-15/4} c (1-a)^{25/12} \tot^{1/12}       n^{1/12} \Big\}    \, .
$$
\end{lemma}

We now derive~(\ref{e.fminusinfinityr}) using this proposition and lemma. 
 The interval $(-\infty,-r]$ may be partitioned into a {\em far} interval $(-\infty,-2\hat{n}^{1/36}]$
 and a collection of disjoint unit intervals of the form $(-t-1,-t]$
 where  $t - r$ varies over a finite initial sequence of integers; i.e., $t = -r$, then $t=-r-1$, and so on. (In fact, one of the unit intervals may overlap with the far interval.) The unit intervals to the left of $-\tfrac{1}{18} (1-a)^{-2/3}$ will be called {\em middle-distant}, and the remaining ones, closer to the origin, will be called {\em near}. 
  To each of these intervals~$I$ corresponds an event $\fluc \, I$.
 In the case of the far interval, the event's probability is bounded above by  Lemma~\ref{l.feventhigh}. Proposition~\ref{p.feventfirstbound}(2) provides the bound for the middle-distant intervals, and Proposition~\ref{p.feventfirstbound}(1)  for the near intervals. Note that Proposition~\ref{p.feventfirstbound}(2) is applicable to the middle-distant intervals because the condition  $\tfrac{1}{18} (1-a)^{-2/3} \geq  5(1-a)^{-1/2}$, i.e., $a \geq 1 - 90^{-6}$, is implied by Proposition~\ref{p.polyfluc}'s hypothesis that  $a \geq 1 - 10^{-11} c_1^2$  alongside $c_1 \leq 1/8$.
 
 We are summing the right-hand sides of  both parts of  Proposition~\ref{p.feventfirstbound}. The next two lemmas provide an upper bound on these sums.  The  calculational proofs  appear in the online Appendix~$D$.
%arxiv/submit
 
\begin{lemma}\label{l.calcone} 
Suppose that $r \geq  10^7 c_1^{-4/5}$. 
%was $r \geq  10^9 c_1^{-4/5}$. 
We have that
$$
 \sum   (10178 + 3r) \, C   \exp \big\{ - c_1 2^{-26} 3^{-3} 5^{3/2}  s^{9/4} \big\} 
 \leq  8  C r   \exp \big\{ - c_1 10^{-8}  r^{9/4} \big\}
  \, ,
$$
% was $$
% \sum   (25210 + 3r) \, C   \exp \big\{ - c_1 2^{-34} 3^{-3} 5^{3/2}  s^{9/4} \big\} 
% \leq  8  C r   \exp \big\{ - c_1 10^{-11}  r^{9/4} \big\}
%  \, ,
%$$
where the sum  is taken over values of $s \in \R$ that satisfy $s \in r + \N$.
\end{lemma}

 \begin{lemma}\label{l.calctwo}
 Suppose now that $r \geq  \tfrac{1}{18} (1-a)^{-2/3} \vee 15 C^{1/2}$ as well as $a \geq 1 - 10^{-11} c_1^2$. Setting $\alpha = (5/3)^{3/2} 2^{-43/4}$, we have that
 \begin{eqnarray*}
 & &  \sum 
 \big( 148 C +  9s/4 \big)  \exp \big\{ - c_1  (1-a)^{3/2}  (5/3)^{3/2} \cdot 2^{-43/4} s^3 \big\} \\
 &  \leq &   
 \exp \big\{ - 10^{-2} c_1 \alpha r^{9/4} \big\}
 \, ,
 \end{eqnarray*}
 where the sum is again over $s \in r + \N$.
 \end{lemma}

The bound~(\ref{e.fminusinfinityr}) may now be obtained using Lemma~\ref{l.feventhigh} alongside the two preceding lemmas.
Indeed, we obtain an upper bound on~$\PP \big( \fluc (-\infty,-r] \big)$ of the form  
$$
9 C r   \exp \big\{ -  10^{-8} c_1   r^{3/4} \big\} \, + \,  13C \exp \Big\{ - 2^{-15/4} c_1 (1-a)^{25/12} \tot^{1/12}       n^{1/12} \Big\} 
$$
because this expression is,  in view of $C \geq 1$, $r \geq 1$ and $c_1 \leq c$, an upper bound on the sum of the right-hand sides in Lemmas~\ref{l.feventhigh},~\ref{l.calcone} and~\ref{l.calctwo}.
That this expression is under the hypotheses of Proposition~\ref{p.polyfluc} at most the right-hand side of~(\ref{e.fminusinfinityr}) is a calculational matter.
This assertion is stated as Lemma~$D.3$
%\ref{l.derive} 
in the online Appendix~$D$ and derived there.
%arxiv/submit

This completes the proof of~(\ref{e.fminusinfinityr}), and  Proposition~\ref{p.polyfluc},
subject to confirming  Proposition~\ref{p.feventfirstbound} and Lemma~\ref{l.feventhigh}.

We now prepare to prove this proposition and lemma. Recall throughout that we are supposing that $a \geq 1/2$.
 
Recalling notation that was used in the overview that preceded the present proof, we write $\Rmac = t_{1,2}^{2/3} (1-a)^{2/3} r$. We also
set $Z = \rho_{n;(x,t_1)}^{(y,t_2)} \big( (1-a)t_1 + a t_2 \big)$ and 
$z = \ell_{(x,t_1)}^{(y,t_2)}  \big(  (1-a) t_1 +  a t_2 \big)$.

%We claim that, when   $r \geq 8(1-a)^{1/3} \tot^{-2/3} \vert x - y \vert$,
%\begin{equation}\label{e.zfst}
% \fluc [-t,-s] \subseteq  \Big\{ \, Z \in y \, - \, \tot^{2/3} (1-a)^{2/3} \cdot \big[ s - r/8 , t + r/8 \big] \, \Big\} \, ,
%\end{equation}
%for any $t > s > 0$.
%To verify this, note that 
%$Z \in z  - t_{1,2}^{2/3} (1-a)^{2/3} \cdot [s,t]$ and that
%$z = y + (1-a)(x - y)$.
%The hypothesised lower bound on $r$  
%implies that $(1-a) \vert x - y \vert \leq \Rmac/8$, confirming~(\ref{e.zfst}).

 Abbreviate $\rho = \rho_{n;(x,t_1)}^{(y,t_2)}$. We now split $\rho$ at the point $(Z,(1-a)t_1 + a t_2)$ according to polymer splitting as described in Section~\ref{s.split}.
Indeed, we may specify $\rho[1]$ and $\rho[2]$ to be the polymers $\rho_{n;(x,t_1)}^{(Z,(1-a)t_1 + a t_2)}$ and $\rho_{n;(Z, (1-a) t_1 +  a t_2)}^{(y,t_2)}$. Thus, $\rho = \rho[1] \circ \rho[2]$. 
 Let $\rho[3]$ denote the polymer $\rho_{n;(Z, (1-a) t_1 +  a t_2)}^{(z - \Rmac,t_2)}$. Let $\hat\rho = \rho[1] \circ \rho[3]$. Thus,  $\rho$ and~$\hat\rho$ 
are two polymers that begin from $(x,t_1)$ and follow a shared course until time $(1-a) t_1 +  a t_2$. They continue on possibly separate trajectories until a shared ending time $t_2$, with $\rho$ then reaching~$y$ and~$\hat\rho$ reaching $z - \Rmac = y + (1-a)(x-y) - \rmac (1-a)^{2/3} \tot^{2/3}$.  

The next lemma is a tool to prove the near and middle-distant  Proposition~\ref{p.feventfirstbound}.  It asserts that, when $\fluc \, [-(r+1),-r]$ occurs, the polymer $\rho[3]$ 
is unlikely to have a significantly negative weight, while $\rho[2]$ is unlikely to have a weight much exceeding a certain constant multiple of $- r^2$; thus, when $r$ is chosen so that either the near and middle-distant case obtains, $\rho[3]$
is typically much heavier than~$\rho[2]$ on the event $\fluc \, [-(r+1),-r]$.

\begin{lemma}\label{l.combine}
Suppose that the hypotheses ~(\ref{e.firstbound.h1}) and~(\ref{e.firstbound.h2}) hold.
\begin{enumerate}
\item We have that
$$
 \PP \Big( \, \fluc \, [-(r+1),-r] \, , \, t_{1,2}^{-1/3} (1-a)^{-1/3} \weight\big(\rho[3]\big) \leq - 2^{-1/2} - s \Big) \leq 140 C \, \exp \big\{ - c_1 s^{3/2} \big\} \, .
 $$
 \item We also have
 \begin{eqnarray*}
 & &  \PP \Big(   \, \fluc \, [-(r+1),-r]  \, , \,  t_{1,2}^{-1/3} (1-a)^{-1/3} \weight\big(\rho[2]\big) \geq - 2^{-1/2} (7r/8)^2 + s \Big)  \\
  &  \leq &  \big(3r/8 + 1 \big) \cdot 6 C  \exp \big\{ - 2^{-11/2} c s^{3/2} \big\}  \, .
 \end{eqnarray*}
  \end{enumerate}
 \end{lemma}
\noindent{\bf Proof.}
 Note that 
 $$
 \weight\big( \rho[3] \big) = \weight_{n;\big(Z, (1-a) t_1 +  a t_2 \big)}^{(z - \Rmac,t_2)} = \mc{L}_{n;(1-a) t_1 +  a t_2}^{\downarrow;(z-\Rmac,t_2)}\big( 1, Z \big) =  (1-a)^{1/3} \tot^{1/3} \scaledle_{n;(1-a) t_1 +  a t_2}^{\downarrow;(z-\Rmac,t_2)}\big( 1, V \big)
 $$ 
 where $V = (1-a)^{-2/3} \tot^{-2/3} \big( Z- (z-\Rmac) \big)$. Note that the event  $\fluc \, [-(r+1),-r]$ is characterized by $- (r+1) \leq (1-a)^{-2/3} \tot^{-2/3} \big( Z- z) \leq -r$ and thus also by the condition that~$V \in [-1,0]$.
Thus,
\begin{eqnarray}
 & & \PP \Big( \, \fluc \, [-(r+1),-r] \, , \, t_{1,2}^{-1/3} (1-a)^{-1/3} \weight\big(\rho[3]\big) \leq - 2^{-1/2} - s \Big) \nonumber \\
 & = &
\PP \Big( \, \fluc \, [-(r+1),-r] \, , \, \scaledle_{n;(1-a) t_1 +  a t_2}^{\downarrow;(z-\Rmac,t_2)}\big( 1, V \big)
 \leq - 2^{-1/2} - s \Big) \nonumber \\
 & \leq &
\PP \Big( \, \fluc \, [-(r+1),-r] \, , \,  \scaledle_{n;(1-a) t_1 +  a t_2}^{\downarrow;(z-\Rmac,t_2)}\big( 1, V \big)
 + 2^{-1/2} V^2 \leq - s \Big)  \nonumber \\
 & \leq &
\PP \bigg(  \inf_{-1 \leq v \leq 0} \Big( \scaledle_{n;(1-a) t_1 +  a t_2}^{\downarrow;(z-\Rmac,t_2)}\big( 1, v \big)
 + 2^{-1/2} v^2  \Big) \leq - s  \bigg)  \nonumber \\
 & \leq &  \Big( 1/2 \vee 5 \vee  5^{1/2} (3 - 2^{3/2})^{-1}  \Big) \, 10 C \, \exp \big\{ - c_1 s^{3/2} \big\} \leq  140 C \, \exp \big\{ - c_1 s^{3/2} \big\} \, ,  \nonumber
\end{eqnarray}
%  5^{1/2} (3 - 2^{3/2})^{-1}  = 13.032759
so that we obtain Lemma~\ref{l.combine}(1).
The final inequality was obtained by applying  
Proposition~\ref{p.mega}(2) 
%Proposition~\ref{p.strongothercurves} 
with the parameter choices ${\bf k} = 1$, ${\bf y} = -1/2$ , ${\bf t} = 1/2$ and ${\bf r} = s$ to the $(c,C)$-regular ensemble $\scaledle_{n;(1-a) t_1 +  a t_2}^{\downarrow;(z-r,t_2)}$. The ensemble has $\hat{n} + 1$ curves, where recall that we denote $\hat{n} = n (1-a) \tot$. As such, this application of 
Proposition~\ref{p.mega}(2) may be made provided that 
$$
1/2 \leq c/2 \cdot \hat{n}^{1/18} \, , \,  
 1/2 \leq \hat{n}^{1/18}  \, , \, s  \in \big[ 2^{3/2} \, , \, 2 \hat{n}^{1/18} \big] \, \, \textrm{ and } \, \, \hat{n}  \geq 1 \vee  (c/3)^{-18} \vee  6^{36} \, .
 $$

We have the equality
\begin{equation}\label{e.vprime.ensemble}
 \weight\big( \rho[2] \big) =  (1-a)^{1/3} \tot^{1/3} \scaledle_{n;(1-a) t_1 +  a t_2}^{\downarrow;(y,t_2)}\big( 1, V' \big)
\end{equation}
 where we set $V'$ equal $(1-a)^{-2/3} \tot^{-2/3} ( Z- y )$. 
 Indeed, the polymer weight $\weight\big( \rho[2] \big)$ equals  $\weight_{n;\big(Z, (1-a) t_1 +  a t_2 \big)}^{(y,t_2)} = \mc{L}_{n;(1-a) t_1 +  a t_2}^{\downarrow;(y,t_2)}\big( 1, Z \big)$, whence we obtain~(\ref{e.vprime.ensemble}). 

We claim that, when  $r \geq 2 \vee 8(1-a)^{1/3} \tot^{-2/3} \vert x - y \vert$,
\begin{equation}\label{e.vprimeclaim}
 \fluc \, [-(r+1),-r]  \subseteq \Big\{ \, 
  V' \in r \cdot [-13/8,-7/8] \,  \Big\} \, . 
\end{equation} 
 To see this, note that, when $\fluc [-(r+1),-r]$ occurs,
$$
 Z  \in  z - (1-a)^{2/3} \tot^{2/3} \cdot [r,r+1]
 = y + (1-a)(x-y) -  (1-a)^{2/3} \tot^{2/3} \cdot [r,r+1]
$$ 
and, in light of our lower bound on $r$,
$$
 Z \in y -  (1-a)^{2/3} \tot^{2/3} \cdot [r -r/8, r + r/8 + 1] \subseteq y \, - \, \Rmac \cdot [7/8,13/8]
$$
where we used $r \geq 2$ in the form $r+1 \leq 3r/2$.
Thus, we confirm~(\ref{e.vprimeclaim}).

We thus obtain Lemma~\ref{l.combine}(2):
\begin{eqnarray}
 & & \PP \Big(   \, \fluc \, [-(r+1),-r]  \, , \,  t_{1,2}^{-1/3} (1-a)^{-1/3} \weight\big(\rho[2]\big) \geq - 2^{-1/2} (7r/8)^2 + s \Big)  \nonumber \\
 & = &
\PP \Big(  \, \fluc \, [-(r+1),-r]  \, , \,  \scaledle_{n;(1-a) t_1 +  a t_2}^{\downarrow;(y,t_2)}\big( 1, V' \big)
 \geq - 2^{-1/2} (7r/8)^2 + s \Big) \nonumber \\
 & \leq &
\PP \Big(   \, \fluc \, [-(r+1),-r]  \, , \,  \scaledle_{n;(1-a) t_1 +  a t_2}^{\downarrow;(y,t_2)}\big( 1, V' \big)
 + 2^{-1/2} (V')^2 \geq  s \Big)  \nonumber \\
  & \leq & 
\PP \bigg(   \sup_{-13r/8 \leq v \leq -7r/8} \Big( \scaledle_{n;(1-a) t_1 +  a t_2}^{\downarrow;(y,t_2)}\big( 1, v  \big)
 + 2^{-1/2} v^2 \Big) \geq  s \bigg)  \nonumber \\
 & \leq & \big(3r/8 + 1 \big) \cdot 6 C  \exp \big\{ - 2^{-11/2} c s^{3/2} \big\} 
 \, .  \nonumber
\end{eqnarray}
Here, the final  inequality is due to an application of 
Proposition~\ref{p.mega}(3)
%Proposition~\ref{p.nobigmax.gen}
 to the $(c,C)$-regular ensemble $\scaledle_{n;(1-a) t_1 +  a t_2}^{\downarrow;(y,t_2)}$ which shares a curve cardinality of $\hat{n} + 1$ with the ensemble that we considered a few moments ago. Our parameter choice for the application is
 ${\bf y} = - 5r/4$, ${\bf r} = 3r/8$ and ${\bf t} = s$.
This application may be made provided that
$$ 
5r/4  \leq c/2 \cdot \hat{n}^{1/9} \, , \, 3r/8 \leq \rsc/4 \cdot \hat{n}^{1/9}  \, , \, s \in \big[ 2^{7/2} , 2 \hat{n}^{1/3} \big] \, \, \textrm{ and } \, \, \hat{n} \geq c^{-18} \, .
$$ 
 This completes the proof of Lemma~\ref{l.combine}. \qed
 
\medskip

\noindent{\bf Proof of Proposition~\ref{p.feventfirstbound}(1).} Here we implement the {\em near} case argument outlined before the proof of Proposition~\ref{p.polyfluc}. 
Consider the event 
\begin{eqnarray}
  & &  \fluc \, [-(r+1),-r]  \, \cap \, \Big\{ t_{1,2}^{-1/3} (1-a)^{-1/3} \weight\big(\rho[3]\big) > - 2^{-1/2}  - s \Big\} \label{e.ftriple} \\
  & & \qquad  \cap \, \, \Big\{  t_{1,2}^{-1/3} (1-a)^{-1/3} \weight\big(\rho[2]\big) < - 2^{-1/2} (7r/8)^2 + s \Big\} \, . \nonumber
\end{eqnarray}
Note that $ \weight\big(\rho[3]\big)  -  \weight\big(\rho[2]\big) =  \weight\big(\hat\rho\big) -  \weight\big(\rho \big)$.

Note that $- 2^{-1/2}  - s = \big( - 2^{-1/2} (7r/8)^2 + s \big) + s$ is solved by $s =  3^{-1} \cdot 2^{-1/2}  \big( (7r/8)^2 - 1 \big)$. Setting the value of $s$ in this way,
note that $s \geq 5/3 \cdot 2^{-7/2} r^2$ since $r \geq 4$.
We find then that the $\PP$-probability of the above event is bounded above by
\begin{equation}\label{e.rhohatrho}
 \PP \bigg(  t_{1,2}^{-1/3} (1-a)^{-1/3} \Big( \weight\big(\hat\rho\big) -  \weight\big(\rho \big)  \Big) \geq  d r^2  \bigg)  \, ,
\end{equation}
where here we write $d = 5/3 \cdot 2^{-7/2}$.

As a temporary notation, we denote the {\em maximum weight difference} $\mdeltaweight_{n;(I,t_1)}^{(J,t_2)}$ to be
 the supremum over $x_1,x_2 \in I$
and $y_1,y_2 \in J$ of
$\big\vert \weight_{n;(x_1,t_1)}^{(y_1,t_2)} -  \weight_{n;(x_2,t_1)}^{(y_2,t_2)} \big\vert$.

Set $I =  \{ x \}$ and $J = [y-9\Rmac/8,y]$.  We claim that, when  $r \geq 8(1-a)^{1/3} \tot^{-2/3} \vert x - y \vert$,
\begin{equation}\label{e.incverify}
 \Big\{ t_{1,2}^{-1/3} (1-a)^{-1/3} \big(  \weight(\hat\rho)  -  \weight(\rho) \big)  \geq d r^2 \Big\} \, \subseteq \, \Big\{ \mdeltaweight_{n;(I,t_1)}^{(J,t_2)} \geq d t_{1,2}^{1/3} (1-a)^{1/3}  r^2 \Big\} \, .
\end{equation}
Indeed, note that
$\weight_{n;(x,t_1)}^{(z - \Rmac,t_2)} \geq \weight(\hat{\rho})$
while 
$\weight_{n;(x,t_1)}^{(y,t_2)} = \weight(\rho)$. Moreover, 
 the lower bound on~$r$ ensures that $z - \Rmac \geq y - 9\Rmac/8$, 
so that  $\mdeltaweight_{n;(I,t_1)}^{(J,t_2)} \geq  \weight(\hat\rho)  -  \weight(\rho)$. This implies~(\ref{e.incverify}).

Since $\Rmac = t_{1,2}^{2/3} (1-a)^{2/3} r$, the right-hand event in~(\ref{e.incverify}) equals 
\begin{equation}\label{e.mdw.bound}
t_{1,2}^{-1/3} \mdeltaweight_{n;(I,t_1)}^{(J,t_2)} \geq  (8/9)^{1/2} d  r^{3/2}  \cdot  \big(9 
t_{1,2}^{-2/3}\Rmac/8 \big)^{1/2}   \, .
\end{equation}
We seek to bound above the probability of this event 
by using 
%Corollary~\ref{c.ordweight}.
Theorem~\ref{t.differenceweight}. 
Two obstacles are that 
 this result is stated for the special case where the start time is zero and the end time is one, and the presence of parabolic terms $Q$ in its statement.
 Regarding the first difficulty, we use the scaling principle from Section~\ref{s.scalingprinciple} to obtain a general form for  
%Corollary~\ref{c.ordweight}.
Theorem~\ref{t.differenceweight} which will be applicable to the time pair $(t_1,t_2)$ with which we work. 
 The conclusion of this general form of the result is that
$$
\PP \left( \sup_{\begin{subarray}{c} u_1,u_2 \in [x,x+\e \tot^{2/3}] \\
    v_1,v_2 \in [y,y+\e \tot^{2/3}]  \end{subarray}} \left\vert \, \tot^{-1/3} \weight_{n;(u_2,t_1)}^{(v_2,t_2)} +  Q\left(\frac{v_2 - u_2}{\tot^{2/3}}\right)  - \tot^{-1/3} \weight_{n;(u_1,t_1)}^{(v_1,t_2)} - Q\left(\frac{v_1 - u_1}{\tot^{2/3}}\right) \, \right\vert  \, \geq \, \e^{1/2}
  R  \right)
$$
  is at most  $10032 \, C  \exp \big\{ - c_1 2^{-21}   R^{3/2}   \big\}$; 
% 
%$$
% \PP \Big( \neg \, \lwr_{n; (  [x,x+\e\tot^{2/3}],t_1 )}^{( [y,y+\e\tot^{2/3}],t_2 )}\big( \e , R  \big)    \Big) \leq 
% 10012 \, C  \exp \big\{ - c_1 2^{-21} R^{3/2}     \big\}
% \, ,
%$$
while the new result's  hypotheses are those of  
%Corollary~\ref{c.ordweight}
Theorem~\ref{t.differenceweight} 
after making the replacements $n \to n\tot$ and $(x,y) \to \tot^{-2/3} \cdot (x,y)$.

We will find an upper bound on the event~(\ref{e.mdw.bound}) by applying this  
general version of 
%Corollary~\ref{c.ordweight}
Theorem~\ref{t.differenceweight}, with parameter choices
with ${\bm \e} =  9
t_{1,2}^{-2/3} \Rmac/8$, which also equals $9(1-a)^{2/3}r/8$; as well as ${\bf x} = x$, ${\bf y} = y    -  9\Rmac/8$, and ${\bf R} = 
  2^{-1}(8/9)^{1/2} d r^{3/2}$.
In order that this indeed provide such an upper bound, we must first overcome the second obstacle just mentioned, confirming that  parabolic curvature is suitably controlled. The relevant bound is  
$$
  \sup_{\begin{subarray}{c} u_1,u_2 \in  [{\bf x},{\bf x}+ {\bm \e} \tot^{2/3}] \\
    v_1,v_2 \in [{\bf y}, {\bf y}+ {\bm \e}\tot^{2/3}]  \end{subarray}} \big\vert Q \big( \tot^{-2/3} (v_2 - u_2) \big) -  Q \big( \tot^{-2/3} (v_1 - u_1) \big) \big\vert \leq 2^{-1} (8/9)^{1/2} d  r^{3/2}  \cdot  \big(9 
t_{1,2}^{-2/3}\Rmac/8 \big)^{1/2} \, .
$$
This bound holds because this supremum is at most 
    $$
    2{\bm \e}\tot^{2/3} \cdot 2^{1/2} \big( \vert {\bf y} - {\bf x} \vert  + 2{\bm \e} \tot^{2/3} \big) \tot^{-4/3} = 2^{-3/2} \cdot 9 (1-a)^{2/3} r \Big( \vert x - y \vert \tot^{-2/3} + 9 (1 - a)^{2/3} r/4 \Big) \, .
    $$ 
Expanding the right-hand bracket, the resulting two terms are both at most $4^{-1} (8/9)^{1/2} d  r^{3/2}  \cdot  \big(9 
t_{1,2}^{-2/3}\Rmac/8 \big)^{1/2}$ (which suffices for our purpose) respectively provided that
    $$
    r \geq 3^3 2^4 5^{-1} (1-a)^{1/3}  \vert x - y \vert \tot^{-2/3}  \, \, \, \textrm{and} \, \, \, 1 - a \leq 5 \cdot 3^{-5} 2^{-2} 
    $$
    since $\tot^{-2/3} \Rmac = (1-a)^{2/3}r$ and $d = 5/3 \cdot 2^{-7/2}$.
    % so that 9 \cdot 2^{1/2} d^{-1} = 3^3 2^4 5^{-1} 

 We further mention that, for  the application in question to be made, it suffices that 
\begin{eqnarray*}
 & & 9(1-a)^{2/3}r/8 \leq 2^{-4} \, , \, 
 n\tot \geq 10^{32} c^{-18} \, , \, \\
 & &  \big\vert x - y  \big\vert  \tot^{-2/3} +  9(1-a)^{2/3}r/8 \leq 2^{-2} 3^{-1} \rsc  (n\tot)^{1/18} \, \, \textrm{ and } \, \,   2^{1/2}/3 \cdot d r^{3/2} \in \big[ 10^4 \, , \,   10^3 (n \tot)^{1/18} \big] \, ,
\end{eqnarray*}
as well as $n \tot \in \N$.

With these pieces in place, we apply the general version of Theorem~\ref{t.differenceweight}, learning that the event in~(\ref{e.mdw.bound}) has probability at most
   $$
 10032 \, C  \exp \big\{ - c_1 2^{-26}   3^{-3} 5^{3/2} r^{9/4}        \big\}   \, .
   $$
 Recall that this quantity is thus an upper bound on the probability of the event~(\ref{e.ftriple}). Combining this information with the two parts of Lemma~\ref{l.combine}, we find that 
\begin{eqnarray*}
 \PP \Big( \, \fluc \, [-(r+1),-r] \, \Big) & \leq  &  
 140 C \, \exp \big\{ - c_1 s^{3/2} \big\}  \,   + \, 
   \big(3r/8 + 1 \big) \cdot 6 C  \exp \big\{ - 2^{-11/2} c s^{3/2} \big\} \\
   & & \qquad \qquad \qquad \qquad
   +  \, 
 10032 \, C  \exp \big\{ - c_1 2^{-26}   3^{-3} 5^{3/2} r^{9/4}        \big\}   \, .
\end{eqnarray*}

Recalling that  $s \geq 5/3 \cdot 2^{-7/2} r^2$, we obtain
\begin{eqnarray*}
 \PP \Big( \, \fluc \, [-(r+1),-r] \, \Big) & \leq   & 
 \big( 146  +  9r/4 \big) C \exp \big\{ - c_1  2^{-43/4} 3^{-3/2} 5^{3/2}  r^3 \big\} \\
  & & \qquad  + \,   
 10032 \, C  \exp \big\{ - c_1 2^{-26} 3^{-3} 5^{3/2}  r^{9/4}     \big\}  \\
  & \leq &  
 (10178 + 3r) \, C   \exp \big\{ - c_1 2^{-26} 3^{-3} 5^{3/2}  r^{9/4} \big\} \, ,
\end{eqnarray*}
where we used $c_1 \leq c$ in the first inequality  and  $r \geq 1$ in the second.
This completes the proof of Proposition~\ref{p.feventfirstbound}(1). \qed

%\begin{theorem}\label{t.differenceweight.revise}
%Let $n \in 2\N$ and   $x,y \in \R$  satisfy 
% $n \geq 10^{29} c^{-18}$ and   $\big\vert x - y  \big\vert \leq 2^{-5/3} 3^{-1} \rsc  n^{1/18}$.
%Let $\e \in (0,2^{-4}]$ and 
% $R \in \big[10^4 \, , \,   10^3 n^{1/18} \big]$.
%Then
%$$
%\PP \left( \sup_{\begin{subarray}{c} u_1,u_2 \in [x,x+\e] \\
%    v_1,v_2 \in [y,y+\e]  \end{subarray}} \Big\vert \weight_{n;(u_2,0)}^{(v_2,1)} + Q(v_2 - u_2) - \weight_{n;(u_1,0)}^{(v_1,1)} - Q(v_1 - u_1) \Big\vert  \, \geq \, \e^{1/2}
%  R  \right)
%$$
%  is at most  $10012 \, C  \exp \big\{ - c_1 2^{-19}   R^{3/2}   \big\}$.
%\end{theorem}

%\begin{corollary}\label{c.ordweight}
%Let 
%$n \in 2\N$ and  $x,y \in \R$ satisfy 
%$n \geq 10^{29} c^{-18}$ and   $\big\vert x - y  \big\vert \leq \e^{-1/2} \wedge 2^{-5/3} 3^{-1} \rsc  n^{1/18}$.
%Let  $\e \in (0,2^{-4}]$ and
% $R \in \big[2 \cdot 10^4 \, , \,   10^3 n^{1/18} \big]$.
%Then
%\begin{equation}\label{e.ordweight}
% \PP \Big( \neg \, \lwr_{n;([x,x+\e],0)}^{([y,y+\e],1)}\big( \e , R  \big)    \Big) \leq 
% 10012 \, C  \exp \big\{ - c_1 2^{-20 - 1/2} R^{3/2}     \big\}
% \, .
%\end{equation}
%\end{corollary}
%*** 

\noindent{\bf Proof of Proposition~\ref{p.feventfirstbound}(2).} We now provide an alternative upper bound on the probability of the event in the display beginning at~(\ref{e.ftriple}), one that treats the middle-distant case. Recalling that this probability is at most~(\ref{e.rhohatrho}), we will find a new upper bound on the probability~(\ref{e.rhohatrho}). Recalling that we set $d = 5/3 \cdot 2^{-7/2}$, the occurrence of the event whose probability is in question entails that 
\begin{eqnarray*}
\textrm{either} & & t_{1,2}^{-1/3} \weight\big(\hat\rho\big)  \geq  2^{-1} d  (1-a)^{1/3}  r^2 \, - \, 2^{-1/2}(x-y)^2 \tot^{-4/3}  \\
\textrm{or} & &    
 t_{1,2}^{-1/3}  \weight(\rho)  \leq  -  2^{-1} d  (1-a)^{1/3}  r^2 \, - \, 2^{-1/2}(x-y)^2 \tot^{-4/3} \, .
\end{eqnarray*}

We now state and prove a claim regarding the probabilities of these two outcomes.

\noindent{\bf Claim: (1).} The former outcome has $\PP$-probability at most 
$\rsC \exp \big\{ - \rsc  (1-a)^{1/2}  5^{3/2} 3^{-3/2}  2^{-33/4} r^3 \big\}$.
\noindent{\bf (2).} The latter outcome has  $\PP$-probability at most
$\rsC \exp \big\{ - \rsc  (1-a)^{1/2}   5^{3/2} 3^{-3/2}  2^{-27/4} r^3 \big\}$.

First we prove Claim~(2). Since 
$\weight_{n;(x,t_1)}^{(y,t_2)} = \weight(\rho)$, the latter outcome entails 
$$ \scaledle_{n;t_1}^{\downarrow;(y,t_2)}\big( 1, \omega' \big) +2^{-1/2} \omega'^2  \leq -  2^{-1} d (1 - a)^{1/3}  r^2
$$ 
where we set $\omega' =  \tot^{-2/3} ( x - y )$. 
The probability of the outcome is thus bounded above by
\begin{eqnarray*}
&  &  \PP \Big( \scaledle_{n;t_1}^{\downarrow;(y,t_2)}\big( 1, \omega' \big)  + 2^{-1/2} \omega'^2 \leq - 5/3 \cdot 2^{-9/2} (1 - a)^{1/3}   r^2 \Big)  \\
 & \leq  & \rsC \exp \big\{ - \rsc  (1-a)^{1/2}  (5/3)^{3/2} \cdot 2^{-27/4} r^3 \big\} \, ,
\end{eqnarray*}
where the displayed inequality is a consequence of an application of one-point lower tail $\rmreg(2)$  with the ensemble equalling the  $(n\tot + 1)$-curve  $\scaledle_{n;t_1}^{\downarrow;(y,t_2)}$, and where the parameters are chosen to be ${\bf z} = \omega'$ and  ${\bf s} =   5/3 \cdot 2^{-9/2} (1-a)^{1/3}  r^2$. This application may be made when 
$\vert \omega' \vert \leq c (n \tot)^{1/9}$
and $5/3 \cdot 2^{-9/2} (1-a)^{1/3}  r^2 \in \big[ 1,(n \tot)^{1/3} \big]$.

We now derive Claim~$(1)$. This is an  assertion about $\hat{\rho}$, an $n$-zigzag with starting and ending points $(x,t_1)$ and $(z - \Rmac,t_2)$. Note that the quadratic correction term in the assertion is 
 $2^{-1/2}(y - x)^2 \tot^{-4/3}$, whereas we would like to instead work with  $2^{-1/2}(z-\Rmac - x)^2 \tot^{-4/3}$, this being the natural quadratic correction associated to the starting and ending points of $\hat\rho$.
 We begin then by noting that
 %, under the hypotheses of Lemma~\ref{l.feventfirstbound}(2), 
 the difference is suitably small:
\begin{equation}\label{e.diffsuit}
 \Big\vert \, 2^{-1/2}(y-x)^2 \tot^{-4/3} \, - \, 2^{-1/2}(z-\Rmac - x)^2 \tot^{-4/3} \, \Big\vert \, \leq \, 4^{-1} d (1-a)^{1/3} r^2 \, . 
\end{equation}
An explicit condition that suffices to ensure this bound is that $r \geq 11  (1-a)^{1/3} \vert y -x \vert   \tot^{-2/3}$ and  $a \geq 1 - 10^{-2}$.
See Lemma~$D.4$ in the online Appendix~$D$.
%\ref{l.fcalc}.
%arxiv/submit

We see from~(\ref{e.diffsuit}) and 
$\weight_{n;(x,t_1)}^{(z - \Rmac,t_2)} \geq \weight(\hat{\rho})$ that the outcome that Claim~$(1)$ concerns entails that
$$
 \scaledle_{n;t_1}^{\downarrow;(z-\Rmac,t_2)} \big( 1, \omega \big)   + 2^{-1/2} \omega^2  \geq  4^{-1} d (1-a)^{1/3}    r^2 \, \, 
$$
where now we set $\omega =  \tot^{-2/3} ( x - z + \Rmac)$.
Next we apply one-point upper tail $\rmreg(3)$ to the $(n\tot + 1)$-curve  ensemble $\scaledle_{n;t_1}^{\downarrow;(z-\Rmac,t_2)}$ with parameter settings ${\bf z} = \omega$ and  ${\bf s} =  4^{-1} d (1-a)^{1/3}    r^2$. We learn that 
$$
\PP \Big( \scaledle_{n;t_1}^{\downarrow;(z-\Rmac,t_2)} \big( 1, \omega \big)  + 2^{-1/2} \omega^2  \geq  4^{-1} d (1-a)^{1/3}    r^2 \Big)  \, \leq  \, \rsC \exp \Big\{ - \rsc  (1-a)^{1/2}  (5/3)^{3/2} \cdot 2^{-33/4} r^3 \Big\} \, .
$$
 This application of $\rmreg(3)$ 
may be made when  $\vert \omega \vert \leq c (n \tot)^{1/9}$
and $(1-a)^{1/3} \cdot  5/3 \cdot 2^{-11/2} r^2 \geq 1$.
We have proved the Claim.

The sum of the two expressions in Claim~$(1)$ and~$(2)$
the sought alternative upper bound on the probability of the event in the display beginning at~(\ref{e.ftriple}). Combining again  with the two parts of Lemma~\ref{l.combine}, 
 we find that
\begin{eqnarray*}
 \PP \Big( \fluc \,  [-(r+1),-r] \Big) & \leq  &  
 140 C \, \exp \big\{ - c_1 s^{3/2} \big\}  \,  + \, 
   \big(3r/8 + 1 \big) \cdot 6 C  \exp \big\{ - 2^{-11/2} c s^{3/2} \big\} \\
   & & \qquad \qquad \qquad \qquad \,  + \, 
2\rsC \exp \big\{ - \rsc  (1-a)^{1/2}  5^{3/2} 3^{-3/2}  2^{-33/4} r^3 \big\}   
      \, . 
\end{eqnarray*}

Recalling that  $s \geq 5/3 \cdot 2^{-7/2} r^2$, we obtain
\begin{eqnarray*}
  \PP \Big( \fluc \,  [-(r+1),-r] \Big) & \leq  &  
 \big( 146  +  9r/4 \big) C \exp \big\{ - c_1  2^{-43/4} 3^{-3/2} 5^{3/2}  r^3 \big\}  \\
  & & \qquad   \,  + \,  2 \rsC \exp \big\{ - \rsc  (1-a)^{1/2}  3^{-3/2} 5^{3/2}  2^{-33/4} r^3 \big\} \\
 & \leq &  
 \big( 148 C +  9r/4 \big)  \exp \big\{ - c_1  (1-a)^{1/2}   3^{-3/2} 5^{3/2} 2^{-43/4} r^3 \big\} \, ,
\end{eqnarray*}
where we used $c_1 \leq \rsc$ in both inequalities. 
This completes the proof of Proposition~\ref{p.feventfirstbound}(2). \qed

\medskip

\noindent{\bf Proof of Lemma~\ref{l.feventhigh}.} Now we treat the far case.
In outline, the one-point lower tail $\rmreg(2)$ shows that neither the weight of $\rho[1]$ nor of $\rho[2]$ can be extremely low, except with a tiny probability. In the opposing case, the polymer weight $\weight(\rho) = \weight\big( \rho[1] \big) + \weight \big( \rho[2] \big)$ is not highly negative. However, this eventuality is unlikely in the far case, where $Z$ is very much less than $z$, due to collapse near infinity  
Proposition~\ref{p.mega}(4).

Naturally, here we are continuing to denote  $Z = \rho\big( (1-a) t_1 + a t_2 \big)$ and $z = (1-a)x + ay$.
To begin the formal argument, note that  the far case event $\fluc \big(-\infty,-2\hat{n}^{1/36} \big]$ occurs 
 if and only if $Z - z \leq  - \tot^{2/3} (1 -a)^{2/3} \cdot 2 \hat{n}^{1/36}$. 
Note that  $\hat{n}^{1/36} \geq (1 -a)^{1/3}  \tot^{-2/3} \vert x - y \vert$ implies that $\vert z - y \vert \leq \hat{n}^{1/36} \cdot  \tot^{2/3} (1 -a)^{2/3}$, so that the occurrence of $\fluc \big(-\infty,-2\hat{n}^{1/36} \big]$ entails that $Z - y \leq  -   \tot^{2/3} (1 -a)^{2/3}  \hat{n}^{1/36}$. 
Using again the notation $V'  = (1-a)^{-2/3} \tot^{-2/3} (Z-y)$, the latter condition is given by $V' \leq - \hat{n}^{1/36}$. 

Recall that in~(\ref{e.vprime.ensemble}), the polymer weight $\weight \big( \rho[2] \big)$ is expressed in terms of $\scaledle_{n;(1-a) t_1 + a t_2}^{\downarrow;(y,t_2)} \big( 1,V' \big)$. We thus find that
\begin{eqnarray}
  & &  \PP \Big(  \fluc \big(-\infty,-2\hat{n}^{1/36} \big] \, , \,  \tot^{-1/3} (1-a)^{-1/3}  \weight\big( \rho[2] \big) \geq  \big( - 2^{-1/2} + 2^{-5/2} \big)  (1-a)^{1/18} \tot^{1/18} n^{1/18}
   \Big) \nonumber \\
  & \leq &  \PP \Big( \sup_{x \leq - \hat{n}^{1/36}} \scaledle_{n;(1-a) t_1 + a t_2}^{\downarrow;(y,t_2)} \big( 1,x \big) \geq  - \big(  2^{-1/2} - 2^{-5/2} \big)  (1-a)^{1/18} \tot^{1/18} n^{1/18}
  \Big) \, . \label{e.ensemblecollapse}
\end{eqnarray}

 Recall that $\hat{n} = n (1-a) \tot$. To find an upper bound on the probability~(\ref{e.ensemblecollapse}), we now apply collapse near infinity 
 Proposition~\ref{p.mega}(4)
% $\rmreg(4)$ 
 to the $(c,C)$-regular $(\hat{n} + 1)$-curve ensemble  $\scaledle_{n;(1-a) t_1 + a t_2}^{\downarrow;(y,t_2)}$. We set the parameter ${\bm \eta}$ so that  ${\bm \eta}( \hat{n} + 1) ^{1/9} = \hat{n}^{1/36}$. For the application to be made, it is enough that
 $$
 \big(n (1-a) \tot \big)^{-1/12}  \leq \rsc \, \, \textrm{ and } \, \, n (1-a) \tot \geq \big( 2^{5/4} \rsc^{-1} \big)^{9} \, .
 $$
Note then that 
$$ 
\ell\big(-\hat{n}^{-1/12} (\hat{n} + 1)^{1/9}   \big) =  \big( - 2^{-1/2} + 2^{-5/2} \big) \big(   \hat{n}^{-1/12} \big)^2 (\hat{n} + 1)^{2/9} = \big( - 2^{-1/2} + 2^{-5/2} \big)     ( \hat{n}^{-1/12})^2 (\hat{n} + 1)^{2/9} \, .
$$
Since $\ell(x)$ is increasing for $x \leq 0$, we find that 
$\ell\big(-\hat{n}^{-1/36}   \big) \leq  - \big(  2^{-1/2} - 2^{-5/2} \big)  \hat{n}^{1/18}$. Thus, the expression~(\ref{e.ensemblecollapse}) is found to be at most
\begin{equation}\label{e.firstub}
6C \exp \Big\{ - c \eta^3  2^{-15/4}   \big( \hat{n} + 1 \big)^{1/3} \Big\}  = 
6C \exp \Big\{ - 2^{-15/4} c  \hat{n}^{1/12} \Big\} \, .
\end{equation}
  
Set   $\tilde{V} =  a^{-2/3} \tot^{-2/3} ( Z - x)$.
  Note that 
 \begin{equation}\label{e.weightrhoone}
   \weight \big( \rho[1] \big)  =  a^{1/3} \tot^{1/3} \scaledle_{n;(x,t_1)}^{\uparrow;  (1-a)t_1 + at_2} \big( 1, \tilde{V} \big) \, ,
 \end{equation}
 since both of these quantities equal $\weight_{n;(x,t_1)}^{\big(Z, (1-a)t_1 + at_2 \big)}$.
 
  We claim that 
  \begin{equation}\label{e.vtilde.claim}
  \fluc \big(-\infty,-2\hat{n}^{1/36} \big]  \,\subseteq  \, \Big\{ \tilde{V} \leq -  (1 -a)^{2/3}  \hat{n}^{1/36} \Big\} \, .
  \end{equation}
  Indeed, note that the left-hand event entails that
  $$
  Z - x  \leq  a \vert x -y \vert  - \tot^{2/3} (1 -a)^{2/3} \cdot 2 \hat{n}^{1/36} \leq  -  \tot^{2/3} (1 -a)^{2/3}  \hat{n}^{1/36} 
  $$
  since $\vert x - y \vert \leq \tot^{2/3} (1 -a)^{2/3}  \hat{n}^{1/36}$ and $a \leq 1$. From another use of $a \leq 1$, 
  we obtain~(\ref{e.vtilde.claim}).

From~(\ref{e.weightrhoone}) and~(\ref{e.vtilde.claim}), we find that
\begin{eqnarray}
  & &  \PP \Big(  
  \fluc \big(-\infty,-2\hat{n}^{1/36} \big]   \, , \,  \tot^{-1/3} a^{-1/3}  \weight\big( \rho[1] \big) \geq  - \big(  2^{-1/2} - 2^{-5/2} \big)   (1-a)^{25/18} \tot^{1/18} n^{1/18}
   \Big) \nonumber \\
  & \leq &  \PP \Big( \sup_{x \leq - (1 -a)^{2/3} \hat{n}^{1/36}} \scaledle_{n;(x,t_1)}^{\uparrow;  (1-a)t_1 + at_2} \big( 1,x \big) \geq - \big(  2^{-1/2} - 2^{-5/2} \big) 
   (1-a)^{25/18} \tot^{1/18} n^{1/18}   \Big) \, . \label{e.hatprob}
\end{eqnarray}

To find an upper bound on the latter probability,
we  apply Proposition~\ref{p.mega}(4)
 % $\rmreg(4)$ 
  to the $(c,C)$-regular ensemble $\scaledle_{n;(x,t_1)}^{\uparrow;  (1-a)t_1 + at_2}$.
This ensemble has $\tilde{n} + 1$ curves, where here we introduce the shorthand $\tilde{n} = na\tot$.
 We set the Proposition~\ref{p.mega}(4) parameter ${\bm \eta}$
 so that  ${\bm \eta} ( \tilde{n} + 1 )^{1/9} =  (1 -a)^{2/3} \hat{n}^{1/36}$. Recalling that $\hat{n} = n (1-a) \tot$, we see that $\eta \leq n^{-1/12} a^{-1/9} (1-a)^{25/36} \tot^{-1/12}$.
 For  Proposition~\ref{p.mega}(4)
  to be applied in this way, it suffices then that
 $$
  n^{-1/12} a^{-1/9} (1-a)^{25/36} \tot^{-1/12} \leq c \, \, \textrm{ as well as } \, \,  n a \tot \geq \big( 2^{5/4} \rsc^{-1} \big)^{9} \, .
 $$
Note then that 
\begin{eqnarray*}
\ell\big(- (1 -a)^{2/3}  \hat{n}^{1/36} \big) & = &  - \big(  2^{-1/2} - 2^{-5/2} \big) \big(   (1 -a)^{2/3} \hat{n}^{1/36} (  \tilde{n} + 1 )^{-1/9} \big)^2 (\tilde{n} + 1)^{2/9} \\
 & = & - \big(  2^{-1/2} - 2^{-5/2} \big)  (1-a)^{25/18} \tot^{1/18} n^{1/18} \, .
\end{eqnarray*}
Also using that  $\ell(x)$ is increasing for $x \leq -  (1 -a)^{2/3} \hat{n}^{1/36}$, we see that our application of Proposition~\ref{p.mega}(4)
 implies that the probability~(\ref{e.hatprob}) is at most
\begin{equation}\label{e.secondub}
6C \exp \Big\{ - c \eta^3  2^{-15/4}   (\tilde{n} + 1)^{1/3} \Big\}  = 
6C \exp \Big\{ -   2^{-15/4}  c   (1-a)^{25/12} \tot^{1/12} n^{1/12}   \Big\}   \, .
\end{equation}
  Note that $\tot^{1/3} \scaledle_{n;(x,t_1)}^{\uparrow;t_2}\big(1, \tot^{-2/3}(y-x)\big) = \mc{L}_{n;(x,t_1)}^{\uparrow;t_2}(1,y) = \weight (\rho) =  \weight \big(\rho[1] \big) +  \weight \big(\rho[2] \big)$.
 Thus, 
\begin{eqnarray}
  & &  \PP \Big(  \tot^{-1/3} a^{-1/3}  \weight\big( \rho[1] \big) \leq  - \big(  2^{-1/2} - 2^{-5/2} \big)   (1-a)^{25/18} \tot^{1/18} n^{1/18} \, , \,  \nonumber \\
   & & \qquad \qquad  \tot^{-1/3} (1-a)^{-1/3}  \weight\big( \rho[2] \big) \leq  - \big(  2^{-1/2} - 2^{-5/2} \big)  (1-a)^{1/18} \tot^{1/18} n^{1/18}
   \Big) \nonumber \\
   & \leq & \PP \bigg( \scaledle_{n;(x,t_1)}^{\uparrow}\big(1, \tot^{-2/3}(y-x)\big) \leq  - \big( 2^{-1/2} - 2^{-5/2} \big)  \Big( (1-a)^{25/18} a^{1/3}  +  (1-a)^{1/18 + 1/3} \Big) \tot^{1/18} n^{1/18} \bigg) \nonumber \\
   & \leq & C \exp \bigg\{  - c  \big(  2^{-1/2} - 2^{-5/2} \big)^{3/2}  \Big( (1-a)^{25/18} a^{1/3} +  (1-a)^{17/18} \Big)^{3/2} \tot^{1/12} n^{1/12}   \bigg\} \, , \label{e.thirdub}
\end{eqnarray}
where in the latter inequality, we applied one-point lower tail $\rmreg(2)$ to the $(n\tot + 1)$-curve ensemble $\scaledle_{n;(x,t_1)}^{\uparrow;t_2}$.
In this application, we take 
$$
{\bf z} = \tot^{-2/3}(y-x) \, \, \textrm{ and }  \, \, {\bf s} =  \big(  2^{-1/2} - 2^{-5/2} \big)  \big( (1-a)^{25/18} a^{1/3}  +  (1-a)^{17/18} \big) \tot^{1/18} n^{1/18} \, .
$$
The application may be made provided that
$$
 \vert y - x \vert \tot^{-2/3} \leq \rsc (n \tot )^{1/9} \, \,
 \textrm{ and } \, \,
  \big(  2^{-1/2} - 2^{-5/2} \big)  \big( (1-a)^{25/18} a^{1/3}  +  (1-a)^{17/18} \big) \tot^{1/18} n^{1/18}  \in \big[ 1, (n\tot)^{1/3} \big] \, .
$$

Recalling the upper bounds on the probabilities in the first lines of the displays (\ref{e.ensemblecollapse}) and (\ref{e.hatprob}) offered by (\ref{e.firstub}) and (\ref{e.secondub}), 
and combining these with the bound (\ref{e.thirdub}), 
we find that  $\PP \big( \fluc (-\infty,-2\hat{n}^{1/36} ]  \big)$ is at most
\begin{eqnarray*}
 & & 
6C \exp \Big\{ - 2^{-15/4} c (1-a)^{1/12} \tot^{1/12}   n^{1/12} \Big\}
    \,  + \, 
6C \exp \Big\{ -   2^{-15/4}  c   (1-a)^{25/12} \tot^{1/12} n^{1/12}   \Big\}   \\
 & & + \,   C \exp \bigg\{  - c  \big(  2^{-1/2} - 2^{-5/2} \big)^{3/2}  \Big( (1-a)^{25/18} a^{1/3} +  (1-a)^{17/18} \Big)^{3/2} \tot^{1/12} n^{1/12}   \bigg\}  \\
 & \leq & 13C \exp \Big\{ - 2^{-15/4} c (1-a)^{25/12} \tot^{1/12}       n^{1/12} \Big\} \, .
\end{eqnarray*}
% $\big(  2^{-1/2} - 2^{-5/2} \big)^{3/2} = 0.386$ and $2^{-15/4} = 0.074$.
This completes the proof of Lemma~\ref{l.feventhigh}. \qed

\newpage

\appendix

\section{Glossary of notation}\label{s.glossary}

This article uses quite a lot of notation. Each line of the list that we now present recalls one of the principal pieces of notation; provides a short summary of its meaning; and gives the number of the page at which the notation is either introduced or formally defined. The summaries are, of course, imprecise: phrases in quotation marks, such as `disjointly travel', are merely verbal approximations of a precise meaning. When a quantity is said to be `$\leq r$', it is in fact the absolute value of the roughly recalled quantity which must be at most $r$.

\bigskip
\def\qq{&}

\begin{center}
\halign{
\!\!\!\!\!\!\!\!\!\!\!\!#\quad \!\!\! \! \hfill&#\quad\hfill&\!\!\!\!\!\quad\hfill#\cr

staircase \qq (a geometric depiction of) an unscaled Brownian LPP path  \fff{staircase}
energy \qq the value assigned to a staircase by Brownian LPP  \fff{energy} 
geodesic \qq a staircase of maximum energy given its endpoints  \fff{geodesic}     
$R_n$ \qq the linear {\em scaling map} \fff{scalingmap}
zigzag \qq the image of a staircase under the scaling map  \fff{zigzag}  
polymer \qq a zigzag that is the image of a geodesic -- and thus of maximum {\em weight}  \fff{polymer}  
compatible triple \qq a triple $(n,t_1,t_2)$ so that $[t_1,t_2]$
is the lifetime of a zigzag  \fff{comptriple}  
$\tot$ \qq the difference $t_2 - t_1$; the lifetime of a given polymer in most applications \fff{tot}
$\rho_{n;(x,t_1)}^{(y,t_2)}$ \qq the polymer of journey $(x,t_1) \to (y,t_2)$  \fff{polynot}
$\maxpoly_{n;(I,t_1)}^{(J,t_2)}$ \qq the maximum cardinality of a set of polymers travelling $(I,t_1) \to (J,t_2)$ \fff{maxcard}
$\pdr_{n;(x,t_1)}^{(y,t_2)}\big(a,r\big)$ \qq $\big\{ \!\!$ $\rho_{n;(x,t_1)}^{(y,t_2)}$ has `normalized' fluctuation $\leq r$ at lifetime fractions $a$ and $1-a$ $\!\!\big\}$  \fff{polydevreg}
weight  \qq the scaled energy, assigned to any zigzag  \fff{weight}  
$\weight_{n;(x,t_1)}^{(y,t_2)}$ \qq the weight of the polymer $\rho_{n;(x,t_1)}^{(y,t_2)}$  \fff{maxweight}  
 $\nearpoly_{n,k;(x,0)}^{(y,1)}(\eta)$  \qq $\big\{ \! \!$ $k$ zigzags `disjointly travel' $(x,0) \to (y,1)$
with weight shortfall at most $\eta$ $\!\!\big\}$ \fff{nearpoly}
$\weight_{n,k;(\bar{x},t_1)}^{(\bar{y},t_2)}$  \qq the maximum weight of $k$ `disjointly travelling' zigzags $(\bar{x},0) \to (\bar{y},1)$ \fff{collectiveweight}
 $\rho_{n,k;(\bar{x},t_1)}^{(\bar{y},t_2)}$  \qq the multi-polymer, the maximizer that achieves 
 $\weight_{n,k;(\bar{x},t_1)}^{(\bar{y},t_2)}$  \fff{multipolymer}
 $\rho_{n,k,i;(\bar{x},t_1)}^{(\bar{y},t_2)}$ \qq the zigzag that is 
 the $i\textsuperscript{th}$ component of  $\rho_{n,k,i;(\bar{x},t_1)}^{(\bar{y},t_2)}$ \fff{multipolymercomp}
$\mc{L}_{n;(x,t_1)}^{\uparrow;t_2}$ \qq the scaled forward line ensemble rooted at $(x,t_1)$
of duration $\tot$ \fff{scaledforwardle}
$\mc{L}_{n;t_1}^{\downarrow;(y,t_2)}$ \qq  the scaled backward line ensemble rooted at $(y,t_2)$
of duration $\tot$  \fff{scaledbackwardle}
 $\scaledle_{n;(x,t_1)}^{\uparrow;t_2}$,
 $\scaledle_{n;t_1}^{\downarrow;(y,t_2)}$ \qq the counterpart {\em normalized} forward and backward ensembles \fff{normalizedle}$Q$ \qq $Q(x) = 2^{-1/2} x^2$, so that $-Q$
is any normalized ensemble's global shape \fff{parabola}   
regular ensemble \qq Brownian Gibbs ensemble with one-point control, hewing to $Q$ globally   \fff{regular}
the scaling principle \qq a basic result that reduces polymer arguments to case of lifetime $[0,1]$ \fff{scaling}
$\maxmin_{n;(I,t_1)}^{(J,t_2)}( r )$ \qq $\big\{ \!\!$ all `normalized' $Q$-adjusted polymer weights $(I,t_1) \to (J,t_2)$
are $\leq r$ $\!\!\big\}$  \fff{pwr}
$\lwr_{n;(I,0)}^{(J,1)}\big( \e , r  \big)$ \qq $\big\{ \!\!$ `normalized' polymer weight differences $(I,0) \to (J,1)$ are $\leq r$ $\!\!\big\}$ \fff{lwr} 
$\fbr_{n,k;(x,t_1)}^{(\bar{u},t_2)}(  r )$ \qq $\big\{ \!\!$   `normalized' $Q$-adjusted forward bouquet weight $(x,t_1) \to (\bar{u},t_2)$ is $\leq r$ $\!\!\big\}$   \fff{fbr}
$\bbr_{n,k;(\bar{v},t_1)}^{(y,t_2)}(  r )$ \qq $\big\{ \!\!$  the counterpart backward bouquet weight $(\bar{v},t_1) \to (y,t_2)$ is $\leq r$ $\!\!\big\}$ \fff{bbr}
$\fsc$ \qq the favourable surgical conditions event \fff{fsc}
$\fluc_{n;(x,t_1)}^{(y,t_2)}\big(a;K \big)$ \qq $\big\{ \!\!$  $\rho_{n;(x,t_1)}^{(y,t_2)}$ has `normalized' fluctuation $\in \! K \!$ at lifetime fractions $a$ and $1-a$ $\!\!\big\}$ \fff{fluc}
}\end{center}

%DisjtIndex. HorSepIndex. Forward proper. Backward proper.

\section{Multi-geodesic ordering}\label{s.mgo}

The 
%$\ell$-tuple 
multi-geodesic with given endpoints was defined in Section~\ref{s.staircase} and further discussed in Section~\ref{s.maxunique}, where its almost surely uniqueness was stated in Lemma~\ref{l.severalpolyunique}.
The principal aim of this appendix is to establish in Lemma~\ref{l.severalpolyorder} that this multi-geodesic satisfies a natural monotonicity property, 
moving to the right in response to rightward displacement of endpoints. The multi-polymer counterpart, Lemma~\ref{l.tworelations}(1), will follow immediately.

We introduce two ordering relations on sets of staircases that share their starting and ending heights.

Let $(i,j) \in \Z^2_<$ and $(x_1,y_1),(x_2,y_2) \in \R^2_<$ satisfy $x_1 \leq x_2$ and $y_1 \leq y_2$.
Consider  two staircases. $\phi_1 \in  \staircase_{i,j}(x_1,y_1)$
and $\phi_1 \in  \staircase_{i,j}(x_2,y_2)$.
We define two relations, $\prec$ and $\preceq$, that $\phi_1$ and $\phi_2$ may satisfy.
\begin{itemize}
\item We say that $\phi_1 \prec \phi_2$ if, whenever $(z,t) \in \R \times [i,j]$
is an element of $\phi_2$, the open planar line segment that emanates rightwards from $(z,t)$, namely $(z,\infty) \times \{ t \}$,
is disjoint from $\phi_1$.
\item We say that $\phi_1 \preceq \phi_2$ if, whenever $(z,t) \in \R \times [i,j]$
is an element of $\phi_1$, the closed planar line segment running rightwards from $(z,t)$, $[z,\infty) \times \{ t \}$,
intersects $\phi_2$.
\end{itemize}
In fact, each of these two conditions is equivalent to the conditions seen in Section~\ref{s.polyorder} when the statements there are made for zigzags, rather than staircases. 

Two basic properties are readily verified: 
\begin{itemize}
\item
$\phi_1 \prec \phi_2$ implies $\phi_1 \preceq \phi_2$;
\item  and 
$\phi_1 \preceq \phi_2$ and 
$\phi_2 \prec \phi_3$, or for that matter
$\phi_1 \prec \phi_2$ and 
$\phi_2 \preceq \phi_3$, imply that 
$\phi_1 \prec \phi_3$.
\end{itemize}

%These definitions may be extended to staircase $k$-tuples.
%Let $\bar{\phi} = \big(\phi_1,\cdots,\phi_k \big)$ and  $\bar{\phi}' = \big(\phi'_1,\cdots,\phi'_k \big)$ 
%be two $k$-vectors of staircases each of whose components have a common start and end height.
%We say that  $\bar{\phi} \preceq \bar{\phi}'$ precisely when $\phi_i \preceq \phi_i'$ for each $i \in \intint{k}$; and similarly for $\prec$.

When two staircases $\phi_1$
and $\phi_2$ share their pair  $(i,j) \in \Z^2_<$ of starting and ending heights, we may define their maximum $\phi_1 \vee \phi_2$.
Indeed, suppose that $\phi_1 \in  \staircase_{i,j}(x_1,y_1)$
and $\phi_2 \in  \staircase_{i,j}(x_2,y_2)$ where $(x_1,y_1),(x_2,y_2) \in \R^2_\leq$; contrary to a moment ago, we permit the cases that $x_2 < x_1$ or $y_2 < y_1$.
The maximum $\phi_1 \vee \phi_2$ is a staircase from $\big(  x_1 \vee x_2 , i \big)$ to $\big(  y_1 \vee y_2 , j \big)$.
It is equal to the unique staircase with these starting and ending points which is contained in the union of the horizontal and vertical segments of $\phi_1$ and $\phi_2$
and which dominates these two staircases in the $\preceq$ ordering. 
To these two staircases is also associated a minimum  $\phi_1 \wedge \phi_2$, a staircase that makes its way from $\big(  x_1 \wedge x_2 , i \big)$ to $\big(  y_1 \wedge y_2 , j \big)$.
It is analogously defined and is instead dominated by $\phi_1$ and $\phi_2$ according to $\preceq$.

\begin{lemma}\label{l.severalpolyorder}
For $k \in \N$, let $\bar{u},\bar{v},\bar{x},\bar{y} \in \R_\leq^k$ be four non-decreasing lists
such that $\bar{v} - \bar{u} \in [0,\infty)^k$ and $\bar{y} - \bar{x} \in [0,\infty)^k$. 
 Let $(j_1,j_2) \in \Z^2_<$. 
 Suppose that these parameters are such that the two maximizers
 $$
  \Big( P^k_{( \bar{u} , j_1 {\bf 1}) \to (\bar{x},j_2 {\bf 1}) ; i} : i \in \intint{k} \Big) 
   \, \, \, \textrm{ and } \, \, \, \Big(  P^k_{( \bar{v} ,  j_1 {\bf 1}) \to (\bar{y},j_2 {\bf 1}); i} : i \in \intint{k} \Big) 
  $$
  exist uniquely.
 Then
$$
  P^k_{( \bar{u} , j_1 {\bf 1}) \to (\bar{x},j_2 {\bf 1}) ; i} \preceq  P^k_{( \bar{v} ,  j_1 {\bf 1}) \to (\bar{y},j_2 {\bf 1}); i }  \, \, \textrm{for $i \in \intint{k}$}  \, .
$$
\end{lemma}
(In this appendix, we write ${\bf 1}$, rather than $\bar{\bf 1}$, for a $k$-tuple whose components equal one.)

This result permits us to derive Lemma~\ref{l.tworelations}.

\noindent{\bf Proof of Lemma~\ref{l.tworelations}: (1).}
This assertion is the multi-polymer counterpart to  Lemma~\ref{l.severalpolyorder}
and follows immediately from this lemma.

\noindent{\bf (2).} This result is a direct consequence of the definitions of the relations $\prec$ and $\preceq$. \qed

\noindent{\bf Proof of Lemma~\ref{l.severalpolyorder}.} For convenience, denote the two maximizers by $\big( \rhomac_{1,i}: i \in \intint{k} \big)$ and $\big( \rhomac_{2,i}: i \in \intint{k} \big)$ in such a way that we seek to show that $\rhomac_{1,i} \preceq \rhomac_{2,i}$ for $i \in \intint{k}$.

Suppose that this ordering property fails. 
Our plan is to argue that this assumption leads to a violation of the almost sure uniqueness of multi-geodesics with given endpoints asserted by Lemma~\ref{l.severalpolyunique}.

Let $i  \in \intint{k}$ be maximal such that $\rhomac_{1,i} \npreceq \rhomac_{2,i}$. Let $j \in \intint{i}$ be minimal such that $\rhomac_{1,j} \npreceq \rhomac_{2,i}$. 
(The stroke through the symbol, $\preceq$, and later $\prec$, is used to indicate that the relation does not hold.)

We now set the value of a parameter $L$ in the interval $\llbracket 0, j-1 \rrbracket$.
If $\rhomac_{1,j-1} \prec \rhomac_{2,i}$, then $L$ is set equal to zero. In the other case, we set
$L \in \intint{j-1}$ so that 
$$
\rhomac_{1,j-m-1} \nprec \rhomac_{2,i-m} \, \, \, \textrm{for each $m \in \llbracket 0, L - 1 \rrbracket$,}
$$
and, if $L \not= j-1$, then $\rhomac_{1,j-L-1} \prec \rhomac_{2,i-L}$.

Whatever the value of $L$, we now specify two $k$-vectors of staircases, $\big( \rhomac'_{1,\ell} : \ell \in \intint{k} \big)$ and  $\big( \rhomac'_{2,\ell} : \ell \in \intint{k} \big)$, according to the formulas
$$
\rhomac'_{1,\ell} = 
\begin{cases} 
 \, \rhomac_{1,\ell} &  \ell \in \intint{k} \setminus \llbracket j - L, j \rrbracket \, , \\
 \, \rhomac_{1,j - m} \wedge \rhomac_{2,i-m} &  \ell = j - m \, \, \textrm{with $m \in \llbracket 0,L \rrbracket$} \, ,
\end{cases}
$$ 
and
$$
\rhomac'_{2,\ell} = 
\begin{cases} 
 \, \rhomac_{2,\ell} &  \ell \in \intint{k} \setminus \llbracket i - L, i \rrbracket \, , \\
 \, \rhomac_{1,j - m} \vee \rhomac_{2,i-m} &  \ell = i - m \, \, \textrm{with $m \in \llbracket 0,L \rrbracket$} \, .
\end{cases}
$$ 
Using such notation as $P'_{1,\cdot} = \big( P'_{1,i}: i \in \intint{k} \big)$,
we now make a claim with four parts.
 
 \noindent{\em Claim.}
\begin{enumerate}
\item
$\rhomac'_{1,\cdot}$ does not equal $\rhomac_{1,\cdot}$, and nor does $\rhomac'_{2,\cdot}$  equal $\rhomac_{2,\cdot}$. 
\item The vector $\rhomac'_{1,\cdot}$
is an element of
$D^k_{(x {\bf 1},j_1 {\bf 1}) \to (\bar{u},j_2{\bf 1})}$. 
\item The vector  $\rhomac'_{2,\cdot}$ belongs to  
$D^k_{(x {\bf 1},j_1 {\bf 1}) \to (\bar{v},j_2{\bf 1})}$.
\item These two vectors are maximizers of their respective sets. 
\end{enumerate}
We verify the four assertions in turn.

\begin{figure}
\begin{center}
\includegraphics[width=0.9\textwidth]{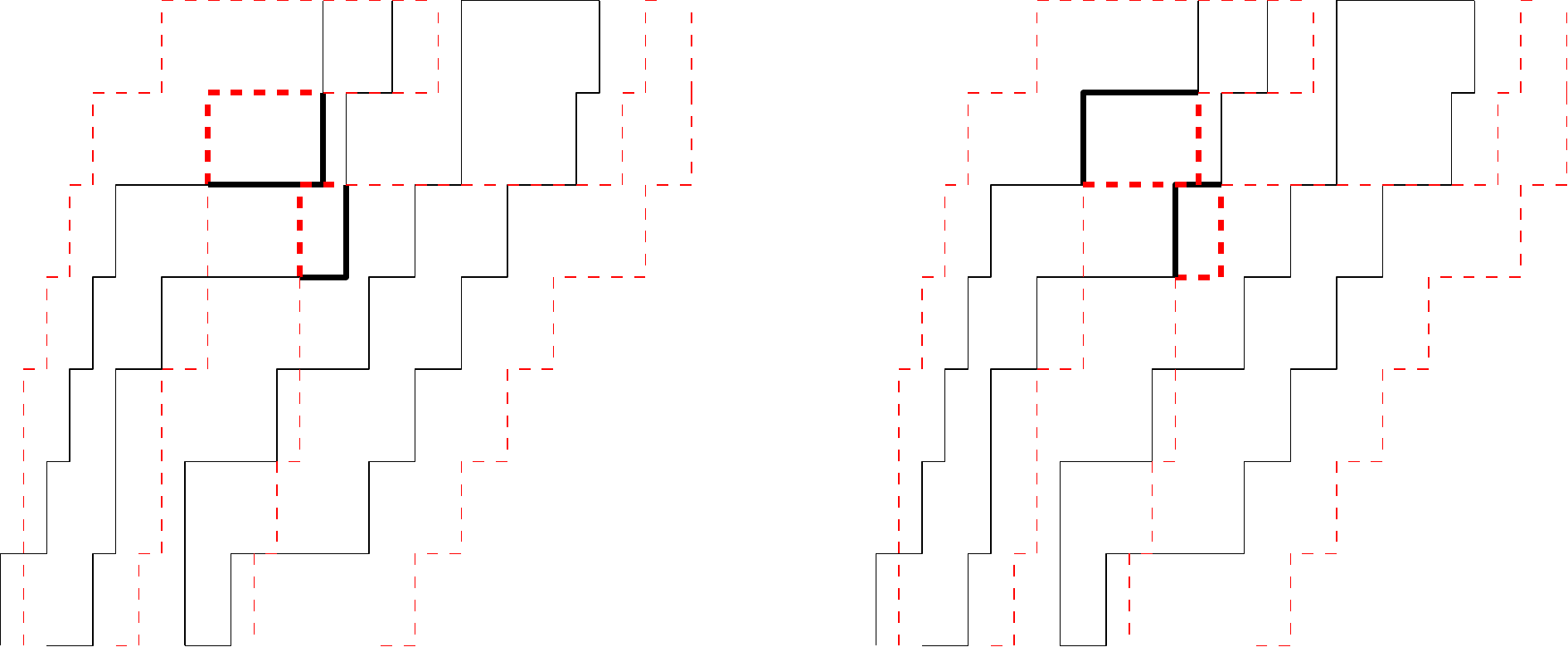}
\end{center}
\caption{Illustrating the proof of Lemma~\ref{l.severalpolyorder}. 
In this example, $\ell = 4$. On the left, the consecutive staircases $\big\{ \rhomac_{1,i}: i \in \intint{4} \big\}$ are indicated in unbroken black.
The four components of $\rhomac_{2,\cdot}$  are shown in dashed red. 
On the right, the $\rhomac'$ counterparts are similarly indicated.
We have $i=3$, $j=2$ and $L =1$.  
Substaircases that change hands are highlighted.
There are two such rearrangements. One leads to $P'_{1,1} = P_{1,1} \wedge P_{2,2}$ and $P'_{2,2} = P_{1,1} \vee P_{2,2}$, and the other to
$P'_{1,2} = P_{1,2} \wedge P_{2,3}$ and $P'_{2,3} = P_{1,2} \vee P_{2,3}$.
Otherwise,  $P'$ components coincide with their $P$ counterparts.  
}\label{f.switch}
\end{figure}

\noindent{\em Verifying Claim (1):} Since $\rhomac'_{1,j} = \rhomac_{1,j} \wedge \rhomac_{2,i}$ and $\rhomac_{1,j} \not\preceq \rhomac_{2,i}$, we see that, although $\rhomac'_{1,j} \preceq \rhomac_{1,j}$, we also have 
 $\rhomac'_{1,j} \not= \rhomac_{1,j}$.
 Similarly, since $\rhomac'_{2,i} = \rhomac_{1,j} \vee \rhomac_{2,i}$
 while $\rhomac_{1,j} \not\preceq \rhomac_{2,i}$, we see that, although $\rhomac'_{2,i} \preceq P_{2,i}$, we also have $\rhomac'_{2,i} \not= \rhomac_{2,i}$.

\noindent{\em Verifying Claim (2):}
For each $i \in \intint{k}$, the staircase $\rhomac_{2,i}$ staircase begins and ends
at or to the right of $\rhomac_{1,i}$. For this reason, each staircase $\rhomac'_{1,i}$ shares its starting and ending points with $\rhomac_{1,i}$. 
In light of this, it is enough in order to show that $\rhomac'_{1,\cdot}\in
 D^k_{(x {\bf 1},j_1 {\bf 1}) \to (\bar{u},j_2{\bf 1})}$ that it be demonstrated that 
$\rhomac'_{1,\ell} \prec \rhomac'_{1,\ell+1}$ for each $\ell \in \llbracket 0,k-1 \rrbracket$. This statement follows from the ordering under $\prec$ of the components of $\rhomac_{1,\cdot}$ in the case that $\ell \not\in \llbracket j-L-1,j \rrbracket$.

We have that $\rhomac'_{1,j} \prec \rhomac'_{1,j+1}$ since $\rhomac_{1,j} \prec \rhomac_{1,j+1}$; and, for $m \in \intint{L}$, that  $\rhomac'_{1,j - m} \prec \rhomac'_{1,j-m+1}$ since $\rhomac_{1,j-m} \prec \rhomac_{1,j-m+1}$ and $\rhomac_{2,i-m} \prec \rhomac_{2,i-m+1}$. Finally, when $L < j-1$, we have   that  $\rhomac'_{1,j - L - 1} \prec \rhomac'_{1,j-L}$ since $\rhomac_{1,j-L-1} \prec \rhomac_{1,j-L}$ and $\rhomac_{1,i-L-1} \prec \rhomac_{2,i-L}$ (an inequality which is indeed valid when $L \not= j-1$). Note that when $L = j-1$, $j - L - 1$ equals zero, so that the verification in the preceding sentence is not made in this case.   

Thus, we verify that  $\rhomac'_{1,\cdot}\in 
 D^k_{(x {\bf 1},j_1 {\bf 1}) \to (\bar{u},j_2{\bf 1})}$.

\noindent{\em Verifying Claim (3):}
For each $i \in \intint{k}$, 
 $\rhomac'_{2,i}$ shares its starting and ending points with $\rhomac_{2,i}$, for the same reason that 
 $\rhomac'_{1,i}$ and $\rhomac_{1,i}$ do.

It is thus enough to verify that $\rhomac'_{2,\cdot}\in D^k_{(x {\bf 1},j_1 {\bf 1}) \to (\bar{v},j_2{\bf 1})}$. This is confirmed in view of:
\begin{itemize}
\item $\rhomac'_{2,\ell} \prec \rhomac'_{2,\ell+1}$ for $\ell \in \intint{k} \setminus \llbracket i-L , i \rrbracket$ since $\rhomac_{2,\ell} \prec \rhomac_{2,\ell+1}$;
\item $\rhomac'_{2,i} \prec \rhomac'_{2,i+1}$  since $\rhomac_{2,i} \prec \rhomac_{2,i+1}$ and, by $j < i+1$ and then $i$'s maximality, $\rhomac_{1,j} \prec \rhomac_{1,i+1} \preceq \rhomac_{2,i+1}$;
\item $\rhomac'_{2,i - m - 1} \prec \rhomac'_{2,i - m}$ for $m \in \intint{L-1}$ since $\rhomac_{1,j - m -1} \prec \rhomac_{1,j -m}$ and $\rhomac_{2,i - m -1} \prec \rhomac_{2,i -m}$;
\item and $\rhomac'_{2,i - L - 1} \prec \rhomac'_{2,i - L}$ since  $\rhomac_{2,i - L - 1} \prec \rhomac_{2,i - L}$.  
\end{itemize}

\noindent{\em Verifying Claim (4):}
To argue that the two $k$-vectors $\rhomac'_{1,\cdot}$ and $\rhomac'_{2,\cdot}$ are maximizers of their respective sets, we first note that
$$
W \big( \rhomac'_1 \big) + 
W \big( \rhomac'_2 \big) = 
W \big( \rhomac_1 \big) + 
W \big( \rhomac_2 \big)  \, . 
$$
Since $\rhomac'_1 \in  D^k_{(x {\bf 1},j_1 {\bf 1}) \to (\bar{u},j_2{\bf 1})}$
and $\rhomac_1$ is a maximizer of this set, $W\big( \rhomac'_1  \big) \leq W\big( \rhomac_1  \big)$. Likewise, $W \big( \rhomac'_2  \big) \leq W \big( \rhomac_2  \big)$.  This circumstance forces $W \big( \rhomac'_1 \big)$ to equal 
$W \big( \rhomac_1 \big)$ and 
$W \big( \rhomac'_2 \big)$ to equal 
$W \big( \rhomac_2 \big)$. Thus, $\rhomac'_1$ and $\rhomac'_2$ is each a maximizer of the set to which it belongs. 

The proof of the claim complete, we need now merely invoke the
 maximizer uniqueness Lemma~\ref{l.severalpolyunique}
to learn  that $\rhomac'_1 = \rhomac_1$, and also that
$\rhomac'_2 = \rhomac_2$. However, we have seen that neither of these equalities holds. From this contradiction, we learn that  $\rhomac_{1,i} \preceq \rhomac_{2,i}$ for $i \in \intint{k}$, as we sought to do. \qed

\section{Normalized ensembles are regular: deriving Proposition~\ref{p.scaledreg}}\label{s.normal}

In order to prove Proposition~\ref{p.scaledreg},
it is convenient to expand a little on the review in Section~\ref{s.encode}  of the definition of the scaled forward ensemble $\mc{L}_{n;(x,t_1)}^{\uparrow;t_2}$ and its normalized counterpart. It is natural to view the first object as the scaled counterpart to an unscaled forward ensemble which we now specify.

To do so, let $(m_1,m_2) \in \N^2_\leq$
and $u \in \R$.
Define the unscaled forward ensemble
$$
L_{(u,m_1)}^{\uparrow;m_2}: \intint{m_2 - m_1 + 1} \times [u,\infty) \to \R
$$
 with base-point $(u,m_1)$ and end height $m_2$
by insisting that, for each $k \in \intint{m_2 - m_1 + 1}$
and $v \geq u$,
$$
 \sum_{i=1}^k L_{(u,m_1)}^{\uparrow;m_2}(i,v) = M^k_{(u,m_1) \to (v,m_2)} \, .
$$

The relation of this object to its scaled counterpart may be described by
letting $(n,t_1,t_2) \in \N \times \R^2_\leq$ be a compatible triple,  and $x \in \R$. The scaled  forward ensemble
$$
\mc{L}_{n;(x,t_1)}^{\uparrow;t_2}: \intint{n \tot + 1} \times \big[ x - 2^{-1} n^{1/3} \tot  , \infty\big) \to \R
$$
is seen to 
satisfy, for each $k \in  \intint{n \tot + 1}$
and $y \geq x - 2^{-1} n^{1/3} \tot$, the identity
\begin{equation}\label{e.ln}
 \mc{L}_{n;(x,t_1)}^{\uparrow;t_2}(k,y) = 2^{-1/2} n^{-1/3} \Big( L_{(n t_1 + 2n^{2/3}x,n t_1)}^{\uparrow;n t_2}\big( k,  n t_2 + 2n^{2/3} y \big) - 2 n t_{1,2}  - 2n^{2/3}(y-x) \Big) \, .
\end{equation}
%Alternatively, we may note that 
%$$
% \sum_{i=1}^k \mc{L}_{n;(x,t_1)}^{\uparrow;t_2}(i,y) = 2^{-1/2} n^{-1/3} \bigg( M^k_{(n t_1 + 2n^{2/3}x,n t_1) \to (n t_2 + 2n^{2/3}y,n t_2)} - k \Big( 2 n t_{1,2}   + 2n^{2/3}(y-x) \Big) \bigg) \, .
%$$
A further, parabolic, change of coordinates~(\ref{e.scaledln})
then specifies the normalized forward ensemble 
 $\scaledle_{n;(x,t_1)}^{\uparrow;t_2}$ in terms of the scaled one~$\mc{L}_{n;(x,t_1)}^{\uparrow;t_2}$.

A similar story holds for the backward ensemble. Indeed,
considering now $(m_1,m_2) \in \N^2_\leq$ 
and $v \in \R$, we may
define the unscaled backward ensemble
$$
L_{m_1}^{\downarrow;(v,m_2)}: \intint{m_2 - m_1 + 1} \times (-\infty, v ] \to \R
$$
 with base-point $(v,m_2)$ and end height $m_1$
by insisting that, for each $k \in \intint{m_2 - m_1 + 1}$
and $u \leq v$,
$$
 \sum_{i=1}^k L_{m_1}^{\downarrow;(v,m_2)}(i,u) = M^k_{(u,m_1) \to (v,m_2)} \, .
$$

Again fixing a  compatible triple $(n,t_1,t_2) \in \N \times \R^2_\leq$, and now letting $y \in \R$,
we note that the scaled backward ensemble
$$
\mc{L}_{n;t_1}^{\downarrow;(y,t_2)}: \intint{n \tot + 1} \times \big( - \infty ,  y + 2^{-1} n^{1/3} \tot \big] \to \R
$$
satisfies, for each $k \in  \intint{n \tot + 1}$
and $x \leq y + 2^{-1} n^{1/3} \tot$,
$$
\mc{L}_{n;t_1}^{\downarrow;(y,t_2)}(k,x) = 2^{-1/2} n^{-1/3} \Big( L_{n t_1}^{\downarrow;( n t_2 + 2n^{2/3}y , n t_2)}\big( k,  n t_1 + 2n^{2/3} x \big) - 2 n \tot  - 2n^{2/3}(y-x) \Big) \, .
$$
%Equally,
%$$
% \sum_{i=1}^k 
%\mc{L}_{n;t_1}^{\downarrow;(y,t_2)}(i,x) 
%   = 2^{-1/2} n^{-1/3} \bigg( M^k_{(n t_1 + 2n^{2/3}x,n t_1) \to (n t_2 + 2n^{2/3}y,n t_2)} - k \Big( 2 n \tot  + 2n^{2/3}(y-x) \Big) \bigg) \, .
%$$

Now consider again a compatible triple $(n,t_1,t_2)$ as well as $x \in \R$, and write $N = n \tot + 1 \in \N$.  Consider the ensemble $L_N: \intint{N} \times [   0 , \infty ) \to \R$, given by
$$
 L_N(i,u) =  L_{(n t_1 + 2n^{2/3}x, n t_1)}^{\uparrow; n t_2}\big( i , 
 n t_1 + 2n^{2/3}x + u   \big) \qquad \textrm{for $(i,u) \in  \intint{N} \times [   0 , \infty )$} \, .
$$
By~\cite[Theorem~$7$]{O'ConnellYor}, this ensemble has the law of an $N$-curve  Dyson's Brownian motion each of whose curves begins at the origin at time zero. 
 The ensemble has a scaled counterpart 
 $\mathcal{L}_N^\scal: \intint{N} \times [-2^{-1} N^{1/3},\infty )  \to \R$, given by   
$$
 \mathcal{L}_N^\scal(i,u) =  2^{-1/2} N^{-1/3} \Big(   L_N\big( i,N + 2N^{2/3} u \big) - 2N - 2 N^{2/3} u \Big) \, .
$$
By~\cite[Proposition~$2.5$]{BrownianReg}, 
%p.lereg
we have the 

\noindent{\em Key fact:} there exist constants $\rsC,\rsc \in (0,\infty)$ such that, for all choices of the parameters $(n,t_1,t_2,x)$, the ensemble 
 $\mathcal{L}_N^\scal$ is $(c,C)$-regular.

There is only one, very mundane, difference between the ensemble $\mathcal{L}_N^\scal$ and the normalized forward ensemble $\scaledle_{n;(x,t_1)}^{\uparrow;t_2}$ that is the subject of  Proposition~\ref{p.scaledreg}.
With our choice of $N = n\tot+ 1$, the latter ensemble is simply equal to  $\mathcal{L}_{N-1}^\scal$. The key fact thus bears the essence of the proof of  Proposition~\ref{p.scaledreg}. The formal proof, which will now be given, is merely a matter of verifying that the adjustment of the parameter value $N \to N -1$
is inconsequential for the purpose of checking the conditions enjoyed by a $(c,C)$-regular ensemble.  

\noindent{\bf Proof of Proposition~\ref{p.scaledreg}.} 
Note from~(\ref{e.ln}) and~(\ref{e.scaledln}) that  $\scaledle_{n;(x,t_1)}^{\uparrow;t_2}(k,z)$ equals
$$
   2^{-1/2} (n \tot)^{-1/3} \Big( L_{(n t_1 + 2n^{2/3}x, n t_1 )}^{\uparrow; n t_2}\big( k,  n t_2 + 2n^{2/3} x + 2(n \tot)^{2/3} z  \big) - 2 n t_{1,2}  - 2  (n \tot)^{2/3}z  \Big) \, .
$$
Note that 
\begin{eqnarray*}
 & & L_{(n t_1 + 2n^{2/3}x, n t_1)}^{\uparrow; n t_2}\big( k,  n t_2 + 2n^{2/3} x + 2(n \tot)^{2/3} z  \big) - 2 n t_{1,2}  - 2  (n \tot)^{2/3}z  \\
  & =  & L_N \big( k  , n \tot   +  2(n \tot)^{2/3} z  \big) -  2 n t_{1,2}  - 2  (n \tot)^{2/3}z \\
  & =  & L_N \big( k  , N    +  2 N^{2/3} z (n \tot N^{-1})^{2/3} + n \tot - N  \big) \\
   & & \qquad  -  2N  - 2 N^{2/3}z -  2 (n t_{1,2} - N)  - 2 \big(   (n \tot)^{2/3} - N^{2/3}\big) z \\
%     & =  & L_N \big( k  , N    +  2 N^{2/3} \tilde{z}  \big)  -  2N  - 2 N^{2/3} \tilde{z} +  2   - 2 \big(   (n \tot)^{2/3} - N^{2/3}\big) z  - 2 N^{2/3}(z - \tilde{z}) \\
    & =  & L_N \big( k  , N    +  2 N^{2/3} \tilde{z}  \big) -  2N  - 2 N^{2/3} \tilde{z} + \phi_1 \, , 
 \end{eqnarray*}
 where $\tilde{z} =    z (n \tot N^{-1})^{2/3} - 2^{-1} N^{-2/3}$ satisfies $\tilde{z} = z \big( 1 + O(N^{-1}) \big) + O(N^{-2/3})$.
Recalling that $N = n \tot + 1$,  the error~$\phi_1$ is seen to equal  $2   - 2 \big(   (N - 1)^{2/3} - N^{2/3}\big) z  - 2 N^{2/3}(z - \tilde{z})$ and thus to satisfy
 $$
 \phi_1 = O(1) + \vert z \vert O(N^{-1/3}) + \vert z \vert O(N^{-1/3}) + O(1) =  \vert z \vert O(N^{-1/3}) + O(1) \, .
 $$

We see then that $\scaledle_{n;(x,t_1)}^{\uparrow;t_2}(k,z)$  equals
$$
   2^{-1/2} (n \tot)^{-1/3}  \Big( L_N \big( k  , N    +  2 N^{2/3} \tilde{z}  \big) -  2N  - 2 N^{2/3} \tilde{z} + \phi_1 \Big) \, =  \,     \mathcal{L}_N^\scal \big( k, \tilde{z} \big)   \phi_2  \,     - \,    2^{-1/2} (n \tot)^{-1/3}  \phi_1 \, ,
$$
 where $\phi_2 =   (n \tot)^{-1/3} N^{1/3} = 1 +  O(N^{-1})$.
We find then that
$$
\scaledle_{n;(x,t_1)}^{\uparrow;t_2}(k,z) =  \mathcal{L}_N^\scal \big( k, \tilde{z} \big)  \big(  1 +  O(N^{-1}) \big) +    \vert z \vert O(N^{-2/3}) + O(N^{-1/3}) \, ,
$$
where $\tilde{z}$ differs from $z$ by $\vert z \vert O(N^{-1}) + O(N^{-2/3})$.
We may thus note that
\begin{eqnarray*}
 & & \scaledle_{n;(x,t_1)}^{\uparrow;t_2}(k,z) + 2^{-1/2} z^2 \\
 &  = &   \mathcal{L}_N^\scal \big( k, \tilde{z} \big)   \big(  1 +  O(N^{-1}) \big)   + 2^{-1/2} \tilde{z}^2  +  \vert z \vert^2 O(N^{-1}) +  \vert z \vert O(N^{-2/3})  +    \vert z \vert O(N^{-2/3}) + O(N^{-1/3}) \\
 &  = &   \mathcal{L}_N^\scal \big( k, \tilde{z} \big)   \big(  1 +  O(N^{-1}) \big)   + 2^{-1/2} \tilde{z}^2  +  \vert z \vert^2 O(N^{-1}) +  \vert z \vert O(N^{-2/3})  + O(N^{-1/3}) \, . 
\end{eqnarray*}
When $z = O(N^{1/2})$, the triple sum of error terms in the last expression is $O(1)$. The three $\rmreg$ properties hold for $\mathcal{L}_N^\scal$
by the key fact stated before the proof. We want to learn that these hold also for 
$\scaledle_{n;(x,t_1)}^{\uparrow;t_2}$. The bound on the triple sum error permits the passage of $\rmreg(2)$ and $\rmreg(3)$ from   $\mathcal{L}_N^\scal$  to
$\scaledle_{n;(x,t_1)}^{\uparrow;t_2}$, provided that $\rsC < \infty$ is increased, and $\rsc > 0$ decreased, in order that the case of small $N$ is included. The property $\rmreg(1)$ also makes the passage, because the left endpoint modulus $\xnmac$ differs by a factor of $\big( n \tot N^{-1} \big)^{1/3} = 1 +  O(N^{-1})$ between the two ensembles; again, $\rsc > 0$ is being decreased in order to capture the case of small $N$. Thus, the ensemble $\scaledle_{n;(x,t_1)}^{\uparrow;t_2}$ is $(\rsc,\rsC)$-regular, as we sought to show. \qed

%arxiv/submit

\bibliographystyle{plain}

\bibliography{airy}

\end{document}